\documentclass[a4paper,12pt]{amsart}
\usepackage[frenchb,english]{babel}
\usepackage[latin1]{inputenc}
\usepackage{pdfsync}
\usepackage[T1]{fontenc}
\usepackage{amsmath}
\usepackage{amssymb}
\usepackage{amsxtra}
\usepackage{amscd}
\usepackage{amsthm}
\usepackage{amsfonts}
\usepackage{eucal}
\usepackage[all]{xy}
\usepackage{graphicx}
\usepackage{comment}
\usepackage{epsfig}
\usepackage{psfrag}
\usepackage{mathrsfs}
\usepackage{amscd}
\usepackage{rotating}
\usepackage{lscape}
\usepackage{amsbsy}
\usepackage{verbatim}
\usepackage{moreverb}
\usepackage{url}

\newcommand{\nc}{\newcommand}
\nc{\renc}{\renewcommand}

\nc\restr[2]{{ 
  \left.\kern-\nulldelimiterspace    #1  
  \vphantom{\big|}  
  \right|_{#2}  
  }}

\newtheorem{thm}{Theorem}[section]
\newtheorem{prop}[thm]{Proposition}
\newtheorem{lem}[thm]{Lemma}
\newtheorem{sous-lem}[thm]{Sub-lemma}
\newtheorem{cor}[thm]{Corollary}

\newtheorem{def-prop}[thm]{Definition-Proposition}
\theoremstyle{definition}
\newtheorem{defi}[thm]{Definition}
\newtheorem{rem}[thm]{Remark}

\newtheorem{observation}[thm]{Observation}
\numberwithin{equation}{section}

\renc{\sec}{\section}
\nc{\ssec}{\subsection}
\nc{\sssec}{\subsubsection}

\nc{\thmref}[1]{theorem~\ref{#1}}
\nc{\secref}[1]{section~\ref{#1}}
\nc{\lemref}[1]{lemma~\ref{#1}}
\nc{\defiref}[1]{definition~\ref{#1}}
\nc{\propref}[1]{proposition~\ref{#1}}
\nc{\corref}[1]{corollary~\ref{#1}}
\nc{\constructionref}[1]{construction~\ref{#1}}
\nc{\conjref}[1]{conjecture~\ref{#1}}
\nc{\remref}[1]{remark~\ref{#1}}
\nc{\questref}[1]{question~\ref{#1}}
\nc\Omegasour{\hbox{$\buildrel\smile\over{\vrule height 6pt depth 0pt width 0pt \smash \Omega}$}}

\nc{\on}{\operatorname}

\nc\wt{\widetilde}
\nc\wh{\widehat}
\nc\ol{\ov}
\nc{\oc}[1]{{\overset{\circ}{#1}}}
\nc{\ov}[1]{{\overline{#1}}}
\nc{\isor}[1]{{\xrightarrow[\raisebox{0.25 em}{\smash{\ensuremath{\sim}}}]{#1}}}
\nc{\modmod}{/ \! \! /}

\nc{\mc}{\mathcal}
\nc{\mf}{\mathfrak}
\nc{\mr}{\mathrm}
\nc{\mb}{\mathbb}
\nc{\mbf}{\mathbf}

\nc{\R}{{\mathbb R}}
\nc{\Z}{{\mathbb Z}}
\nc{\N}{{\mathbb N}}
\nc{\C}{{\mathbb C}}
\nc{\Q}{{\mathbb Q}}

\nc{\Fq}{{\mathbb F}_q}
\nc{\Fl}{{\mathbb F}_\ell}
\nc{\Fqbar}{\ol{{\mathbb F}_q}}
\nc{\Flbar}{\ol{{\mathbb F}_\ell}}
\nc{\Zl}{{\mathbb Z}_\ell}
\nc{\Zlbar}{\ol{{\mathbb Z}_\ell}}
\nc{\Ql}{{\mathbb Q}_\ell}
\nc{\Qlbar}{\ol{{\mathbb Q}_\ell}}
\nc{\hl}{\overset{\leftarrow}h{}}
\nc{\hr}{\overset{\rightarrow}h{}}
\nc{\Gr}{{\on{Gr}}}
\nc{\Hecke}{\on{Hecke}}
 \nc{\Hom}{\on{Hom}}
 \nc{\Coker}{\on{Coker}}
 \nc{\Ker}{\on{Ker}}
 \nc{\Lie}{\on{Lie}}
\nc{\Loc}{\on{Loc}}
\nc{\Pic}{\on{Pic}}
\nc{\Bun}{\on{Bun}}
\nc{\IC}{\on{IC}}
\nc{\Aut}{\on{Aut}}
\nc{\Perv}{\on{Perv}}
\nc{\pos}{{\on{pos}}}
\nc{\Sym}{\on{Sym}}

\nc{\ta} {{}^\tau}
\nc {\tu}[1]{{}^{\tau^{#1}}\!}

\nc{\Id}{\on{Id}}
\nc{\Fil}{\on{Fil}}
\nc{\pr}{\on{pr}}
\nc{\Res}{\on{Res}}
\nc{\cusp}{\on{cusp}}
\nc{\Frob}{\on{Frob}}
\nc{\diag}{\Delta}
\nc{\gr}{\on{gr}}
\nc{\Inj}{\on{Inj}}
\nc{\Bl}{\on{Bl}}
\nc{\dem}{\noindent {\bf Proof. }}
\nc{\cqfd}{{\ }\hfill $\square$ \vskip 1mm}
\nc{\s}[1]{\langle #1 \rangle}
\nc{\Cht}{\on{Cht}}
\nc{\isom}{\overset {\thicksim}{\to}}
\nc{\sm}{\smallsetminus}

\begin{document}

\title[Chtoucas and Langlands parameterization]{
Introduction to  chtoucas for   reductive groups and to the  global  Langlands parameterization}

\author{Vincent Lafforgue}
\thanks{The author is supported by ANR-13-BS01-0001-01}
\address{Vincent Lafforgue: CNRS and Institut Fourier, UMR 5582, Université Grenoble Alpes, 
 100 rue des Maths, 38610 Gières, France.}

\date{\today}
\maketitle

      \section*{Introduction}
  This is a   translation in English\footnote{This is an experiment, see \url{http://vlafforg.perso.math.cnrs.fr/bilingual.pdf} for general ideas about bilingual versions.} of     version 5 of \cite{intro}.  
  
  We explain  the  result of 
      \cite{coh} and 
    give   {\it all the  ideas}
     of the proof. This article corresponds essentially to the introduction of \cite{coh}  (slightly expanded) and to subsections 12.1 and 12.2 (shortened).       
      
         We show  the direction ``automorphic to Galois'' of the global Langlands correspondence    \cite{langlands67} for all  reductive  group  $G$ over a function field.              
      Moreover  we construct a {\it canonical} decomposition  of the space  of  cuspidal automorphic forms, indexed by  global Langlands parameters.
  We do not obtain any new result when  $G=GL_{r}$ since everything was already known by  Drinfeld \cite{drinfeld78,Dr1,drinfeld-proof-peterson,drinfeld-compact} for $r=2$  and Laurent Lafforgue \cite{laurent-inventiones} for $r$ arbitrary. 

          The  method  is completely  independent of Arthur-Selberg trace formulas. It uses the following ingredients 
        \begin{itemize}
     \item classifying stacks of chtoucas, introduced by   Drinfeld for $GL_{r}$  
     \cite{drinfeld78,Dr1} and generalized  to  all reductive groups  by Varshavsky
     \cite{var} 
   \item the geometric Satake equivalence of   Lusztig, Drinfeld, Ginzburg, and  Mirkovic--Vilonen 
   \cite{lusztig-satake,ginzburg,hitchin,mv}. 
     \end{itemize}

\noindent {\bf Acknowledgements.  }
This work  comes from a work in progress with  Jean-Benoît Bost, 
and would not  exist either without many discussions with  
Alain Genestier. I thank them also very much for their help.  

This article relies on works of  Vladimir Drinfeld,  both on  chtoucas and  on the  geometric Langlands program, and I acknowledge my huge debt to him. 
 
I thank Vladimir Drinfeld, Dennis Gaitsgory,  Jochen Heinloth, Laurent Lafforgue, Sergey Lysenko and Yakov Varshavsky for   their explanations. I thank  
  Gebhard B\"ockle,  Michael Harris, Chandrashekhar Khare and Jack Thorne who found the statement of  \lemref{prop-harris} (which improves the redaction and is also used in their  article \cite{boeckle-harris...}). 
I thank also 
    Ga\"etan Chenevier, Pierre Deligne,  Ed Frenkel,  Nicholas Katz,  Erez Lapid, Gérard Laumon,  Colette Moeglin, Sophie Morel, Ngô Bao Châu, Ngo Dac Tuan, Jean-Pierre Serre,  Jean-Loup Waldspurger   and  Cong Xue.

 I am very grateful to CNRS.        Langlands program is far from my    first subject and I would never  have been able to devote myself to it without the great freedom 
      given  to researchers for their works. 
 
\noindent {\bf Contents. } In section \ref{para-enonces} we state the main result (for  $G$ split), the idea of the proof and two intermediate statements which give the structure of the proof. Section \ref{para-commentaires}, which is a summary of  subsections 12.1 and 12.2   of \cite{coh}, discusses the  non split case,  some complements and a number of open problems. Section \ref{defi-chtou-intro} introduces the stacks of chtoucas. Section \ref{section-esquisse-abc} shows that the ``Hecke-finite'' part of their  cohomology satisfies some properties which were stated in  \propref{prop-a-b-c}  of section \ref{para-enonces}.  With the help of these properties we show in  section \ref{intro-idee-heurist} the main theorem,  which had been stated in  section \ref{para-enonces}. In sections \ref{subsection-crea-annihil-intro} to \ref{subsection-intro-decomp} we show some properties of the  cohomology of the stacks  of chtoucas,  which had been admitted in  section \ref{section-esquisse-abc} and at the end of section \ref{intro-idee-heurist}.   
Section \ref{subsection-link-langl-geom}
explains the link with the geometric Langlands program. 
Lastly  section \ref{intro-previous-works}
discusses the relation with previous works.   

Depending on the time at his disposal, the reader may restrict himself to section \ref{para-enonces}, to sections \ref{para-enonces} and \ref{para-commentaires}, or to  sections \ref{para-enonces} to  \ref{intro-idee-heurist}. 
Sections \ref{para-enonces}, \ref{para-commentaires} and \ref{intro-idee-heurist} do not require any background in algebraic geometry. 

\section{Statement and main  ideas}\label{para-enonces}

   \subsection{Preliminaries}         Let $\Fq$ be a finite field. 
      Let $X$ be a smooth projective geometrically irreducible curve   over $\Fq$ and  let $F$ be its  field of functions. 
       Let $G$ be a connected   reductive  group    over $F$. 
    Let $\ell$ be a prime number not dividing  $q$.

    To state  the main theorem we assume that $G$ is split (the non split case will be explained in  section \ref{cas-non-deploye}).  
        We denote by   $\wh G$ the Langlands dual group  of  $G$,  considered as  a split group over  $\Ql$. 
        Its roots and   weights are the coroots and  coweights of $G$, and reciprocally (see \cite{borel-corvallis}  for more details).

For  topological  reasons we have to work with  finite extensions of $\Ql$ instead of  $\Qlbar$.  Let $E$ be a finite extension of $\Ql$ containing a square root of $q$ and let $\mc O_{E}$ be its ring of integers.

         Let $v$ be  a place of $X$. We denote by   $\mc O_{v}$ the completed local ring at  $v$ and  by $F_{v}$ its field  of fractions. 
   We have  the  Satake isomorphism    $[V]\mapsto h_{V,v}$ from the ring of  representations of  $\wh G$ (with coefficients in $E$) to  the Hecke algebra  
      $C_{c}(G(\mc O_{v})\backslash G(F_{v})/G(\mc O_{v}),E)$ (see  
      \cite{satake,cartier-satake,gross}). In fact the $h_{V,v}$  for $V$ irreducible form a basis over  $\mc O_{E}$ of $C_{c}(G(\mc O_{v})\backslash G(F_{v})/G(\mc O_{v}),\mc O_{E})$. 
      We denote by    $\mb A=\prod_{v\in |X|} ' F_{v}$ the ring of adèles  of $F$ and we write  $\mb O=
      \prod_{v\in |X|} \mc O_{v}$. 
              Let  $N$ be a finite subscheme of  $X$.   
   We denote by   
   $\mc O_{N}$ the ring of functions of  $N$, and 
by    \begin{gather}\label{def-K-N}K_{N}=\on{Ker}(G(\mb O)\to G(\mc O_{N}))\end{gather} the compact open subgroup  of $G(\mb A)$ associated to the  level  $N$. We fix a lattice    $\Xi\subset Z(F)\backslash Z(\mb A)$ (where $Z$ is the center of $G$).  
   A function $f\in C_{c} (G(F)\backslash G(\mb A)/K_N \Xi,E)$ 
   is said to be   cuspidal if for all 
    parabolic subgroup $P\subsetneq G$, of  Levi $M$ and of unipotent radical  $U$,    the constant  term  $f_{P}: g\mapsto \int_{U(F)\backslash U(\mb A)}f(ug)$ vanishes as a  function on $U(\mb A)M(F)\backslash G(\mb A)/K_{N}\Xi$.  
 We recall that $C_{c}^{\rm{cusp}}(G(F)\backslash G(\mb A)/K_N \Xi,E)$ is a $E$-vector space  of finite dimension. It is endowed with a  structure of module over the Hecke algebra $C_{c}(K_{N}\backslash G(\mb A)/K_{N},E)$ : we ask that  the characteristic function  of $K_{N}$ is a unit and acts by  the identity and for $f\in   C_{c}(K_{N}\backslash G(\mb A)/K_{N},E)$ we denote by $T(f)\in \on{End}(C_{c}^{\rm{cusp}}(G(F)\backslash G(\mb A)/K_N \Xi,E))$ the corresponding Hecke operator. 
    
 \subsection{Statement of  the main theorem} 
   We will construct  the following ``excursion operators''. 
    Let $I$ be a finite set, $f$ be a function over 
$\wh G \backslash (\wh G)^{I}/\wh G    $ 
(the coarse quotient of $(\wh G)^{I}$ by  the left and right translations by the diagonal  $\wh G $), and 
$(\gamma_{i})_{i\in I}\in \on{Gal}(\ov F/F)^{I}$. 
We will construct 
the ``excursion operator''  $$S_{I,f,(\gamma_{i})_{i\in I}}\in \mr{End}_{C_{c}(K_{N}\backslash G(\mb A)/K_{N},E)}( C_{c}^{\mr{cusp}}(G(F)\backslash G(\mb A)/K_{N}\Xi,E)). 
    $$ We will show that these operators  generate a commutative subalgebra  of $\mc B$. 
    
    We do not know if  $\mc B$ is reduced but 
by spectral decomposition   we nevertheless obtain  a canonical  decomposition 
 \begin{gather}\label{intro1-nu-dec-canonical}
 C_{c}^{\mr{cusp}}(G(F)\backslash G(\mb A)/K_{N}\Xi,\Qlbar)=\bigoplus_{\nu}
 \mf H_{\nu} \end{gather}
    where the direct sum  in the RHS is indexed by the characters $\nu$ of $\mc B$, and where 
   $\mf H_{\nu}$ is the generalized eigenspace     associated to  $\nu$. We will show later that to every character $\nu$ of $\mc B$ corresponds a {\it unique}  Langlands parameter  $\sigma$ (in the sense of the  following theorem), characterized by \eqref{relation-fonda} below. 
   Setting $ \mf H_{\sigma} =\mf H_{\nu}$, we will deduce the following theorem.

 \begin{thm}  \label{intro-thm-ppal}  We have   a  canonical  decomposition of  
   $C_{c}(K_{N}\backslash G(\mb A)/K_{N},\Qlbar)$-modules
 \begin{gather}\label{intro1-dec-canonical}
 C_{c}^{\mr{cusp}}(G(F)\backslash G(\mb A)/K_{N}\Xi,\Qlbar)=\bigoplus_{\sigma}
 \mf H_{\sigma},\end{gather}
 where the direct sum  in the RHS is indexed by global  Langlands parameters, i.e.  $\wh G(\Qlbar)$-conjugacy classes    of  morphisms 
       $\sigma:\on{Gal}(\ov F/F)\to \wh G(\Qlbar)$ 
       defined  over a finite extension of  $\Ql$, continuous,  semisimple  and unramified outside   $N$.
       
     This  decomposition  is characterized by  the following property : 
       $ \mf H_{\sigma}$ is equal to the generalized eigenspace   $\mf H_{\nu}$ 
        associated to  the character $\nu$ of $\mc B$ defined  by  
       \begin{gather}\label{relation-fonda}\nu(S_{I,f,(\gamma_{i})_{i\in I}})=f((\sigma(\gamma_{i}))_{i\in I}. \end{gather}
       
      It  is compatible with the Satake isomorphism at every place  $v$  of $X\sm N$, i.e.    for every irreducible representation $V$ of $\wh G$, 
       $T(h_{V,v})$ acts on    $\mf H_{\sigma}$
       by  multiplication by  the scalar  $\chi_{V}(\sigma(\Frob_{v}))$, where $\chi_{V}$ is the character of $V$ and $\Frob_{v}$ is a   arbitrary lifting of a  Frobenius element   at $v$. 
     It is also compatible with the limit  over  $N$.   
    \end{thm}

The compatibility with the Satake isomorphism at the places     of $X\sm N$ shows that this theorem realizes the global Langlands ``correspondence''   in the direction ``from automorphic to  Galois''. 
In fact, except in the case of $GL_{r}$ \cite{laurent-inventiones},   the conjectures  of Langlands rather consist of 
\begin{itemize}
\item a parameterization, obtained in the theorem above, 
\item the Arthur multiplicity formulas for the  $\mf H_{\sigma}$, which we are not able to handle with   the methods of this article. \end{itemize}

\subsection{Main ideas of the proof} 
To construct the excursion operators and prove this  theorem the strategy will be the following. The  stacks of chtoucas,  which play a role analoguous to the Shimura varieties over number fields, exist in a much greater generality. Indeed, while  the Shimura varieties are defined   over an  open subscheme of the spectrum of the ring of integers of a number field and are associated to  {\it a  minuscule coweight } of the dual group, we possess 
for every finite set $I$, every  level $N$ and every irreducible representation $W$  of $(\wh G)^{I}$  a stack of chtoucas $\on{Cht}_{N,I,W}$  which is defined   over $(X\sm N)^{I}$. 

We will then construct  a $E$-vector space $H_{I,W}$
(where the letter $N$ is omitted to shorten the formulas) 
 as {\it a  subspace}  of the  intersection cohomology 
 of the fiber of $\on{Cht}_{N,I,W}$   over a geometric  generic point   of 
$(X\sm N)^{I}$.  By  intersection cohomology we mean here  intersection cohomology 
 with compact support, with coefficients in $E$ and  in degree $0$ (for the perverse normalization). Using the ``partial Frobenius morphisms'' introduced by Drinfeld,    we will endow  $H_{I,W}$  with an  action of $\on{Gal}(\ov F/F)^{I}$.

 \begin{rem} \label{rem-HIW-dim-finie} In this article we will define this {\it subspace} $H_{I,W}$ by a technical condition   (of finiteness  under the action of  the Hecke operators). 
 In fact  Cong Xue proved in \cite{these-cong} 
 that it has an equivalent definition as ``cuspidal'' subspace of the intersection cohomology and that it is finite dimensional. \end{rem}
 
Geometric Satake equivalence will enable us to refine this construction of  
  $H_{I,W}$  (which was defined above for any isomorphism class of irreducible representation $W$ of $(\wh G)^{I}$) into the construction of functors 
    $W\mapsto  
    H_{I,W}$, equipped with the {\it data} of the  isomorphisms 
    \eqref{isom-chi-zeta-b} below. In other words  $H_{I,W}$ is  {\it canonical} 
    and behaves well when we change $W$ and $I$.

   \begin{prop}\label{prop-a-b-c}    The $H_{I,W}$ satisfy   the 
  following properties : 
        \begin{itemize}
    \item[] {\bf a)} for every  finite set   $I$,      $$W\mapsto  
    H_{I,W},  \ \ u\mapsto \mc H(u)$$  is  a $E$-linear functor   from  the  category of  finite-dimensional  $E$-linear  representations  of  $(\wh G)^{I}$ to the  category of inductive   limits of continuous finite-dimensional $E$-linear  representations    of     $\on{Gal}(\ov F/F)^{I}$,    
              \item[] {\bf b)} for every map   $\zeta: I\to J$, 
 we possess an isomorphism 
      \begin{gather}\label{isom-chi-zeta-b}
           \chi_{\zeta}: H_{I,W}\isom 
 H_{J,W^{\zeta}},\end{gather} 
  which is 
 \begin{itemize}
 \item  functorial in   $W$, where  $W$ is a  representation of $(\wh G)^{I}$ and  $W^{\zeta}$ denotes  the   representation of $(\wh G)^{J}$ on  $W$ obtained by composition  with the  diagonal  morphism   $$  (\wh G)^{J}\to (\wh G)^{I}, (g_{j})_{j\in J}\mapsto (g_{\zeta(i)})_{i\in I} $$ 
 \item $\on{Gal}(\ov F/F)^{J}$-equivariant, where $\on{Gal}(\ov F/F)^{J}$ acts on  the LHS by the   diagonal morphism  
 \begin{gather}\label{morp-diag-Gal-intro0}\on{Gal}(\ov F/F)^{J}\to \on{Gal}(\ov F/F)^{I},  \ (\gamma_{j})_{j\in J}\mapsto (\gamma_{\zeta(i)})_{i\in I}, 
\end{gather}
 \item   and compatible with the  composition, i.e.   for  $I\xrightarrow{\zeta} J\xrightarrow{\eta} K$ we have 
 $\chi_{\eta\circ \zeta}=\chi_{\eta}\circ\chi_{\zeta}$,
  \end{itemize}
     \item[] {\bf c)} for $I=\emptyset$ and  $W=\mbf  1$, we have an  isomorphism     \begin{gather}\label{c-de-la-prop}
       H_{\emptyset,\mbf  1}=C_{c}^{\mr{cusp}}(G(F)\backslash G(\mb A)/K_{N}\Xi,E). \end{gather}
    \end{itemize}
    
  Moreover  the $H_{I,W}$ are modules over 
  $C_{c}(K_{N}\backslash G(\mb A)/K_{N},E)$, in a way  compatible  with 
the properties a), b), c) above. 
  \end{prop}

\begin{rem}
For every finite set $J$, applying b) to  the obvious map $\zeta:\emptyset \to J$, we get an isomorphism 
$\chi_{\zeta}: H_{\emptyset,\mbf  1}\isom H_{J,\mbf  1}$ 
(where $\mbf  1$ is the trivial representation  of  $(\wh G)^{J}$) 
and therefore   the action of $\on{Gal}(\ov F/F)^{J}$ over $H_{J,\mbf  1}$ is trivial. 
\end{rem}

Thanks to \eqref{c-de-la-prop}  the decomposition \eqref{intro1-dec-canonical} we are looking for  is equivalent to a decomposition 
\begin{gather}\label{dec-emptyset-1-sigma}H_{\emptyset,\mbf  1}=\bigoplus_{\sigma}
 \mf H_{\sigma}\end{gather}
 (where, increasing $E$ if necessary, we assume that the $\sigma$ and $ \mf H_{\sigma}$ are defined  over $E$).

The  following definition, where  we construct  the excursion operators,  will be repeated in  section \ref{intro-idee-heurist}. 
The reader may consult now, if he wishes, the subsection \ref{descr-conj-HIW} below for a heuristic description of the $H_{I,W}$,  which enlightens  {\it a posteriori}
the  definition of the  excursion operators. 

 We need to consider a set with one  element and we denote it by  $\{0\}$. 
For every finite set $I$ we denote by  $\zeta_{I}:I\to \{0\}$ the obvious map,  so that  $W^{\zeta_{I}}$ is nothing but  $W$ equipped with  the diagonal action   of $\wh G$. 

Thanks to the remark above (applied to $J=\{0\}$, so that  $\zeta:\emptyset\to \{0\}$ is now written  $\zeta_{\emptyset}$),and to   \eqref{c-de-la-prop}, we have  
 \begin{gather}\label{egalite-vide-0-cusp}H_{\{0\},\mbf 1}\isor{\chi_{\zeta_{\emptyset}}^{-1}} H_{ \emptyset ,\mbf 1}= C_{c}^{\mr{cusp}}(G(F)\backslash G(\mb A)/K_{N}\Xi,E).\end{gather}    
Now we will define the excursion operators as endomorphisms of \eqref{egalite-vide-0-cusp}  by using  $H_{\{0\},\mbf 1}$.

\begin{defi}  \label{defi-constr-excursion-intro} 
For every  function $f\in\mc O(\wh G\backslash (\wh G)^{I}/\wh G)$ we can find  
    a representation $W$ of $(\wh G)^{I}$,  and $x\in W$ and $\xi\in W^{*}$ invariant  by the diagonal action   of $\wh G$,  such that 
    \begin{gather}\label{f-W-x-xi-intro}    f ((g_{i})_{i\in I})=\s{\xi, (g_{i})_{i\in I}\cdot x}.\end{gather}
    Then  the endomorphism  
    $S_{I,f,(\gamma_{i})_{i\in I}}$ of  \eqref{egalite-vide-0-cusp} 
    is defined  as the composition 
  \begin{gather}\label{excursion-def-intro1}
  H_{\{0\},\mbf  1}\xrightarrow{\mc H(x)}
 H_{\{0\},W^{\zeta_{I}}}\isor{\chi_{\zeta_{I}}^{-1}} 
  H_{I,W}
  \xrightarrow{(\gamma_{i})_{i\in I}}
  H_{I,W} \isor{\chi_{\zeta_{I}}} H_{\{0\},W^{\zeta_{I}}}  
  \xrightarrow{\mc H(\xi)} 
  H_{\{0\},\mbf  1} 
    \end{gather}
    where $x:\mbf 1\to W^{\zeta_{I}}$ and $\xi: W^{\zeta_{I}}\to \mbf 1$ are considered here as   morphisms of representations of $\wh G$. 
\end{defi}

We will show   in  section \ref{intro-idee-heurist}, using properties a) and b) of   proposition  \ref{prop-a-b-c}, that 
$S_{I,f,(\gamma_{i})_{i\in I}}$   does not depend on the  choice of  $W,x,\xi$ satisfying 
\eqref{f-W-x-xi-intro} and is therefore well-defined. 
We will show in  \lemref{prop-harris} 
(found by B\"ockle,    Harris,  Khare and  Thorne) 
      that    $S_{I,f,(\gamma_{i})_{i\in I}}$ depends only on the image of  
 $(\gamma_{i})_{i\in I}$ in    $\pi_{1}(X\sm N, \ov\eta)^{I}$ (where $   \ov\eta=\on{Spec} \ov F$).  
We will show  that these   excursion operators  commute to each other, and if $\mc B$ denotes the commutative sub-algebra of  
$\on{End}(C_{c}^{\mr{cusp}}(G(F)\backslash G(\mb A)/K_{N}\Xi,E))$ which they generate, that they satisfy the following properties  (which are the expected ones because they are tautologically satisfied by the  RHS of \eqref{relation-fonda}). 

 \begin{prop} \label{prop-SIf-i-ii-iii} The  excursion operators  
   $S_{I,f,(\gamma_{i})_{i\in I}}$ satisfy the following properties:   
  \begin{itemize}
  \item [] (i) for every    $I$ and 
 $(\gamma_{i})_{i\in I}\in  \pi_{1}(X\sm N, \ov\eta) ^{I}$, 
  $$f\mapsto 
  S_{I,f,(\gamma_{i})_{i\in I}}$$ is a  morphism 
  of  commutative algebras  $\mc O(\wh G\backslash (\wh G)^{I}/\wh G)\to \mc B$, 
  \item [] (ii) for every  map 
  $\zeta:I\to J$, every  function  $f\in \mc O(\wh G\backslash (\wh G)^{I}/\wh G)$ and  every 
  $(\gamma_{j})_{j\in J}\in  \pi_{1}(X\sm N, \ov\eta) ^{J}$, we have 
  $$S_{J,f^{\zeta},(\gamma_{j})_{j\in J}}=S_{I,f,(\gamma_{\zeta(i)})_{i\in I}}$$
   where $f^{\zeta}\in \mc O(\wh G\backslash (\wh G)^{J}/\wh G)$ is defined   by    $$f^{\zeta}((g_{j})_{j\in J})=f((g_{\zeta(i)})_{i\in I}),$$
   \item [] (iii) 
  for every  $f\in \mc O(\wh G\backslash (\wh G)^{I}/\wh G)$
  and  $(\gamma_{i})_{i\in I},(\gamma'_{i})_{i\in I},(\gamma''_{i})_{i\in I}\in  \pi_{1}(X\sm N, \ov\eta) ^{I}$ we have     $$S_{I\cup I\cup I,\wt f,(\gamma_{i})_{i\in I}\times (\gamma'_{i})_{i\in I}\times (\gamma''_{i})_{i\in I}}=
  S_{I,f,(\gamma_{i}(\gamma'_{i})^{-1}\gamma''_{i})_{i\in I}}$$
   where  $I\cup I\cup I$ is a disjoint union and 
   $\wt f\in \mc O(\wh G\backslash (\wh G)^{I\cup I\cup I}/\wh G)$ is defined   by  
   $$\wt f((g_{i})_{i\in I}\times (g'_{i})_{i\in I}\times (g''_{i})_{i\in I})=f((g_{i}(g'_{i})^{-1}g''_{i})_{i\in I}).$$
      \item [] (iv) for every  $I$ and  $f$, the morphism 
      \begin{gather}\label{mor-excursion-f}\pi_{1}(X\sm N, \ov\eta) ^{I}\to \mc B, \ \ (\gamma_{i})_{i\in I}\mapsto S_{I,f,(\gamma_{i})_{i\in I}}\end{gather} is continuous, when  $\mc B$ is equipped with the  $E$-adic topology, 
      \item [] (v) for every place $v$ of $X\sm N$ and  every irreducible representation $V$ of $\wh G$, the Hecke operator $$T(h_{V,v})\in \on{End}(C_{c}^{\mr{cusp}}(G(F)\backslash G(\mb A)/K_{N}\Xi,E))$$ is equal to the excursion operator $ S_{\{1,2\}, f,(\Frob_{v},1)}$ where $f\in \mc O(\wh G\backslash (\wh G)^{2}/\wh G)$ is given  by $f(g_{1},g_{2})=\chi_{V}(g_{1}g_{2}^{-1})$, and $\Frob_{v}$ is a   Frobenius element   at $v$. 
            \end{itemize}
         \end{prop}

In fact  properties (i), (ii), (iii) and (iv) will follow formally from   proposition  \ref{prop-a-b-c} and (v) will be obtained by a geometric argument  (the computation of the composition of two  cohomological correspondences between  stacks of chtoucas). 
    
  \begin{rem} \label{rem-gamma-en-plus}
For every $\gamma\in \pi_{1}(X\sm N, \ov\eta) $, we have  
  $S_{I,f,(\gamma_{i})_{i\in I}}=S_{I,f,(\gamma_{i}\gamma)_{i\in I}}$.   
 This  is easily checked from the definition of excursion operators, or deduced from 
  the previous proposition. 
  Similarly we have $S_{I,f,(\gamma_{i})_{i\in I}}=S_{I,f,(\gamma\gamma_{i})_{i\in I}}$. 
\end{rem}
    
    In   section \ref{intro-idee-heurist} we will deduce quite easily from the previous proposition that to every   character $\nu$ of $\mc B$   corresponds a Langlands parameter
    $\sigma:\pi_{1}(X\sm N, \ov\eta) \to \wh G(\Qlbar)$ satisfying \eqref{relation-fonda}, unique up to conjugation  by 
    $\wh G(\Qlbar)$. Indeed the knowledge  of $\nu(S_{I,f,(\gamma_{i})_{i\in I}})$ (which has to be  equal to $f((\sigma(\gamma_{i}))_{i\in I}$)   for every  function $f$ determines the image of the $I$-uplet 
    $(\sigma(\gamma_{i}))_{i\in I}\in (\wh G(\Qlbar))^{I}$ as a  point defined  over $\Qlbar$ of the coarse quotient $\wh G\backslash (\wh G)^I/\wh G$. Taking $I=\{0,..,n\}$ we see that  
       $$  (\wh G)^{n}\modmod \wh G\isom \wh G\backslash (\wh G)^{\{0,...,n\}}/\wh G ,  (g_{1},...,g_{n})\mapsto 
   (1,g_{1},...,g_{n})$$ is  an  isomorphism, where  $ (\wh G)^{n}\modmod \wh G$ denotes the coarse quotient of $ (\wh G)^{n}$ by diagonal conjugation. Therefore for every integer  $n$ and every $n$-uplet 
   $(\gamma_{1},...,\gamma_{n})\in  \pi_{1}(X\sm N, \ov\eta) ^{n}$, the knowledge of $\nu$ determines the image of $(\sigma(\gamma_{1}), ..., \sigma(\gamma_{n}))$ as a  point defined  over $\Qlbar$ of the coarse quotient   $ (\wh G)^{n}\modmod \wh G$.  
    Thanks to  results of \cite{richardson} based on   geometric invariant theory, this means that  $\nu$ determines 
    $(\sigma(\gamma_{1}), ..., \sigma(\gamma_{n}))\in (\wh G(\Qlbar))^{n}$ up to  semisimplification and diagonal conjugation. Since we require  $\sigma$ to be  semisimple, it is clear  
    (choosing $n$ and $(\gamma_{1},...,\gamma_{n})$ such that 
    the subgroup generated by $\sigma(\gamma_{1}), ..., \sigma(\gamma_{n})$ is  Zariski dense in the image of $\sigma$) 
    that these date determine  $\sigma$ up to conjugation. Conversely  the relations (i), (ii), (iii) and (iv) satisfied by the excursion operators will allow, in   \propref{intro-Xi-n},  to prove the existence of $\sigma$ satisfying \eqref{relation-fonda} and property 
   (v) will ensure the compatibility with the Satake isomorphism at the places of $X\sm N$. 
    
     \subsection{A heuristic remark}  \label{descr-conj-HIW}
    This subsection suggests a conjectural description  of the 
   $H_{I,W}$,  which justifies {\it a posteriori} the definition of the  excursion operators given  in \eqref{excursion-def-intro1}. Of course this conjectural description of the  $H_{I,W}$ will never appear  in the arguments 
  and   this subsection will be used  nowhere in the rest  of the article.

  We conjecture that there exists a finite set $\Sigma$ (depending on $N$) 
 of semisimple Langlands parameters  (well defined  up to conjugation ), and that, increasing  $E$ if necessary,    
 we have for each  $\sigma\in \Sigma$
     a    $E$-linear   representation  $A_{\sigma}$ 
      of the centralizer    $S_{\sigma}$  of the image of  $\sigma$ in  $\wh G$  (trivial on $Z(\wh G)$), in such a way that for every  
      $I$ and $W$  \begin{gather}\label{intro-dec-I-W}
 H_{I,W}\overset{?}{=}\bigoplus_{\sigma\in \Sigma}
\Big(A_{\sigma} \otimes _{E}W_{\sigma^{I}}\Big)^{S_{\sigma}}, \end{gather}
 where 
$W_{\sigma^{I}}$ denotes the representation of 
$\pi_{1}(X\sm N, \ov\eta) ^{I}$     obtained by composition of the   representation $W$ with the   morphism   
$\sigma^{I}: \pi_{1}(X\sm N, \ov\eta) ^{I} \to (\wh G(E))^{I}$. 
Moreover  $A_{\sigma}$ should be a module over 
$C_{c}(K_{N}\backslash G(\mb A)/K_{N},E)$, 
and  \eqref{intro-dec-I-W} should be an  isomorphism of 
$C_{c}(K_{N}\backslash G(\mb A)/K_{N},E)$-modules. 
In the particular case where  $I= \emptyset$ and $W=1$, 
\eqref{intro-dec-I-W} should be the decomposition   
\eqref{dec-emptyset-1-sigma}  and one should have 
$ \mf H_{\sigma}=(A_{\sigma})^{S_{\sigma}}$. 
  
   There conjectures are well known by  experts, by  extrapolation of the 
conjectures  of    \cite{kottwitz3}    on the multiplicities in the cohomology of Shimura varieties, and thanks to the result of  Cong Xue mentioned in \remref{rem-HIW-dim-finie}.   
In the case of  $GL_{r}$ we expect that $\Sigma$ is the set of irreducible representations of  dimension $r$ of $\pi_{1}(X\sm N, \ov\eta)$ and that for every $\sigma\in \Sigma$, 
$S_{\sigma}=\mathbb G_{m}=Z(\wh G)$  and  $A_{\sigma}=(\pi_{\sigma})^{K_{N}}$ where  $\pi_{\sigma}$ is the cuspidal automorphic  representation corresponding to  $\sigma$ (see \cite{laurent-inventiones} and the  conjecture 2.35 of  \cite{var}). In general if $\sigma$ is associated to an elliptic Arthur parameter $\psi$ (as in subsection \ref{para-param-Arthur} below), $A_{\sigma}$ should be induced from a finite dimensional representation of the subgroup of $S_{\sigma}$ generated by the centralizer of $\psi$ and  the diagonal $\mathbb{G}_{m}\subset SL_{2}$ (because we consider only the cohomology in degree $0$). 

We conjecture moreover that    \eqref{intro-dec-I-W} is functorial in    $W$ and that for every  map $\zeta:I\to J$ it intertwines   $\chi_{\zeta}$ with 
  \begin{gather}\nonumber \Id:  \bigoplus_{\sigma} \Big(A_{\sigma} \otimes _{E}W_{\sigma^{I}}\Big)^{S_{\sigma}}\to
 \bigoplus_{\sigma} 
 \Big(A_{\sigma} \otimes _{E}(W^{\zeta})_{\sigma^{J}}\Big)^{S_{\sigma}}
   \end{gather}
   (since $W_{\sigma^{I}}$ and  $(W^{\zeta})_{\sigma^{J}}$ are both equal to  $W$ as $E$-vector spaces, the meaning of $\Id$ is clear and it is   $\pi_{1}(X\sm N, \ov\eta) ^{J}$-equivariant). 
 Under these hypotheses,   the composition \eqref{excursion-def-intro1} (which defines $S_{I,f,(\gamma_{i})_{i\in I}}$) 
acts on 
    $ \mf H_{\sigma}= (A_{\sigma})^{S_{\sigma}} \subset  H_{\{0\},\mbf 1}$ by the composition      \begin{gather}
 \nonumber (A_{\sigma})^{S_{\sigma}}
 \xrightarrow{\Id_{A_{\sigma}}\otimes x} 
\Big(A_{\sigma} \otimes _{E}W_{\sigma^{I}}\Big)^{S_{\sigma}}
\xrightarrow{\on{Id}_{A_{\sigma}}\otimes (\gamma_{i})_{i\in I}}
\Big(A_{\sigma} \otimes _{E}W_{\sigma^{I}}\Big)^{S_{\sigma}}
\xrightarrow{\Id_{A_{\sigma}}\otimes \xi} 
 (A_{\sigma})^{S_{\sigma}}
    \end{gather}
   i.e. by the product by the scalar  $\s{\xi, (\sigma(\gamma_{i}))_{i\in I} \cdot x }=f((\sigma(\gamma_{i}))_{i\in I} )$. This justifies
  {\it   a posteriori} the definition of the $S_{I,f,(\gamma_{i})_{i\in I}}$
   (and suggests that these operators are diagonalizable, but we do not know how to prove it).

   The conjecture \eqref{intro-dec-I-W} is not proven but following ideas of Drinfeld we can show that properties a) and b) of \propref{prop-a-b-c}   imply a decomposition 
 in the style of  \eqref{intro-dec-I-W} 
   (but more difficult to state, because the data of  $\Sigma$ and the  $A_{\sigma}$ should be replaced by  a  ``$\mc O$-module on the stack of global  Langlands parameters'').    More details on this construction are given in \cite{texte-ICM}.

              This subsection was  heuristic and from now on we forget  conjecture \eqref{intro-dec-I-W}. 
     
\section{Non split case, complements and  open questions }\label{para-commentaires}

\subsection{Case where $G$ is not necessarily  split}\label{cas-non-deploye}  Here we give only  the statements and we refer to  chapter 12 of \cite{coh} for the proofs,  which do not require new ideas compared to the split case. 
Let $G$ be a connected reductive  group   over $F$. 
Let $\wt F$ be a finite extension of $F$ splitting $G$ and ${}^{L}G=\wh G \rtimes \on{Gal}(\wt F/F)$ (where the semi-direct product is taken for the action of $\on{Gal}(\wt F/F)$ by automorphisms of $\wh G$ preserving a splitting). Let $U$ be an open subscheme of $X$ over which $G$ is reductive. At each point of $X\sm U$ we choose a Bruhat-Tits   parahoric model  \cite{bruhat-tits} for $G$,  so that $G$ is a smooth group scheme over $X$. It is convenient to assume that   the level $N$ is big enough so that  $X\sm N \subset U$. 
We denote by  $\Bun_{G,N}$ the   stack over $\Fq$ classifying the $G$-torsors over $X$ with structure of level $N$, in other words for every scheme $S$ over $\Fq$, $\Bun_{G,N}(S)$ is the groupoid classifying the $G$-torsors $\mc G$ over $X\times S$ equipped with a  trivialization of $\restr{\mc G}{N\times S}$. 
This stack is smooth (\cite{heinloth-unif}). 
We define  $K_{N}$ as previously  in \eqref{def-K-N}. 
We have 
 \begin{gather}\label{dec-alpha-general-ker1}  \Bun_{G,N}(\Fq)=
  \bigcup_{\alpha\in \ker^{1}(F,G)}G_{\alpha}(F)\backslash G_{\alpha}(\mb A)/K_{N}\end{gather} where the union is  disjoint, $\ker^{1}(F,G)$ is finite and $G_{\alpha}$ is the pure inner form of $G$  obtained by torsion  by $\alpha$. For every  $\alpha\in \ker^{1}(F,G)$ we have   $G_{\alpha}(\mathbb A)=G(\mathbb A)$ and therefore the  quotient by $K_{N}$ in the RHS makes sense.     
     
      We fix  a lattice $\Xi\subset Z(\mb A)/Z(F)$. 
   We define     \begin{gather}\label{dec-alpha}C_{c}^{\mr{cusp}}(\Bun_{G,N}(\Fq)/\Xi,E)= \bigoplus_{\alpha\in \ker^{1}(F,G)}C_{c}^{\mr{cusp}}(G_{\alpha}(F)\backslash G_{\alpha}(\mb A)/K_{N}\Xi,E).\end{gather} 
   Then the excursion operators are  endomorphisms 
  $$S_{I,f,(\gamma_{i})_{i\in I}}\in \mr{End}_{C_{c}(K_{N}\backslash G(\mb A)/K_{N},E)}( C_{c}^{\mr{cusp}}(\Bun_{G,N}(\Fq)/\Xi,E))$$
where  $I$ is a finite set, $(\gamma_{i})_{i\in I}\in  \pi_{1}(X\sm N, \ov\eta) ^{I}$   
 and    $f$ is a function over the coarse quotient 
   $\wh G \backslash ({}^{L } G)^{I}/\wh G    $. 
  The method  to construct them  is the same as in  the split case, thanks to a twisted variant over $X\sm N$ of the geometric Satake equivalence   (the unramified case  of \cite{richarz,zhu}). In this variant 
    ${}^{L} G$ intervenes because   the splitting of $\wh G$  appears naturally in the fiber functor of Mirkovic-Vilonen
   (indeed this  fiber functor is given by the  total cohomology  and  the  splitting  is determined by the graduation 
   by the cohomological degree  and by 
   the cup-product by the $c_{1}$ of a very ample line bundle over the affine grassmannian). 
The excursion operators  generate a commutative subalgebra  $\mc B$ and by spectral decomposition  with respect to  the characters of $\mc B$ we get a decomposition 
 \begin{gather}\label{intro2-dec-canonical}
 C_{c}^{\mr{cusp}}
 (\Bun_{G,N}(\Fq)/\Xi,\Qlbar)
 =\bigoplus_{\sigma}  \mf H_{\sigma}. \end{gather} 
 The  direct sum  in the RHS is indexed by global Langlands parameters, 
  i.e.  the $\wh G(\Qlbar)$-conjugacy  classes  of  morphisms 
       $\sigma:\on{Gal}(\ov F/F)\to {}^{L} G(\Qlbar)$ 
       defined  over a finite extension of  $\Ql$, continuous,  semisimple,  unramified outside   $N$ and giving rise to the commutative diagram  
        \begin{gather}\label{diag-sigma}
 \xymatrix{
\on{Gal}(\ov F/F) \ar[rr] ^{\sigma}
\ar[dr] 
&& {}^{L} G(\Qlbar) \ar[dl] 
 \\
& \on{Gal}(\wt F/F) }\end{gather}
As in \thmref{intro-thm-ppal}  
the  decomposition \eqref{intro2-dec-canonical} is characterized by  
\eqref{relation-fonda}, and it is compatible with the (twisted) Satake isomorphism   \cite{satake,cartier-satake,borel-corvallis,blasius-rogawski-pspm} at all  places of $X\sm N$.   

\begin{rem} Usually (for instance in the Arthur multiplicity formulas) we take the quotient of the set of  morphisms $\sigma$ 
by  a weaker  equivalence  relation, which, in addition to the conjugation by $\wh G(\Qlbar)$, allows to twist  $\sigma$ by   elements of $\ker^{1}(F,Z(\wh G)(\Qlbar))$. 
 By  Kottwitz~\cite{kottwitz1,kottwitz2} (and theorem 2.6.1 of  Nguyen Quoc Thang~\cite{thang}   for the adaptation to characteristic $p$), this finite group  is the dual of $\ker^{1}(F,G)$. Therefore it has the same  cardinal  and from the point of view of  Arthur multiplicity formulas our finer equivalence relation on  $\sigma$ compensates exactly the fact that  our space is a sum indexed by $ \alpha\in \ker^{1}(F,G)$. 
 For example, if $G$ is a torus, the subspaces $\mf H_{\sigma}$ of \eqref{intro2-dec-canonical} are of dimension $1$. 
 \end{rem}
 
 \begin{rem} When $G$ is split, $\ker^{1}(F,Z(\wh G)(\Qlbar))$ is $0$ by the theorem of Tchebotarev, hence 
  $\ker^{1}(F,G)$ is $0$ also and  
$\Bun_{G,N}(\Fq)=G (F)\backslash G (\mb A)/K_{N}$. This is why  the quotient $G(F)\backslash G(\mb A)/K_{N}\Xi$ appeared in  section \ref{para-enonces}, and will appear again in 
 section \ref{defi-chtou-intro}  and afterwards (where we will come back to the split case to simplify the redaction).  
\end{rem}

 We refer to    proposition 12.5  
 of \cite{coh} for the fact that the decomposition \eqref{intro2-dec-canonical} is compatible with the  isogenies of $G$ (and more generally with  all  the morphisms $G\to G'$ whose image is normal). 

We hope that the decomposition \eqref{intro2-dec-canonical}  is  compatible with all the known cases of functoriality  given ky an explicit kernel. In particular we should be able to prove it is  compatible with 
  the    theta correspondence, thanks to  the geometrization of the  theta kernel   by Lysenko \cite{sergey-theta,sergey-theta-SO-Sp} and to the link between our  construction and  \cite{brav-var} (explained in   section  \ref{subsection-link-langl-geom}).

\subsection{Arthur parameters} \label{para-param-Arthur} 
We would like to prove that the  Langlands parameters $\sigma$  which appear in the decomposition \eqref{intro1-dec-canonical} (or \eqref{intro2-dec-canonical} in  the non split case)  come from elliptic Arthur parameters.  We recall that an Arthur parameter  is a $\wh G(\Qlbar)$-conjugacy class     of  morphism $$
  \psi : \on{Gal}(\ov F/F) \times SL_{2}(\Qlbar)\to 
 {}^{L} G(\Qlbar) \text{ (algebraic over  $SL_{2}(\Qlbar)$),}$$ whose restriction to   $\on{Gal}(\ov F/F)$ takes its values  in a finite extension of  $\Ql$, is continuous, almost everywhere unramified, semisimple, pure of weight $0$ and  makes a triangle similar to  \eqref{diag-sigma} commute.  
 Moreover  $ \psi$ is said to be elliptic if the centralizer 
   of $\psi$ in  $\wh  G(\Qlbar)$ is finite modulo $(Z(\wh  G)(\Qlbar))^{\on{Gal}(\wt F/F)}$. 
  
The   Langlands parameter associated to $\psi$ is  $\sigma_{\psi}: \on{Gal}(\ov F/F) \to {}^{L} G(\Qlbar)$ defined  by  
 $$\sigma_{\psi}(\gamma)=\psi\Big(\gamma, \begin{pmatrix} |\gamma|^{1/2} & 0 \\
 0 & |\gamma|^{-1/2}
 \end{pmatrix}\Big)$$  where $|\gamma|^{1/2}$ is well defined  thanks to the choice of a square root  of $q$.   We conjecture  that every Langlands  parameter  $\sigma$ occuring   in the  decomposition 
 \eqref{intro1-dec-canonical} (or \eqref{intro2-dec-canonical} in the  non split case)
    is of the form   $\sigma_{\psi}$ 
     with  $\psi$ an elliptic Arthur parameter unramified on $X\sm N$.   By \cite{kostant-betti} the $\wh G(\Qlbar)$-conjugacy class  of $\psi$ is uniquely determined by the $\wh G(\Qlbar)$-conjugacy class of $\sigma$.

    We would like to obtain  a  canonical  decomposition as  \eqref{intro1-dec-canonical} (or \eqref{intro2-dec-canonical} in the non split case)  for {\it the whole discrete part} (and  not only the cuspidal part) and this  decomposition should be indexed by elliptic Arthur parameters. 
       
\subsection{Meaning of the decomposition}   The decomposition 
\eqref{intro1-dec-canonical} 
is finer in general than the one  obtained  by diagonalization  
of the Hecke operators at the unramified places. Even the isomorphism classes  of representations of 
 $C_{c}(K_{N}\backslash G(\mb A)/K_{N},\Qlbar)$  do not allow to determine in general the decomposition \eqref{intro1-dec-canonical}, and although  the  Arthur multiplicity formulas are stated with a sum indexed by Arthur parameters, such a   canonical  decomposition seems unkown in general in the case of number fields.     Indeed after  Blasius \cite{blasius}, Lapid  \cite{lapid} and Larsen 
    \cite{larsen1,larsen2}, for some groups  $ G $ (including split ones)  the same  representation of 
 $C_{c}(K_{N}\backslash G(\mb A)/K_{N},\Qlbar)$ may occur in differents subspaces  $\mf H_{\sigma}$, because of the following phenomenon. There are examples of finite groups  $ \Gamma $ and of morphisms $ \tau, \tau ': \Gamma \to \wh G (\Qlbar) $ such that $ \tau $ and $ \tau' $ are not conjugated but that  for every  $ \gamma \in \Gamma $, $ \tau (\gamma) $ and $ \tau '(\gamma) $ are  conjugated.
 We expect then there could exist a surjective everywhere unramified   morphism   $ \rho: \on {Gal} (\ov F / F) \to \Gamma $   and a representation $ (H_ {\pi}, \pi) $  of $ G (\mb A) $ such that $ (H_ {\pi}) ^ {K_ {N}} $ occurs in both  $ \mf H_ {\tau \circ \rho} $ and $ \mf H_ {\tau' \circ \rho} $.
 The examples of Blasius and Lapid are for $ G = SL_{r} $, $ r \geq 3$
  (in fact in this case we can recover  {\it a posteriori} the decomposition \eqref{intro1-dec-canonical} thanks to the embedding   $SL_{r}\hookrightarrow GL_{r}$). 
   But for some groups  (for instance $E_{8}$)  
 we do not know how to recover  the decomposition \eqref{intro1-dec-canonical} by other means than the  methods of the present  article,  which work only for function fields.

 \subsection{Independence of $\ell$} We hope that  the decomposition   \eqref{intro1-dec-canonical} (or \eqref{intro2-dec-canonical} in the non split case) is defined  over $\ov \Q$, independent of $\ell$ (and of the embedding $\ov \Q\subset \Qlbar$), and indexed by motivic Langlands parameters (the  notion of 
 motivic Langlands parameter is clear if we admit the Standard Conjectures but we note that in \cite{drinfeld-pro-completion} Drinfeld 
 gave an unconditionnal definition). This conjecture seems out of reach for the moment. We refer to   conjecture   12.12  of \cite{coh} for a more precise statement. 
   
  \subsection{Case of number fields} 
It is obviously out of reach to apply the methods of this article to number fields. However one can hope a decomposition analoguous to the {\it canonical}  decomposition  
  \eqref{intro1-dec-canonical} (or \eqref{intro2-dec-canonical} in the non split case) and  an Arthur multiplicity formula for each of the spaces   $\mf H_{\sigma}$. 
When  $F$ is a function field as in this article, 
the limit    $ \varprojlim_{N } \Bun_{G,N}(\Fq)$ is equal to   $\big(G(\ov F)\backslash G(\mb A\otimes_{F}\ov F)   \big)^{\on{Gal}(\ov F/F)}$. This expression still makes sense for number fields and to understand its topology we just notice it is equal 
\begin{itemize}\item   to 
    $\big(G(\check F)\backslash G(\mb A\otimes_{F}\check F)\big)^{\on{Gal}(\check F/F)}$ where  $\check F$ is a  finite Galois extension   of $F$ over which  $G$ is split
    \item to 
    $ \bigcup_{\alpha\in \ker^{1}(F,G)}G_{\alpha}(F)\backslash G_{\alpha}(\mb A)
$ 
where $\ker^{1}(F,G)$ is finite and $G_{\alpha}$ is a inner form of $G$.     \end{itemize} 
  We can hope that if  $F$ is a number field  and if 
$\Xi$ is a lattice  in  $Z(F)\backslash Z(\mb  A)$, 
the discrete part     \begin{gather}\label{L2disc}L^{2}_{\mr{disc}}\Big(\big(G(\ov F)\backslash G(
\mb A\otimes_{F}\ov F)\big)^{\on{Gal}(\ov F/F)}/\Xi,\C\Big)\end{gather} admits a {\it canonical}    decomposition   indexed by the $\wh G(\C)$-conjugacy classes   of 
elliptic Arthur parameters. 
 The particular case which is the most similar to the case of function fields 
 is the case of cohomological automorphic forms. 
  Indeed the cohomological part  of \eqref{L2disc} is defined   over  $\ov\Q$ and we can hope that it admits  a {\it canonical}   decomposition   over  $\ov \Q$ indexed by equivalence classes of elliptic 
  motivic Langlands parameters  (with the subtleties  of \cite{buzzard-gee} about the difference between  $L$-algebricity and $C$-algebricity). 
   
\subsection{Importance of the sum over $\ker^{1}$} 
 In the non split case we don't know  if the inclusion 
 $$
  C_{c}^{\mr{cusp}}(G (F)\backslash G (\mb A)/K_{N}\Xi,\Qlbar)\subset C_{c}^{\mr{cusp}}(\Bun_{G,N}(\Fq)/\Xi,\Qlbar)$$  (corresponding to the case  $\alpha=0$ in the RHS of \eqref{dec-alpha-general-ker1}) 
   is compatible with the decomposition \eqref{intro1-dec-canonical}, even after we regroup together the  $\sigma$ which differ by an element of   $\ker^{1}(F,Z(\wh G)(\Qlbar))$. So, except in the case  when 
   the sum \eqref{dec-alpha-general-ker1} is reduced to one term
   (i.e. when $\ker^{1}(F,Z(\wh G)(\Qlbar))=0$, and in particular when 
    $G$  is an  inner form of a  split group), 
  we do not obtain a    canonical decomposition 
 of  the space  $
  C_{c}^{\mr{cusp}}(G (F)\backslash G (\mb A)/K_{N}\Xi,\Qlbar)$. 
   
\subsection{Coefficients in finite fields} Thanks to the fact that  the geometric Satake equivalence is defined  over $\mc O_{E}$ \cite{mv,ga-de-jong}, we can prove that when the function  $f$ is defined  over $\mc O_{E}$, $S_{I,f,(\gamma_{i})_{i\in I}}$ preserves 
  $ C_{c}^{\mr{cusp}}(\Bun_{G,N}(\Fq)/\Xi,\mc O_{E})$. By spectral  decomposition of the reduction of these operators modulo the maximal ideal  of $\mc O_{E}$ we get a decomposition  of 
   $ C_{c}^{\mr{cusp}}(\Bun_{G,N}(\Fq)/\Xi,\Flbar))$ indexed by the classes of $\wh G(\Flbar)$-conjugation of morphisms $\sigma : \pi_{1}(X\sm N, \ov\eta)\to {}^{L} G(\Flbar)$ defined  over a finite  field, continuous, making the diagram analoguous to       \eqref{diag-sigma} commute, and completely reducible in the sense  of Jean-Pierre Serre \cite{bki-serre,bmr} (i.e.  if the image is included in a  parabolic subgroup of ${}^{L} G$, it is included in an associated  Levi subgroup). We refer to  chapter 13 of \cite{coh} for more details. 
   
\subsection{Local paramerization}  The local  Langlands parameterization (up to semisimplification) and  the local-global compatibility are proven  in an article   with Alain Genestier \cite{genestier-lafforgue}. 
They are deduced from  the following  statement   : if $v$ is a place of $X$   and if all the $\gamma_{i}$ belong to $\on{Gal}(\ov F_{v}/F_{v})$, then  $S_{I,f,(\gamma_{i})_{i\in I}}$ is equal to the action of an element of the ($\ell$-adic completion of the)  Bernstein center  of $G(F_{v})$ (obviously the interesting case is when  $v\in N$ because the unramified case is completely solved by the compatibility with the Satake isomorphism in  \thmref{intro-thm-ppal}).  
   
\subsection{Multiplicities} This work does not determine the multiplicities, and in particular it does not say for which   Langlands parameters $\sigma$ 
   the space  $\mf H_{\sigma}   $ is non zero. 
          
     \subsection{Case of metaplectic groups}
     In chapter  14 de \cite{coh}  we sketch how to extend our results  to metaplectic groups, thanks to the metaplectic variant of the geometric Satake equivalence established in 
     \cite{finkelberg-lysenko, lysenko-red, dennis-sergey}.       
          
    \section{
    Chtoucas of Drinfeld for the  reductive groups, after Varshavsky 
   } \label{defi-chtou-intro}
       
       In the rest of this article $G$ is split
       (the non split case,  which was discussed in  subsection \ref{cas-non-deploye} above,  
      is handled in  chapter 12 of \cite{coh} but no new idea is necessary).        
       The geometric ingredients  of our  construction are explained in this section and in section \ref{subsection-crea-annihil-intro}. 
       Here is short overview. The intersection  cohomology  with compact support of the stacks  of chtoucas provides  for every finite  set $I$, for every  level $N$ and for every  representation $W$ of $(\wh G)^{I}$ an inductive system   $\varinjlim_{\mu}\mc H_{N,I,W}^{\leq\mu}$ of constructible  $E$-sheaves  over  $(X\sm N)^{I}$. The goal of this  section is to construct this inductive system,  functorially in  $W$,  and to equip it with actions of the Hecke operators and of the  partial Frobenius morphisms $F_{\{i\}}$, and to establish 
        the coalescence isomorphisms  \eqref{intro-isom-coalescence}  which describe its  restriction by a diagonal  morphism   $(X\sm N)^{J}\to (X\sm N)^{I}$              (associated to an arbitrary map $I\to J$). 
       In  section \ref{subsection-crea-annihil-intro} we will show 
       that the Hecke operators at unramified places  can be rewritten with the help of   coalescence isomorphisms  and  partial Frobenius morphisms. It is this property that will ensure  the compatibility of our construction with the Satake isomorphism at unramified places.  It will also play a fundamental technical role by allowing to extend  the Hecke operators to  morphisms of sheaves over the whole $(X\sm N)^{I}$, and by providing 
        the Eichler-Shimura relations. These relations   (which will be stated 
         in   \propref{Eichler-Shimura-intro} below) 
claim that for each place  of $X\sm N$ and for each  $i\in I$ the restriction at $x_{i}=v$ of  the partial Frobenius morphism $F_{\{i\}}$ is killed by a polynomial whose  coefficients are   Hecke operators at $v$ (with coefficients in $\mc O_{E}$).  They will be used in  section \ref{section-sous-sheaves-constr} to prove that the property of finiteness under   the action of the  Hecke operators (which gives the definition of  the $H_{I,W}$) implies a 
property of finiteness  under the action of  the partial Frobenius morphisms, whence, thanks to a fundamental lemma  of Drinfeld, the action of $\on{Gal}(\ov F/F)^{I}$ over $H_{I,W}$. We warn the reader that this use of the  Eichler-Shimura relations  
is  completely unusual. 
       
    The chtoucas were introduced by  Drinfeld  \cite{drinfeld78,Dr1} for  $GL_{r}$    and  generalized  to arbitrary  reductive groups (and arbitrary  coweights) by  Varshavsky in  
    \cite{var} (meanwhile the case of division algebras was considered   by  Laumon-Rapoport-Stuhler, Laurent Lafforgue, Ngô Bao Châu and  Eike Lau, see the  references at the beginning  of section \ref{intro-previous-works}).
        Let $I$  be a finite  set 
 and  $W$ be a  $E$-linear irreducible  representation  of $(\wh G)^{I}$. We write $W=\boxtimes_{i\in I}W_{i}$ where $W_{i}$ is an irreducible representation of  $\wh G$.  
 The  stack $\Cht_{N,I,W}^{(I)}$ classifying the  $G$-chtoucas with structure of  level $N$, was studied  in  \cite{var}.  
  Contrary to  \cite{var} we require  it to be reduced in the following definition   
  (this does not matter for étale cohomology). 
  
  \noindent{\bf Notation. }  For every scheme  $S$ over $\Fq$  and for every 
  $G$-torsor $\mc G$ over  $ X\times S$ we write   $\ta \mc G=(\Id_{X}\times \Frob_{S})^{*}(\mc G)$. 
  
\begin{defi} We define  $\Cht_{N,I,W}^{(I)}$ as the {\it reduced}  Deligne-Mumford  stack  over $(X\sm N)^{I}$ whose  points   over a scheme  $S$ over $\Fq$ classify    
  \begin{itemize}
  \item  points  $(x_{i})_{i\in I}: S\to (X\sm N)^{I}$,  
\item a    $G$-torsor $\mc G$ over  $ X\times S$, 
\item  an isomorphism 
$$\phi :\restr{\mc G }{(X\times S)\sm(\bigcup_{i\in I }\Gamma_{x_i})}\isom \restr{\ta \mc G}{(X\times S)\sm(\bigcup_{i\in I }\Gamma_{x_i})}$$ 
 where   $\Gamma_{x_i}$ denotes the graph  of $x_{i}$, such that the  relative position  at   $x_{i}$ is bounded by the   dominant coweight  of  $G$ corresponding to the dominant weight  $\omega_{i}$ of   $W_{i}$, 
 \item a trivialization of   $(\mc G,\phi)$  over  $N\times S$. 
 \end{itemize}
 \end{defi}
 
 This  definition will be  generalized  in   \defiref{defi-Cht-I1-Ik} below, and the condition about the  relative position  will be made more precise  in  \remref{position-relative}. 
   
 We write  $\Cht_{ I,W}^{(I)}$ when $N$ is empty and we note that   $\Cht_{N,I,W}^{(I)}$ is a $G(\mc O_{N})$-torsor 
 over $\restr{\Cht_{ I,W}^{(I)}}{(X\sm N)^{I}}$.

 \begin{rem}  Readers who know the geometric Langlands program will recognize that  $\Cht_{N,I,W}^{(I)}$ is the fiber product over  $\Bun_{G,N}\times \Bun_{G,N}$ of a Hecke stack  (considered as  a correspondence from $\Bun_{G,N}$ to itself) with the graph of the Frobenius morphism of $\Bun_{G,N}$. 
 \end{rem}

The  $x_{i}$ will be called the legs of the  chtouca. We will denote by    $$\mf p_{N,I,W} ^{(I)}: \Cht_{N,I,W} ^{(I)}\to (X\sm N)^{I}$$ the corresponding morphism. 
 
 For every   dominant coweight   $\mu$ of $G^{\mr{ad}}$ we denote by 
$\Cht_{N,I,W}^{(I),\leq\mu}$ the open substack of  $\Cht_{N,I,W}^{(I)}$ defined  by  the condition that  the Harder-Narasimhan polygon    of $\mc G$ 
(or rather, to be precise, of its associated $G^{\mr{ad}}$-torsor)  is $\leq \mu$ (by which we mean that the difference is a  combination with nonnegative rational  coefficients of simple coroots).

We fix a lattice   $\Xi\subset Z(F)\backslash Z(\mb A)$.  
Then  $\Xi$ maps to  $\Bun_{Z,N}(\Fq)$ which acts $\Cht_{N,I,W}^{(I)}$ by torsion, and preserves the open substacks  $\Cht_{N,I,W}^{(I),\leq\mu}$. 

We can show that   $\Cht_{N,I,W}^{(I),\leq\mu}/\Xi$ is a  Deligne-Mumford stack    of finite type. It is equipped with the morphism   
$$\mf p_{N,I,W} ^{(I),\leq\mu}: \Cht_{N,I,W}^{(I),\leq\mu}/\Xi\to (X\sm N)^{I}$$
 deduced from $\mf p_{N,I,W} ^{(I)}$. 

 We denote by  $\on{IC}_{\Cht_{N,I,W}^{(I),\leq \mu}/\Xi}$ the IC sheaf   of  $\Cht_{N,I,W}^{(I),\leq \mu}/\Xi$ with coefficients in  $E$, normalized relatively to     $(X\sm N)^{I}$. 
The  following definition will be made more  canonical (and in particular functorial in $W$) in   \defiref{defi-HNIW-can}. 

\begin{defi} \label{defi-HNIW-naive}
We set 
$$\mc H_{N,I,W}^{\leq\mu}=R^{0}\big(\mf p_{N,I,W} ^{(I),\leq\mu}\big)_{!}
 \Big(\on{IC}_{\Cht_{N,I,W}^{(I),\leq \mu}/\Xi}\Big).$$
  \end{defi}

 Compared to    \cite{coh} we have shortened the notations, by removing the indices  which recalled that  $\mc H_{N,I,W}^{\leq\mu}$ is a  cohomology   in degree $0$ (for the perverse normalization) and  with coefficients in $E$. 
 The cohomology is taken in the sense  of \cite{laumon-moret-bailly,laszlo-olsson} but in fact the only background we need is the étale cohomology of schemes. 
 Indeed, as soon as  the degree of $N$ is big enough in  function of $\mu$, $\Cht_{N,I,W}^{(I),\leq \mu}/\Xi$ is a scheme of finite type. 
 Therefore for any open subscheme $U\subset X\sm N$ such that $U\subsetneq X$ (in order to be able to increase $N$ without changing  $U$),
  $\restr{\Cht_{N,I,W}^{(I),\leq \mu}/\Xi}{U^{I}}$ is the quotient of a scheme of finite type by a finite group.  
 
 When  $I$ is empty  and $W=\mbf 1$,  we have  
\begin{gather}\label{I-vide-intro} \restr{\varinjlim_{\mu}\mc H_{N,\emptyset,\mbf 1}^{\leq\mu}}{\Fqbar}=C_{c}(G(F)\backslash G(\mb A)/K_N \Xi,E)\end{gather}
 because   $\Cht_{N,\emptyset,\mbf 1}$ is the discrete stack  
  $\Bun_{G,N}(\Fq)$, considered as  a constant stack   over $\Fq$, and  moreover  $\Bun_{G,N}(\Fq)=G(F)\backslash G(\mb A)/K_N$
 (here we use the hypothesis that    $G$ is split, in general  by \eqref{dec-alpha-general-ker1}, $\Bun_{G,N}(\Fq)$ would be a finite union of adélic quotients of $G$).  
 
 \begin{rem} \label{rem-cht-W1} More generally  for every $I$ and  $W=\mbf  1$, the stack  $\Cht_{N,I,{\mbf  1} }^{(I)}/\Xi$ is simply  the constant stack  
  $G(F)\backslash G(\mb A)/K_{N}\Xi$ over $(X\sm N)^{I}$.  \end{rem}
    
We consider  $\varinjlim_{\mu}\mc H_{N,I,W}^{\leq\mu}$ as an inductive system of  constructible $E$-sheaves  over  $(X\sm N)^{I}$.
We will now   define   the actions of the  partial Frobenius morphisms and the  Hecke operators over this inductive system  (we warn the reader that these  actions increase  $\mu$).  
For every subset    $J\subset I$ we denote by   $$
 \Frob_{J}:(X\sm N)^{I}\to (X\sm N)^{I}$$  the  morphism  which sends   
 $(x_{i})_{i\in I}$ to    $(x'_{i})_{i\in I}$ with   $$x'_{i}=\Frob(x_{i})\text{ \  if \  }
  i\in J\text{   \ and  \  }x'_{i}=x_{i} \text{ \ otherwise.}$$  
 Then we have 
 \begin{itemize}
 \item for $\kappa$ big enough  and for every   $i\in I$, a  morphism  
 \begin{gather}\label{intro-action-Frob-partiel}F_{\{i\}}:\Frob_{\{i\}}^{*}(\mc H_{N,I,W}^{\leq\mu})\to 
 \mc H_{N,I,W}^{\leq\mu+\kappa}\end{gather} 
 of  constructible sheaves over   $  (X\sm N)^{I}$, in such a way that  the 
 $F_{\{i\}}$ commute with each other and  their product for $i\in I$ is the natural action  of the total  Frobenius morphism   of $(X\sm N)^{I}$ on the sheaf $\mc H_{N,I,W}^{\leq\mu}$, 
  \item for every   $f\in C_{c}(K_{N}\backslash G(\mb A)/K_{N},E)$  and for $\kappa$ big enough, a  morphism 
 \begin{gather}\label{defi-Tf}T(f):\restr{\mc H_{N,I,W}^{\leq\mu}}{(X\sm \mf P)^{I}}\to 
 \restr{\mc H_{N,I,W}^{\leq\mu+\kappa}}{(X\sm \mf P)^{I}}\end{gather} of  constructible sheaves over   $  (X\sm \mf P)^{I}$ where $\mf P$ is a finite  set of places containing  $|N|$ and outside of which    $f$ is trivial.  
 \end{itemize}
 
 The morphisms $T(f)$ are called   ``Hecke operators'' although they are morphisms of sheaves. 
They are obtained thanks to the   obvious construction of Hecke  correspondences  between the stacks of chtoucas. 
  We will see  after  \propref{prop-coal-frob-cas-part-intro} that  $T(f)$ can be extended naturally to  a  morphism of sheaves over   $(X\sm N)^{I}$,  but this is not trivial. 
  Of course when $I=\emptyset$ and $W=\mbf 1$, the morphisms $T(f)$ are the usual  Hecke operators   on \eqref{I-vide-intro}. 
 
  To construct the actions \eqref{intro-action-Frob-partiel} of the partial Frobenius morphisms, we need a small generalization of the stacks  
 $\Cht_{N,I,W} ^{(I)}$ where we ask a factorization of  $\phi$ as a composition of several  modifications. Let $(I_{1},...,I_{k})$ be  an (ordered)  partition    of $I$. 
 
 \begin{defi}\label{defi-Cht-I1-Ik}
 We define      $\Cht_{N,I,W} ^{(I_{1},...,I_{k})}$ as the {\it reduced }   Deligne-Mumford  stack    whose points over 
 a scheme   $S$ over  $\Fq$  classify   \begin{gather}\label{intro-donnee-chtouca}\big( (x_i)_{i\in I}, (\mc G_{0}, \psi_{0}) \xrightarrow{\phi_{1}}  (\mc G_{1}, \psi_{1}) \xrightarrow{\phi_{2}}
\cdots\xrightarrow{\phi_{k-1}}  (\mc G_{k-1}, \psi_{k-1}) \xrightarrow{ \phi_{k}}    (\ta{\mc G_{0}}, \ta \psi_{0})
\big)
\end{gather}
with 
 \begin{itemize}
\item $x_i\in (X\sm N)(S)$ for $i\in I$, 
\item for $i\in \{0,...,k-1\}$, $(\mc G_{i}, \psi_{i})\in \Bun_{G,N}(S)$ (i.e.   
$\mc G_{i}$ is a   $G$-torsor over  $X\times S$ and 
$\psi_{i} : \restr{\mc G_{i}}{N\times S} 
   \isom 
   \restr{G}{N\times S}$ is a trivialization over   $N\times S$) and we  note  
   $(\mc G_{k}, \psi_{k})=(\ta{\mc G_{0}}, \ta \psi_{0})$
 \item  
for   $j\in\{1,...,k\}$
 $$\phi_{j}:\restr{\mc G_{j-1}}{(X\times S)\sm(\bigcup_{i\in I_{j}}\Gamma_{x_i})}\isom \restr{\mc G_{j}}{(X\times S)\sm(\bigcup_{i\in I_{j}}\Gamma_{x_i})}$$ is an   isomorphism  such that the  relative position  of  $\mc G_{j-1}$ w.r.t. $\mc G_{j}$ at  $x_{i}$ (for  $i\in I_{j}$) is bounded by the   dominant coweight  of  $G$ corresponding to the  dominant weight  of  $W_{i}$, 
 \item the $\phi_{j}$,  which induce isomorphisms over $N\times S$, respect  the level structures, i.e.   
$\psi_{j}\circ \restr{\phi_{j}}{N\times S}=\psi_{j-1}$ for every  $j\in\{1,...,k\}$. 
 \end{itemize}
 \end{defi}
The condition on the  relative position  will be made more precise  in  \remref{position-relative}. We denote by  $\Cht_{N,I}^{(I_{1},...,I_{k})}$ the indstack obtained when we forget this  condition. 
 
For any dominant coweight $\mu$ of  $G^{\mr{ad}}$   we denote  by 
 $\Cht_{N,I,W}^{(I_{1},...,I_{k}),\leq\mu}$   the open substack of  $\Cht_{N,I,W}^{(I_{1},...,I_{k})}$ defined  by  the condition that the Harder-Narasimhan polygon   of  $\mc G_{0}$ is  $\leq\mu$. We denote by    $$\mf p_{N,I,W} ^{(I_{1},...,I_{k})}: \Cht_{N,I,W} ^{(I_{1},...,I_{k})}\to (X\sm N)^{I}$$ the morphism    which associates to a  chtouca the family of its legs. 
   
    \noindent{\bf Exemple.} When  $G=GL_r$,  $I=\{1,2\}$ and $W=\mr{St}\boxtimes \mr{St}^{*}$, the stacks 
 $\Cht_{N,I,W}^{(\{1\}, \{2\})}$, {\it resp.} 
  $\Cht_{N,I,W}^{(\{2\}, \{1\})}$ 
 are the stacks of left, {\it resp.}  right chtoucas   introduced by  Drinfeld (and used  also in  \cite{laurent-inventiones}), 
 and $x_{1}$ and  $x_{2}$ are the pole and the zero.

We  construct now a     smooth morphism  
\eqref{lisse-chtouca-grass-intro} from  $\Cht_{N,I,W}^{(I_{1},...,I_{k})}$ to  the  quotient of a closed stratum of a  Beilinson-Drinfeld affine grassmannian  by a smooth group scheme  (in addition this will allow us  in   \remref{position-relative} to formulate  more precisely  the condition on the relative positions in   \defiref{defi-Cht-I1-Ik}).   
Readers familiar  with    Shimura varieties may consider this morphism as a ``local model'' provided they note 
\begin{itemize}
\item that we are in a situation of good reduction since the $x_{i}$ belong to $X\sm N$, 
\item and that however this local model is not smooth (except if all  the $I_{j}$ are   singletons and all  the coweights are minuscule). 
\end{itemize}

\begin{defi}
The  Beilinson-Drinfeld affine grassmannian is the indscheme  $\mr{Gr}_{I }^{(I_{1},...,I_{k})}$ over $X^{I}$ whose $S$-points classify    
\begin{gather}\label{formule-rem-grassm-intro}\big((x_{i})_{i\in I}, \mc G_{0} \xrightarrow{\phi_{1}}  
\mc G_{1}\xrightarrow{\phi_{2}}
\cdots\xrightarrow{\phi_{k-1}}  \mc G_{k-1} \xrightarrow{ \phi_{k}}   \mc G_{k}\isor{\theta} G_{X\times S} \big)  \end{gather}
where the $\mc G_{i}$ are   $G$-torsors over $X\times S$, $\phi_{i}$ is an   isomorphism over $(X\times S)\sm(\bigcup_{i\in I_{j}}\Gamma_{x_i})$ and $\theta$ is a trivialization of $\mc G_{k}$. 
 The closed stratum    $\mr{Gr}_{I,W}^{(I_{1},...,I_{k})}$ is the reduced closed  subscheme   of 
  $\mr{Gr}_{I }^{(I_{1},...,I_{k})}$ defined  by  the condition that 
  the  relative position  of  $\mc G_{j-1}$ w.r.t. $\mc G_{j}$ at  $x_{i}$ (for  $i\in I_{j}$) is  bounded by the   dominant coweight  of  $G$ corresponding to the   dominant weight  $\omega_{i}$ of  $W_{i}$. More precisely 
  over  the open subscheme $\mc U$ of $X^{I}$ where the $x_{i}$ are pairwise  distinct, $\mr{Gr}_{I }^{(I_{1},...,I_{k})}$ is a product  of usual affine grassmannians   and 
  \begin{itemize}
  \item we define 
  the restriction of $\mr{Gr}_{I,W}^{(I_{1},...,I_{k})}$ over  $\mc U$  as the product of the usual closed strata (denoted  $\ov {\mr{Gr}_{\omega_{i}}}$ in  \cite{brav-gaitsgory,mv}), 
  \item then we define  $\mr{Gr}_{I,W}^{(I_{1},...,I_{k})}$ as the Zariski closure (in $\mr{Gr}_{I }^{(I_{1},...,I_{k})}$) of 
  its  restriction over  $\mc U$. 
  \end{itemize}
  \end{defi}
  
   By Beauville-Laszlo \cite{BL} (see also the first  section   of \cite{coh} for additional references in \cite{hitchin}), 
   $\mr{Gr}_{I }^{(I_{1},...,I_{k})}$ can also be  defined  as the indscheme whose $S$-points classify  
       \begin{gather}\label{formule-rem-grassm-intro-loc}\big((x_{i})_{i\in I}, \mc G_{0} \xrightarrow{\phi_{1}}  
\mc G_{1}\xrightarrow{\phi_{2}}
\cdots\xrightarrow{\phi_{k-1}}  \mc G_{k-1} \xrightarrow{ \phi_{k}}   \mc G_{k}\isor{\theta} G_{\Gamma_{\sum \infty x_i}} \big)  \end{gather}
 where the $\mc G_{i}$ are  $G$-torsors  over the formal neighborhood 
  $\Gamma_{\sum \infty x_i}$ of the union of the graphs of the  $x_{i}$ in 
 $  X\times S$, $\phi_{i}$ is an   isomorphism over $\Gamma_{\sum \infty x_i}\sm(\bigcup_{i\in I_{j}}\Gamma_{x_i})$ and $\theta$ is a trivialization of $\mc G_{k}$. 
The restriction à la  Weil $G_{\sum \infty x_i}$ of $G$ from  $\Gamma_{\sum \infty x_i}$ to  $S$ acts therefore  on  $\mr{Gr}_{I }^{(I_{1},...,I_{k})}$  and on $\mr{Gr}_{I,W}^{(I_{1},...,I_{k})}$ by change  of the trivialization $\theta$. 

We have a natural morphism  
    \begin{gather}\label{lisse-chtouca-grass-intro0}
    \Cht_{N,I,W}^{(I_{1},...,I_{k})}\to \mr{Gr}_{I,W}^{(I_{1},...,I_{k})}/
  G_{\sum \infty  x_i}\end{gather} 
   which associates  to  a chtouca \eqref{intro-donnee-chtouca} the 
  $G_{\sum \infty x_{i}}$-torsor $\restr{\mc G_{k}}{\Gamma_{\sum \infty x_{i}}}$
  and, for any trivialization $\theta$ of it, the  point of $\mr{Gr}_{I,W}^{(I_{1},...,I_{k})}$ equal to  \eqref{formule-rem-grassm-intro-loc}.

 \begin{rem} \label{position-relative} The best way to formulate  the condition on the relative positions in the \defiref{defi-Cht-I1-Ik} is to {\it define}   $\Cht_{N,I,W}^{(I_{1},...,I_{k})}$ as the  inverse image of 
  $\mr{Gr}_{I,W}^{(I_{1},...,I_{k})}/
  G_{\sum \infty x_i}$ by the morphism   
  $\Cht_{N,I}^{(I_{1},...,I_{k})}\to \mr{Gr}_{I}^{(I_{1},...,I_{k})}/
  G_{\sum \infty x_i}$ constructed as in  \eqref{lisse-chtouca-grass-intro0}. 
    \end{rem}

For $(n_{i})_{i\in I}\in \N^{I}$ we denote by $\Gamma_{\sum n_{i} x_i}$ the closed subscheme  of $X\times S$ associated to a Cartier divisor $\sum n_{i} x_i$  which is  effective and relative  over $S$. We denote by  $G_{\sum n_{i} x_i}$ the smooth group scheme over $S$ obtained by restriction à la Weil of $G$ from $\Gamma_{\sum n_{i} x_i}$ to $S$. 
 Then if the integers $n_{i}$ are big enough in function of $W$, the action of $G_{\sum \infty x_i}$ on $\mr{Gr}_{I,W}^{(I_{1},...,I_{k})}$ factorizes through 
  $G_{\sum n_{i} x_i}$. 
Then the   morphism \eqref{lisse-chtouca-grass-intro0} provides   a morphism 
    \begin{gather}\label{lisse-chtouca-grass-intro}
    \Cht_{N,I,W}^{(I_{1},...,I_{k})}\to \mr{Gr}_{I,W}^{(I_{1},...,I_{k})}/
  G_{\sum n_{i} x_i}\end{gather} 
  (which associates to a chtouca \eqref{intro-donnee-chtouca} the 
  $G_{\sum n_{i}x_{i}}$-torsor $\restr{\mc G_{k}}{\Gamma_{\sum n_{i}x_{i}}}$
  and, for any trivialization $\lambda$ of it, the  point of $\mr{Gr}_{I,W}^{(I_{1},...,I_{k})}$ equal to  \eqref{formule-rem-grassm-intro-loc} for any trivialization $\theta$ of $\restr{\mc G_{k}}{\Gamma_{\sum \infty x_i}}$ extending   $\lambda$ of $\Gamma_{\sum n_{i} x_i}$ to 
  $\Gamma_{\sum \infty   x_i}$).

 We show   in   proposition 2.8 of \cite{coh} that  the morphism 
  \eqref{lisse-chtouca-grass-intro} 
    is  smooth of  dimension $\dim G_{\sum n_{i}x_{i}}=(\sum_{i\in I} n_{i})\dim G$.  
    
 We deduce from this that  the  morphism (forgetting intermediate  modifications) 
 \begin{gather}
 \label{oubli-Cht}
 \Cht_{N,I,W}^{(I_{1},...,I_{k})}\to \Cht_{N,I,W}^{(I)} \\ \nonumber 
 \text{    which sends \eqref{intro-donnee-chtouca} to } 
\big( (x_i)_{i\in I}, (\mc G_{0}, \psi_{0}) \xrightarrow{\phi_{k} \cdots \phi_{1}}      (\ta{\mc G_{0}}, \ta \psi_{0})
\big) \end{gather}
 is small. Indeed it is known that the analoguous  morphism  
\begin{gather}\label{mor-Gr-oubli}\mr{Gr}_{I,W}^{(I_{1},...,I_{k})}\to \mr{Gr}_{I,W}^{(I)}  \text{     which sends \eqref{formule-rem-grassm-intro-loc} to }
\big( (x_i)_{i\in I},  \mc G_{0}  \xrightarrow{\phi_{k} \cdots \phi_{1}}       \mc G_{k}\isor{\theta} G_{\Gamma_{\sum \infty x_i}} 
\big)  \end{gather}   
is small, and by the way this fact plays an essential role   in   \cite{mv}. Moreover the inverse image  of $\Cht_{N,I,W}^{(I),\leq\mu}$ 
by  \eqref{oubli-Cht} is exactly  $\Cht_{N,I,W}^{(I_{1},...,I_{k}),\leq\mu}$ 
since the truncatures were defined with the help of the  Harder-Narasimhan polygon of $\mc G_{0}$. 
 Therefore   $$\mc H_{N,I,W}^{\leq\mu}=
   R^{0}(\mf p_{N,I,W} ^{(I_{1},...,I_{k}),\leq\mu})_{!}\Big(\on{IC}_{\Cht_{N,I,W}^{(I_{1},...,I_{k}),\leq\mu}/\Xi}\Big)$$ for {\it every} partition $(I_{1},...,I_{k})$ of $I$ (whereas  \defiref{defi-HNIW-naive} used the coarse partition  $(I)$).  

The   partial Frobenius morphism 
$$\on {Fr}_{I_{1}} ^{(I_{1},...,I_{k})}: \Cht_{N,I,W} ^{(I_{1},...,I_{k})}\to \Cht_{N,I,W} ^{(I_{2},...,I_{k},I_{1})},$$ 
 defined  by  
\begin{gather}\on {Fr}_{I_{1}} ^{(I_{1},...,I_{k})}\big( (x_i)_{i\in I}, (\mc G_{0}, \psi_{0}) \xrightarrow{\phi_{1}}  (\mc G_{1}, \psi_{1}) \xrightarrow{\phi_{2}}
\cdots\xrightarrow{\phi_{k-1}}  (\mc G_{k-1}, \psi_{k-1}) \xrightarrow{ \phi_{k}}    (\ta{\mc G_{0}}, \ta \psi_{0})
\big)\nonumber \\
= \nonumber
\big( \Frob_{I_{1}}\big((x_i)_{i\in I}\big), (\mc G_{1}, \psi_{1}) \xrightarrow{\phi_{2}}  (\mc G_{2}, \psi_{2}) \xrightarrow{\phi_{3}}
\cdots \xrightarrow{ \phi_{k}}    (\ta{\mc G_{0}}, \ta \psi_{0}) \xrightarrow{\ta  \phi_{1}   } (\ta{\mc G_{1}}, \ta \psi_{1})
\big)  
\end{gather}
lies  over  the morphism  $\Frob_{I_{1}}:(X\sm N)^{I}\to (X\sm N)^{I}$. 
The composition of the  morphisms $\on {Fr}_{I_{1}} ^{(I_{1},...,I_{k})}$, $\on {Fr}_{I_{2}} ^{(I_{1},...,I_{k})}$, ..., $\on {Fr}_{I_{k}} ^{(I_{1},...,I_{k})}$ is equal to the total   Frobenius morphism of  $\Cht_{N,I,W} ^{(I_{1},...,I_{k})}$ over  $\Fq$. 
Since   $\on {Fr}_{I_{1}} ^{(I_{1},...,I_{k})}$
is a  completely radicial local homeomorphism, we have a canonical  
isomorphism   
\begin{gather}\label{action-Frob-partiels}
\Big(\on {Fr}_{I_{1}} ^{(I_{1},...,I_{k})}\Big)^{*}\Big(\on{IC}_{\Cht_{N,I,W} ^{(I_{2},...,I_{k},I_{1})}}\Big)=\on{IC}_{ \Cht_{N,I,W} ^{(I_{1},...,I_{k})}}   \end{gather}
which we normalize by a  coefficient $q^{-d/2}$, where $d$ is the relative  dimension of $\mr{Gr}_{I_{1},\boxtimes _{i\in I_{1}}W_{i}}^{(I_{1})}$ over  $X^{I_{1}}$ 
(this normalization will be justified in  \remref{rem-Frob-partiel}). 
Thanks to the  proper base   change isomorphism 
(and to the fact that completely radicial local homeomorphisms do not alter étale topology), or by an easy cohomological correspondence given in \cite{coh}, 
we deduce from \eqref{action-Frob-partiels} a morphism 
\begin{gather}\label{morph-F1}F_{I_{1}}:\Frob_{I_{1}}^{*}(\mc H_{N,I,W}^{\leq\mu})\to 
 \mc H_{N,I,W}^{\leq\mu+\kappa}\end{gather}
  for $\kappa$ big enough
 (indeed, with the notations above, if the Harder-Narasimhan polygon of $\mc G_{1}$ is $\leq \mu$, since the modification between $\mc G_{0}$ and $\mc G_{1}$ is bounded in function of $W$, the Harder-Narasimhan polygon of $\mc G_{0}$ is $\leq \mu+\kappa$ where $\kappa$ depends on  $W$). By taking any  partition $(I_{1},...,I_{k})$  such that   $I_{1}=\{i\}$ we get  $F_{\{i\}}$ in \eqref{intro-action-Frob-partiel}.

 For the moment  we have defined  $\mc H_{N,I,W}^{\leq\mu}$ for the  isomorphisms classes of irreducible representations  $W$ of $(\wh G)^{I}$. We will refine this  construction into the more  {\it canonical}   construction of a $E$-linear {\it  functor} 
   \begin{gather}\label{fonc-W-HIW-intro}W\mapsto \mc H_{N,I,W}^{\leq\mu} \end{gather} 
from  the category of  finite dimensional $E$-linear representations of $(\wh G)^{I}$ to  the category of  constructible $E$-sheaves  over $(X\sm N)^{I}$. In  particular for every    morphism $u:W\to W'$ of  $E$-linear representations  of $(\wh G)^{I}$  we will denote by 
 $$\mc H(u): \mc H_{N,I,W}^{\leq\mu} \to \mc H_{N,I,W'}^{\leq\mu}$$  
the associated  morphism of  constructible sheaves. 

 The  functor \eqref{fonc-W-HIW-intro} will be  compatible with the coalescence of legs, in the following sense. In all this article we call coalescence the situation where legs are glued to each other. 
 We could have used the word fusion instead of  coalescence
 but we prefered to use the word  coalescence for the legs
 (which are just points over the curve) while keeping the word fusion 
 for the fusion product (which occurs in the geometric Satake equivalence and involves the perverse sheaves  over the affine grassmannians of Beilinson-Drinfeld). 
     Let $\zeta: I \to J$ be any map.  We denote by  $W^{\zeta}$ the representation of  $\wh G^{J}$  which is the composition of the  representation $W$ with the diagonal morphism
 $$ \wh G^{J}\to \wh G^{I}, (g_{j})_{j\in J}\mapsto (g_{\zeta(i)})_{i\in I}.$$
   We denote by  
   \begin{gather}\label{morph-giad-X-intro}\Delta_{\zeta}: X^{J}\to X^{I},(x_{j})_{j\in J}\mapsto (x_{\zeta(i)})_{i\in I}\end{gather} the   diagonal morphism. 
We will construct, after   \cite{var} and \cite{brav-var},    an isomorphism  of  constructible sheaves over $(X\sm N)^{J}$, called the  {\it coalescence isomorphism}:  \begin{gather}\label{intro-isom-coalescence}\chi_{\zeta}: \Delta_{\zeta}^{*}(\mc H_{N,I,W}^{\leq\mu})\isom 
 \mc H^{\leq\mu}_{N,J,W^{\zeta}}.\end{gather} 
 This isomorphism will be  {\it canonical} in the sense that it will be an isomorphism of functors in $W$, compatible with composition of $\zeta$. 
 Now let us  explain the construction of the functor \eqref{fonc-W-HIW-intro} and of the coalescence isomorphism  \eqref{intro-isom-coalescence}.

 When $W$ is not  irreducible we denote by  $\mr{Gr}_{I,W }^{(I_{1},...,I_{k})}$  the union of the $\mr{Gr}_{I,V }^{(I_{1},...,I_{k})}\subset \mr{Gr}_{I }^{(I_{1},...,I_{k})}$ for irreducible constituents $V$  of $W$. We do the same  with 
  $\Cht_{N,I,W}^{(I_{1},...,I_{k})}$.

   We recall now the geometric Satake equivalence, of Lusztig, Drinfeld, Ginzburg and Mirkovic-Vilonen. 
   For references we quote   \cite{lusztig-satake,ginzburg,hitchin,mv,ga-iwahori,ga-de-jong,richarz,zhu}. Here we will use it in the form explained   by Gaitsgory in \cite{ga-de-jong}.  Usually the geometric Satake equivalence can be formulated in the following way. 
   For every algebraically closed field $k$ of characteristic prime to  $\ell$, 
   the category of  $G(k[[z]])$-equivariant  perverse sheaves  over the affine grassmannian $G(k((z)))/G(k[[z]])$ is equipped  with 
      \begin{itemize}
    \item
    a   tensor structure by the fusion product (or the  convolution product), 
    with a modification of the signs in the commutativity constraint that we recall in   
     \remref{modif-comm} below, 
\item      a  fiber functor,  given  by the  total cohomology, 
which, thanks to this  modification is a tensor functor with values in the category of vector spaces, and not super-vector spaces   (to be more canonical  we have to  introduce 
    a Tate twist, i.e. tensor by    $E(\frac{i}{2})$ the part of cohomological    degree  $i$ for every $i\in \Z$). 
\end{itemize}
   Thus this tensor category   is   equivalent to the category of  representations 
   of the group of automorphisms of the fiber functor,  which happens to be isomorphic to 
   $\wh G$ (and equipped with a canonical splitting). 
   Moreover its objects  are naturally equivariant under the action of the automorphism group of $k[[z]]$. This allows to replace 
   $\on{Spf}(k[[z]])$ by an arbitrary formal disk, and in particular a formal disk moving on a curve.

  Here we use only one direction of this  equivalence, namely the functor from  the category of  representations of $\wh G$ to  the category of $G(k[[z]])$-equivariant  perverse sheaves  
    over the affine grassmannian. On the other hand we state it  with the help of the  affine grassmannian of Beilinson-Drinfeld. 
  The fact that the fiber     functor 
  which gives rise to the geometric Satake equivalence is given the total cohomology implies,   in the notations of theorem below, that $W$ is canonically equal to the total  cohomology    of $\mc S_{I,W }^{(I_{1},...,I_{k})}$ in the fibers of $\mr{Gr}_{I }^{(I_{1},...,I_{k})}$ over  $X^{I}$.
  However we do not use this property as it is, only through the fact that it gives   the fiber functor of  the geometric Satake equivalence, whence the canonicity of our constructions.

 \begin{thm}  \label{satake-geom-thm} (one direction  of the geometric Satake equivalence \cite{hitchin, mv, ga-de-jong}). We have for every  finite  set $I$ and for every partition $(I_{1},...,I_{k})$ of $I$ a      $E$-linear functor 
   $$W\mapsto \mc S_{I,W }^{(I_{1},...,I_{k})}$$
from  the category of   finite dimensional $E$-linear representations of $(\wh G)^{I}$ to  the  category of  $G_{\sum \infty x_{i}}$-equivariant     perverse $E$-sheaves   
 over $\mr{Gr}_{I }^{(I_{1},...,I_{k})}$ (for the perverse normalization  relative to  $X^{I}$). Moreover $\mc S_{I,W}^{(I_{1},...,I_{k})}$ is supported by 
$\mr{Gr}_{I,W}^{(I_{1},...,I_{k})}$ and is universally locally acyclic  relatively to $X^{I}$. These  functors satisfy the following properties.  

\begin{itemize}
\item [] a) Compatibility with the morphisms which forget the intermediate modifications : $S_{I,W}^{(I)}$ is canonicaly isomorphic to the direct image of 
$\mc S_{I,W}^{(I_{1},...,I_{k})}$ by the forgetful morphism  
$\mr{Gr}_{I}^{(I_{1},...,I_{k})}\to \mr{Gr}_{I}^{(I)}$ (defined  in \eqref{mor-Gr-oubli}).   
\item [] b)  Compatibility with    convolution : if $W=\boxtimes_{j\in \{1,...,k\}} W_{j}$ where $W_{j}$ is a  representation of $(\wh G)^{I_{j}}$, 
  $\mc S_{I,W}^{(I_{1},...,I_{k})}$  is canonicaly isomorphic to the inverse image of $\boxtimes _{j\in \{1,...,k\}} \mc S_{I_{j},W_{j}}^{(I_{j})}$ by the morphism 
  \begin{align*}
  \mr{Gr}_{I}^{(I_{1},...,I_{k})}/G_{\sum_{i\in I} \infty x_{i}} &\to 
  \prod_{j=1}^{k}\Big( \mr{Gr}_{I_{j}}^{(I_{j})}/G_{\sum_{i\in I_{j}} \infty x_{i}}\Big) \\
  (\mc G_{0}\to \mc G_{1} \to  \cdots \to  \mc G_{k})& 
  \mapsto \Big(\big(\restr{\mc G_{j-1}}{\Gamma_{\sum_{i\in I_{j}} \infty x_{i}}}\to \restr{\mc G_{j}}{\Gamma_{\sum_{i\in I_{j}} \infty x_{i}}} \big) \Big)_{j=1,...,k}
  \end{align*}
  where the $\mc G_{i}$ are   $G$-torsors over 
  $\Gamma_{\sum_{i\in I} \infty x_{i}}$. 
   \item [] c)  Compatibility with  fusion: let  $I,J$ be finite sets  and  $\zeta: I\to J$ be any map. Let $(J_{1},...,J_{k})$ be a partition of $J$. Its inverse image 
   $(\zeta^{-1}(J_{1}), ..., \zeta^{-1}(J_{k}))$ is a partition of $I$. 
We denote by   $$\Delta_{\zeta} : X^{J}\to X^{I},  \ \ (x_{j})_{j\in J}\mapsto (x_{\zeta(i)})_{i\in I}$$ the   diagonal morphism associated to 
  $\zeta$. We denote again by  
  $\Delta_{\zeta}$    the inclusion  
 $$\mr{Gr}_{J}^{(J_{1},...,J_{k})} =\mr{Gr}_{I}^{(\zeta^{-1}(J_{1}), ..., \zeta^{-1}(J_{k}))}\times _{X^{I}}X^{J}\hookrightarrow \mr{Gr}_{I}^{(\zeta^{-1}(J_{1}), ..., \zeta^{-1}(J_{k}))}.$$
    Let $W$ be a   finite dimensional $E$-linear representation    of  $\wh G^{I}$. We denote by   $W^{\zeta}$ the   representation of $\wh G^{J}$  which is the composition of the   representation $W$ with the    diagonal morphism $$\wh G^{J}\to \wh G^{I},  \ \ (g_{j})_{j\in J}\mapsto (g_{\zeta(i)})_{i\in I}  .$$   Then we have   a   canonical isomorphism
 \begin{gather}\label{coalescence-Gr-section1-thm} \Delta_{\zeta}^{*}\Big( \mc S_{I,W}^{(\zeta^{-1}(J_{1}), ..., \zeta^{-1}(J_{k}))} \Big)\simeq 
 \mc S_{J,W^{\zeta}}^{(J_{1},...,J_{k})}\end{gather}  
  which is functorial in   $W$ and compatible with the composition for $\zeta$. 
  \item [] d) 
 When  
$W$ is irreducible, the perverse sheaf  $\mc S_{I,W}^{(I_{1},...,I_{k})}$ over 
 $\mr{Gr}_{I,W}^{(I_{1},...,I_{k})}$ is isomorphic to its   IC-sheaf   (with the  perverse normalization relative to   $X^{I}$). 
 \end{itemize}  
\end{thm}
  
  The properties a) and b) could have been stated  with  more general partitions, but at the cost of cumbersome notations. 
  
 In the previous theorem  $\mc S_{I,W}^{(I_{1},...,I_{k})}$ is supported by 
$\mr{Gr}_{I,W}^{(I_{1},...,I_{k})}$ and therefore we can consider it  as a perverse sheaf 
(up to a shift)   over  $\mr{Gr}_{I,W}^{(I_{1},...,I_{k})}/G_{\sum n_{i}x_{i}}$  (if  the integers $n_{i}$ are  big enough).

    \begin{rem}\label{modif-comm}
 The commutativity constraint is defined by the  {\it modified} convention explained in the  discussion before   proposition 6.3 of  \cite{mv} (and also in   \cite{hitchin}). Here is a brief summary. 
  Lemma 3.9 of  \cite{mv}  shows that 
   \begin{itemize}
   \item 
   for any connected component of the   affine grassmannian,
   the strata are all of even dimension, or all of odd dimension 
   (then we say that the component is even or odd)
       \item if a perverse sheaf  $\mc S_{I_{j},W_{j},\mc A}^{(I_{j})}$ is supported on a even ({\it resp.} odd) component,   its total cohomology is concentrated  in even ({\it resp.} odd) cohomological degrees. 
   \end{itemize}
 The modified commutativity constraint 
 is obtained by adding to the usual Koszul signs a minus sign when we permute two sheaves $\mc S_{I_{j},W_{j},\mc A}^{(I_{j})}$ and  $\mc S_{I_{j'},W_{j'},\mc A}^{(I_{j'})}$  both supported on odd components  of the affine grassmannian.  
Otherwise stated it is the commutativity constraint  we would obtain naturally if 
we normalized the sheaves  $\mc S_{I_{j},W_{j},\mc A}^{(I_{j})}$ 
so that their total cohomology   be concentrated in even cohomological degrees.      \end{rem}

 \begin{rem}   The  compatibility between the classical   Satake isomorphism and the geometric Satake equivalence  is expressed by the fact that, if  $v$ is a place in $X$ of residual field $k(v)$ and  $V$ is an  irreductible representation of  $\wh G$ of highest weight  $ \omega$  and
 $\rho$ denotes the half-sum of positive coroots of  $\wh G$, 
 $(-1)^{\s{2\rho, \omega}}h_{V,v}$ is equal to the trace of  $\Frob_{\mr{Gr}_{v}/k(v)} $  (where $\mr{Gr}_{v}$ is the affine grassmannian at  $v$)
 on the perverse sheaf  $\restr{\mc S_{\{0\},V,E}^{(\{0\})}}{\mr{Gr}_{v}}$ 
 (which is the IC-sheaf of the closed stratum  $\mr{Gr}_{v,\omega}$).   The sign  $(-1)^{\s{2\rho, \omega}}$   in the above definition of  $h_{V,v}$ comes from the fact that the stratum  $\mr{Gr}_{v,\omega}$ has  dimension $\s{2\rho, \omega}$. In other words  $h_{V,v}$ would be equal to the trace of  $\Frob_{\mr{Gr}_{v}/k(v)} $ on   $\restr{\mc S_{\{0\},V,E}^{(\{0\})}}{\mr{Gr}_{v}}$ 
if we had normalized it so that  its total cohomology is concentrated in even cohomological degrees (but without changing the  Tate torsion).
This choice is coherent with the fact that the modified commutativity constraint 
 of   \cite{mv} would have been the natural one with such normalizations, as we recalled in the previous remark. 
 The coherence of this choice is the reason why there will be no sign in the equality of   \propref{prop-coal-frob-cas-part-intro}. 
 \end{rem}

    \begin{rem}\label{rem-Satake-Gad}
    In the previous theorem $Z_{\sum_{i\in I} \infty x_{i}}\subset 
G_{\sum_{i\in I} \infty x_{i}}$ acts trivially on  
 $ \mr{Gr}_{I}^{(I_{1},...,I_{k})}$ and therefore (by d)) on the sheaves  $\mc S_{I,W }^{(I_{1},...,I_{k})}$. We write  $G^{\mr{ad}}=G/Z$. Thus we can consider  $\mc S_{I,W}^{(I_{1},...,I_{k})}$  as a $G^{\mr{ad}}_{\sum \infty x_{i}}$-equivariant perverse sheaf  (up to a shift)
 on  $\mr{Gr}_{I }^{(I_{1},...,I_{k})}$
 or also as a  perverse sheaf  (up to a shift)   on  $\mr{Gr}_{I,W}^{(I_{1},...,I_{k})}/G^{\mr{ad}}_{\sum n_{i}x_{i}}$  (for  $n_{i}$   big enough).   
    \end{rem}

Here is  the construction of the  functor \eqref{fonc-W-HIW-intro}. 
The morphism \eqref{lisse-chtouca-grass-intro} does not factorize   through the quotient by $\Xi$ (as noticed by an anonymous referee), but it is true for its composition with the forgetful morphism $ \mr{Gr}_{I,W}^{(I_{1},...,I_{k})}/
  G_{\sum n_{i} x_i}\to  \mr{Gr}_{I,W}^{(I_{1},...,I_{k})}/
  G^{\mr{ad}}_{\sum n_{i} x_i}$ (because $\Xi$ acts by twisting by  $Z$-torsors). In other words we get a  morphism 
  \begin{gather}\label{lisse-chtouca-grass-intro-ad}
    \Cht_{N,I,W}^{(I_{1},...,I_{k})}/\Xi\to \mr{Gr}_{I,W}^{(I_{1},...,I_{k})}/
  G^{\mr{ad}}_{\sum n_{i} x_i} \end{gather} 
and by \remref{rem-Satake-Gad}, $\mc S_{I,W}^{(I_{1},...,I_{k})}$  is a perverse sheaf (up to a shift) on the target space.

\begin{defi}\label{defi-HNIW-can}
We 
 define the  perverse sheaf (with the normalization relative to $(X\sm N)^{I}$) 
 $\mc F_{N,I,W}^{(I_{1},...,I_{k})}$ on 
 $\Cht_{N,I,W}^{(I_{1},...,I_{k}) }/\Xi$ 
 as the inverse image of 
 $\mc S_{I,W}^{(I_{1},...,I_{k})}$ by the morphism 
   \eqref{lisse-chtouca-grass-intro-ad}.    
  Then  we define    the functor \eqref{fonc-W-HIW-intro} by setting 
  \begin{gather}\label{defi-can-H}\mc H_{N,I,W}^{\leq\mu}=
   R^{0}(\mf p_{N,I,W} ^{(I_{1},...,I_{k}),\leq\mu})_{!}\Big(\restr{\mc F_{N,I,W}^{(I_{1},...,I_{k})}}{\Cht_{N,I,W}^{(I_{1},...,I_{k}), \leq\mu}/\Xi}\Big)\end{gather} for any partition $(I_{1},...,I_{k})$ of $I$. 
   \end{defi}
  Thanks to  a) of the previous theorem, 
   the definition  \eqref{defi-can-H} does not depend on the 
    choice of  the partition  $(I_{1},...,I_{k})$. 
 
 When $W$ is irreducible, the smoothness of the morphism \eqref{lisse-chtouca-grass-intro}, and therefore the smoothness  of  the morphism \eqref{lisse-chtouca-grass-intro-ad},  and the computation of its dimension imply that   
   $\mc F_{N,I,W}^{(I_{1},...,I_{k})}$ is isomorphic to the  IC-sheaf of 
   $\Cht_{N,I,W}^{(I_{1},...,I_{k})}/\Xi$ (with the perverse normalization relative to $(X\sm N)^{I}$).   Thus the previous definition is compatible  with 
     definition \ref{defi-HNIW-naive} (which it refines and makes more canonical).

  With the help of  \eqref{defi-can-H} we can reformulate the action of Hecke operators  in an obvious way, as well as the action of the    partial Frobenius morphisms  in a way indicated in the next remark, which the reader can skip. 
    
 \begin{rem}\label{rem-Frob-partiel}
Let  $W=\boxtimes_{j\in \{1,...,k\}} W_{j}$ where $W_{j}$ is a  representation of  $(\wh G)^{I_{j}}$. 
 Then the action of the    partial Frobenius morphisms on  \eqref{defi-can-H}    still takes the form of a morphism like  \eqref{morph-F1}   and comes, in the same way as \eqref{morph-F1},   from an isomorphism 
 \begin{gather}\label{action-Frob-partiels-general}\Big(\on {Fr}_{I_{1}} ^{(I_{1},...,I_{k})}\Big)^{*}\Big(\mc F_{N,I,W} ^{(I_{2},...,I_{k},I_{1})} \Big)
 \isom \mc F_{N,I,W} ^{(I_{1},...,I_{k})}  \end{gather}
which  was already defined in   \eqref{action-Frob-partiels} for irreducible   $W$. Thanks to the more canonical definition of  $\mc F_{N,I,W}^{(I_{1},...,I_{k})}$ given  in  \defiref{defi-HNIW-can}, the isomorphism  \eqref{action-Frob-partiels-general} admits the following more natural alternative definition. 
  
  Thanks to   \defiref{defi-HNIW-can} and to  b) of  \thmref{satake-geom-thm}, $\mc F_{N,I,W}^{(I_{1},...,I_{k})}$ is isomorphic to the inverse image of  $\boxtimes _{j\in \{1,...,k\}} \mc S_{I_{j},W_{j}}^{(I_{j})}$ by the natural smooth morphism 
 $$  \Cht_{N,I,W}^{(I_{1},...,I_{k}), \leq\mu}/\Xi \to 
  \prod_{j=1}^{k}\Big( \mr{Gr}_{I_{j}}^{(I_{j})}/G^{\mr{ad}}_{\sum_{i\in I_{j}} n_{i} x_{i}}\Big) 
  $$ (where the integers $n_{i}$ are big enough) and writing 
  $\Frob_{1}$ 
  the  Frobenius  morphism    of   $\mr{Gr}_{I_{1},W_{1}}^{(I_{1})}/G^{\mr{ad}}_{\sum _{i\in I_{1}}n_{i}x_{i}}$ and 
   \begin{gather}\label{Frob-Gr-1} \on F_{1}:  
 \Frob_{1}^{*}\big( \mc S_{I_{1},W_{1},E}^{(I_{1})}
\big) \isom \mc S_{I_{1},W_{1},E}^{(I_{1})} 
\end{gather} 
the natural isomorphism,  the  diagram 
\begin{gather*}
\begin{CD}
 \Cht_{N,I,W} ^{(I_{1},...,I_{k})} /\Xi@>{\on  {Fr}_{I_{1}, N,I} ^{(I_{1},...,I_{k})}}>> \Cht_{N,I,W} ^{(I_{2},...,I_{1})}/\Xi \\
 @VVV @VVV \\
\prod_{j=1}^{k}
\mr{Gr}_{I_{j},W_{j}}^{(I_{j})}/G^{\mr{ad}}_{\sum _{i\in I_{j}}n_{i}x_{i}}
 @>{\Frob_{1}
\times \Id}>>  \prod_{j=1}^{k} \mr{Gr}_{I_{j},W_{j}}^{(I_{j})}/G^{\mr{ad}}_{\sum _{i\in I_{j}}n_{i}x'_{i}}
 \end{CD}
 \end{gather*}   is commutative and makes compatible 
 \eqref{action-Frob-partiels-general} and  
   $$F_{1}\times \Id : (\Frob_{1}\times \Id )^{*}\big(\boxtimes _{j\in \{1,...,k\}} \mc S_{I_{j},W_{j}}^{(I_{j})}
\big) \isom \boxtimes _{j\in \{1,...,k\}} \mc S_{I_{j},W_{j}}^{(I_{j})},  
   $$ which gives a more natural alternative definition of   \eqref{action-Frob-partiels-general} (in particular the factor $q^{-d/2}$ comes from $F_{1}$). 
 \end{rem}

 However the canonicity of  definition \eqref{defi-can-H} is mostly crucial  for the construction of the coalescence isomorphisms  
 \eqref{intro-isom-coalescence}, which will be explained now,  
 because the source and target spaces are not the same.  
     
 \begin{defi}  The {\it canonical} isomorphism \eqref{intro-isom-coalescence} is defined  (thanks to the  proper base change theorem) by the {\it canonical}  isomorphism 
   between 
   $\mc F_{N,J,W^{\zeta}}^{(J_{1},...,J_{k})}$   and the inverse image of 
 $\mc F_{N,I,W}^{(\zeta^{-1}(J_{1}), ..., \zeta^{-1}(J_{k}))}$ by the inclusion 
   $$\Cht_{N,J }^{(J_{1},...,J_{k})}
     =\Cht_{N,I }^{(\zeta^{-1}(J_{1}), ..., \zeta^{-1}(J_{k}))}\times _{(X\sm N)^{I}}(X\sm N)^{J}\hookrightarrow \Cht_{N,I }^{(\zeta^{-1}(J_{1}), ..., \zeta^{-1}(J_{k}))} 
     $$
      which comes from  the isomorphism \eqref{coalescence-Gr-section1-thm} in   c) of \thmref{satake-geom-thm}. 
     \end{defi}
     The previous definition  is  independent on the choice of  the partition $(J_{1},...,J_{k})$.  The fact that we took an arbitrary  partition  allows us to prove
    the following proposition. 
  
   \begin{prop} 
    \label{rem-coalescence-frob-compat} 
    The  actions of the partial Frobenius morphisms and the Hecke operators are  compatible with the  morphisms $\mc H(u)$ and with the coalescence  isomorphisms. 
\end{prop}
  \noindent The compatibility between the partial Frobenius morphisms and the coalescence isomorphism  \eqref{intro-isom-coalescence}  means  that for every  $j\in J$, $\Delta_{\zeta}^{*}(F_{\zeta^{-1}(\{j\})})$ and 
    $F_{\{j\}}$ correspond to each other  by the isomorphism $\chi_{\zeta}$ 
     of \eqref{intro-isom-coalescence}.

       \section{Proof of \propref{prop-a-b-c} }
       \label{section-esquisse-abc}

 We call a geometric point  $\ov x$ of a scheme  $Y$ the data of 
 an algebraically closed field  $k(\ov x)$ and a   morphism $\on{Spec}(k(\ov x))\to Y$.
 In this article $k(\ov x)$ will always be an algebraic closure of the residue  field $k(x)$ of the point $x\in Y$ below   $\ov x$. 
  We denote by   $Y_{(\ov x)}$ the strict localization  (or    strict henselianization) of  $Y$ at   $\ov x$. 
  In other words  $Y_{(\ov x)}$ is the spectrum of the ring   $\varinjlim \Gamma(U,\mc O_{U})$, where the  inductive limit  is taken over the     $\ov x$-pointed étale neighborhoods of  $x$. It is a  local henselian ring   whose  residue field is the  separable closure  of  $k(x)$ in  $k(\ov x)$. 
 If  $\ov x$ and $\ov y$ are two geometric points of $Y$,  we call   a specialization arrow $\on{\mf{sp}}:\ov x\to \ov y$   a  morphism  
 $Y_{(\ov x)}\to Y_{(\ov y)}$, or equivalently a    morphism  
 $\ov x\to Y_{(\ov y)}$ (such an arrow exists if and only if $y$ is in the Zariski closure of $x$). By section 7 of \cite{grothendieck-sga4-2-VIII} a   specialization arrow $\on{\mf{sp}}:\ov x\to \ov y$ induces for every sheaf    $\mc F$ over a open subscheme of $Y$ containing $y$ a specialization homomorphism $\on{\mf{sp}}^{*}: \mc F_{\ov y}\to \mc F_{\ov x}$.

  We fix an algebraic closure $\ov F$ of  $F$ and we  denote by 
      $\ov \eta=\on{Spec}(\ov F)$ the corresponding geometric point over   the generic point  $\eta$    of $X$.  
 
  We denote by  $\Delta:X\to X^{I}$ the   diagonal morphism.  
 We fix a  geometric point $\ov{\eta^{I}}$ over  the generic point $\eta^{I}$ of $X^{I}$ and a specialization arrow   
 $\on{\mf{sp}}: \ov{\eta^{I}}\to \Delta(\ov \eta)$. The  role of  $\on{\mf{sp}}$  is to make  the fiber functor at  $\ov{\eta^{I}}$ more canonical, and in particular  compatible with the  coalescence of legs (the last claim is clear when   $\on{\mf{sp}}^{*}$ is an   isomorphism and in practice we will be essentially  in this  situation). So $ \ov{\eta^{I}}$ and $\on{\mf{sp}}$ go together and below the statements involving $ \ov{\eta^{I}}$  depend on the choice of  $\on{\mf{sp}}$.

 A   fondamental result of Drinfeld (theorem 2.1 of \cite{drinfeld78} and proposition 6.1 of \cite{drinfeld-compact}) is recalled   in the  
 following lemma (see   chapter 8 of \cite{coh} for a reminder of the proof and other  references).
  We will always denote    the  $\mc O_{E}$-modules and  the $\mc O_{E}$-sheaves   by gothic letters. 
   
 \begin{lem} \label{lem-Dr-intro} (Drinfeld)  If 
 $\mf E$ is a  smooth   constructible  $\mc O_{E}$-sheaf  
 over a  dense  open subscheme   of  
 $(X\sm N)^{I}$, equipped with an  action of the   partial Frobenius morphisms, i.e.   with isomorphisms
 $$F_{\{i\}}:\restr{\Frob_{\{i\}}^{*}(\mf E)}{\eta^{I}}\to \restr{\mf E}{\eta^{I}}$$
 commuting to each other and whose  composition is  the natural isomorphism  
    $\Frob^{*}(\mf E)\isom \mf E$, 
then it extends to a smooth  sheaf  $\wt {\mf E}$ over  $U^{I}$, where  $U$ is  a small enough open dense subscheme of  $X\sm N$, and the fiber   $\restr{\wt {\mf E}}{\Delta(\ov \eta)}$  is equipped with an  action of   $\pi_{1}(U,\ov\eta)^{I}$.

Moreover, if we fix   $U$,  the  functor
 $\mf E\mapsto \restr{\mf E}{\Delta(\ov \eta)}$ 
 provides an  equivalence 
 \begin{itemize}
 \item from  the  category of smooth constructible    $\mc O_{E}$-sheaves
   over   $U^{I}$ equipped with an  action of the  partial Frobenius morphisms
   \item  to  the  category of continuous  
 representations of   $\pi_{1}(U,\ov\eta)^{I}$ on   $\mc O_{E}$-modules of   finite type, 
 \end{itemize}
  in a way  compatible with  coalescence (i.e. inverse   image by morphisms $U^{J}\to U^{I}$ associated to any map $I\to J$). 
 \end{lem}
       
       \begin{rem} \label{rem-lem-Drinfeld} In the situation of the beginning of the previous lemma, 
       $$\on{\mf{sp}}^{*}: 
       \restr{\wt {\mf E}}{\Delta(\ov \eta)}\to \restr{\wt {\mf E}}{\ov{\eta^{I}}}=\restr{\mf E}{\ov{\eta^{I}}}$$  is an isomorphism, hence $\restr{\mf E}{\ov{\eta^{I}}}$ is also equipped with an  action of $\pi_{1}(U,\ov\eta)^{I}$. 
       \end{rem}
          \dem   The idea of the proof of  \lemref{lem-Dr-intro} is that it is enough to handle the case where      $\wt {\mf E}=\boxtimes _{i\in I} \mf E_{i}$ with  $ \mf E_{i}$ smooth over  $U$, and in this case  the action of  $\pi_{1}(U,\ov\eta)^{I}$
    on  $ \restr{\wt {\mf E}}{\Delta(\ov \eta)}=  \otimes _{i\in I} \restr{ \mf E_{i}}{ \ov \eta}$ is obvious.  \cqfd

 We cannot apply directly this lemma, because the action of the   partial Frobenius morphisms increases $\mu$, and on the other hand  the inductive limit   $\varinjlim_{\mu}\mc H_{N,I,W}^{\leq\mu}$ is not  constructible  (because its fibers are of   infinite dimension). 
Nevertheless we will be able to apply Drinfeld's lemma to the  ``Hecke-finite'' part, in the following sense.

 \begin{defi}
 Let $\ov x$ be a  geometric point of  $(X\sm N)^{I}$. 
An  element of  $\varinjlim _{\mu}\restr{\mc H _{N, I, W}^{\leq\mu}}{\ov{x}}$ is said to be   Hecke-finite if it belongs to a   $\mc O_{E}$-submodule  of finite type  of  $\varinjlim _{\mu}\restr{\mc H _{N, I, W}^{\leq\mu}}{\ov{x}}$  which is 
stable by  $T(f)$ for every  $f\in C_{c}(K_{N}\backslash G(\mb A)/K_{N},\mc O_{E})$. 
\end{defi}
We denote by   $\Big( \varinjlim _{\mu}\restr{\mc H _{N, I, W}^{\leq\mu}}{\ov{x}}\Big)^{\mr{Hf}}$ the set of all the  Hecke-finite elements. 
It is  a $E$-vector subspace of   $ \varinjlim _{\mu}\restr{\mc H _{N, I, W}^{\leq\mu}}{\ov{x}}$ and it is stable 
 by  $\pi_{1}(x,\ov{x})$ and  $C_{c}(K_{N}\backslash G(\mb A)/K_{N},E)$. 

\begin{rem}
The  definition above will be used    with  $\ov x$ equal to 
$\Delta(\ov{\eta})$ or  $\ov{\eta^{I}}$. In this case the action of the  Hecke operators $T(f)$ on $\varinjlim _{\mu}\restr{\mc H _{N, I, W}^{\leq\mu}}{\ov{x}}$ is obvious (and does not require their extension to morphisms of sheaves over $(X\sm N)^{I}$  which will be obtained  after   \propref{Eichler-Shimura-intro}). 
\end{rem}
 
 We have  the specialization homomorphism
 \begin{gather}\label{sp*-sans-Hf-intro}\on{\mf{sp}}^{*}: 
 \varinjlim _{\mu}\restr{\mc H _{N, I, W}^{\leq\mu}}{\Delta(\ov{\eta})} \to
  \varinjlim _{\mu}\restr{\mc H _{N, I, W}^{\leq\mu}}{\ov{\eta^{I}}}\end{gather} 
 where both sides are considered as   $E$-vector spaces
 (inductive limits  of $E$-vector spaces of finite  dimension). 

We admit now  two results,  which will be justified  in  
sections \ref{subsection-crea-annihil-intro}, \ref{section-sous-sheaves-constr} and \ref{para-homom-spe-Hf}.

 \noindent{\bf First result that we admit for the moment} (\lemref{lem-Hf-union-stab} and \propref{cor-action-Hf-intro}). The space 
   $\Big( \varinjlim _{\mu}\restr{\mc H _{N, I, W}^{\leq\mu}}{\ov{\eta^{I}}}\Big)^{\mr{Hf}}$ is a union of $\mc O_{E}$-submodules of finite type stable  by  the partial Frobenius morphisms and by   Drinfeld's lemma it is endowed with an  action  of $\on{Gal}(\ov F/F)^{I}$, more precisely  $\Big( \varinjlim _{\mu}\restr{\mc H _{N, I, W}^{\leq\mu}}{\ov{\eta^{I}}}\Big)^{\mr{Hf}}$ is an inductive limit of    finite  dimensional  continuous representations of  $\on{Gal}(\ov F/F)^{I}$. 
   
  \noindent{\bf Second result that we admit for the moment}  (\corref{bijectivite-Hecke-fini}). The restriction of the homomorphism $\on{\mf{sp}}^{*}$ of \eqref{sp*-sans-Hf-intro} to the Hecke-finite parts is  an  isomorphism 
   \begin{gather}\label{isom-sp*-avec-Hf-intro} 
  \Big( \varinjlim _{\mu}\restr{\mc H _{N, I, W}^{\leq\mu}}{\Delta(\ov{\eta})}\Big)^{\mr{Hf}} \isor{\on{\mf{sp}}^{*}}
 \Big(  \varinjlim _{\mu}\restr{\mc H _{N, I, W}^{\leq\mu}}{\ov{\eta^{I}}}\Big)^{\mr{Hf}}. \end{gather}

 Thanks to these two results we can  now define the $E$-vector spaces $H_{I,W}$  (we omit  $N$ in the notation $ H_{I,W}$ to reduce the size of the diagrams in the next section). 
    
   \begin{defi}
    We define  $H_{I,W}$ as the LHS of \eqref{isom-sp*-avec-Hf-intro}. 
  \end{defi}
 
 The action of $\on{Gal}(\ov F/F)^{I}=\pi_{1}(\eta,\ov\eta)^{I}$ on 
  $H_{I,W}$ does not depend on the  choice of  $\ov{\eta^{I}}$ and $\on{\mf{sp}}$. 
 Indeed  by the first admitted result there exist  
  \begin{itemize}
 \item an increasing union 
 (indexed by $\lambda\in \N$) 
 of constructible $\mc O_{E}$-subsheaves 
 $\mf F_{\lambda}\subset \varinjlim _{\mu}\restr{\mc H _{N, I, W}^{\leq\mu}}{\eta^{I}}$ stable  by the partial Frobenius morphisms
 (to which we can apply Drinfeld's lemma) 
 \item a decreasing sequence of open subschemes $U_{\lambda}\subset X\sm N$ such that 
 $\mf F_{\lambda}$ extends to a smooth sheaf  over $(U_{\lambda})^{I}$
 \end{itemize}
 in such a way that  $\bigcup _{\lambda\in \N}  \restr{\mf F_{\lambda}}{\ov{\eta^{I}}} 
 =
  \Big(  \varinjlim _{\mu}\restr{\mc H _{N, I, W}^{\leq\mu}}{\ov{\eta^{I}}}\Big)^{\mr{Hf}}$. 
   Then the second admitted result   
  implies that   the natural morphism  
  \begin{gather}\label{mor-HIW-Flambda}
  H_{I,W}= \Big( \varinjlim _{\mu}\restr{\mc H _{N, I, W}^{\leq\mu}}{\Delta(\ov{\eta})}\Big)^{\mr{Hf}} \to \bigcup _{\lambda\in \N}  \restr{\mf F_{\lambda}}{\Delta(\ov{\eta})} \end{gather}
(which comes from the smoothness of $\mf F_{\lambda}$ over $(U_{\lambda})^{I}\ni 
\Delta(\ov{\eta})$)
  is an isomorphism. But  the action of $\on{Gal}(\ov F/F)^{I}$   on the RHS of \eqref{mor-HIW-Flambda},  which is given by Drinfeld's lemma, does not depend on the choice of  $\ov{\eta^{I}}$ and $\mf{sp}$, and therefore  the action  of $\on{Gal}(\ov F/F)^{I}$ on the LHS does not depend on it either. 
  
 \begin{rem}
In this article we only prove that   $H_{I,W}$ is an inductive limit of $E$-vector spaces of finite  dimension equipped with  continuous representations  of $\on{Gal}(\ov F/F)^{I}$. In fact 
   Cong Xue proved in  \cite{these-cong} that 
 $H_{I,W}$ is of finite dimension. The proof 
is difficult and written only when $G$ is  split.  
 \end{rem}

   For every map  $\zeta:I\to J$, the coalescence isomorphism  
   \eqref{intro-isom-coalescence} obviously  respects  the  Hecke-finite parts and therefore it induces an  isomorphism 
    \begin{gather} \label{isom-chi-zeta-intro}
    H_{I,W}= \Big( \varinjlim _{\mu}\restr{\mc H _{N, I, W}^{\leq\mu}}{\Delta(\ov{\eta})}\Big)^{\mr{Hf}}
    \isom 
    \Big( \varinjlim _{\mu}\restr{\mc H _{N, J,W^{\zeta}}^{\leq\mu}}{\Delta(\ov{\eta})}\Big)^{\mr{Hf}}=H_{J,W^{\zeta}}
    \end{gather}
    where $\Delta$ denotes  the diagonal morphism $X\to X^{I}$ or  $X\to X^{J}$. 
    
    \begin{defi} We define  the isomorphism
    $$ \chi_{\zeta}:    H_{I,W} \isom H_{J,W^{\zeta}}$$ occuring  in  b) of the \propref{prop-a-b-c} to be  \eqref{isom-chi-zeta-intro}. 
    \end{defi}
    
    The isomorphism  $\chi_{\zeta}$ is $\on{Gal}(\ov F/F)^{J}$-equivariant,  where $\on{Gal}(\ov F/F)^{J}$ acts on the LHS by the   diagonal morphism  
\begin{gather}\label{morp-diag-Gal-intro}\on{Gal}(\ov F/F)^{J}\to \on{Gal}(\ov F/F)^{I},  \ (\gamma_{j})_{j\in J}\mapsto (\gamma_{\zeta(i)})_{i\in I}. 
\end{gather}
Indeed, if $\Delta_{\zeta}:X^{J}\to X^{I}$ is the diagonal morphism 
    \eqref{morph-giad-X-intro} and if the sequence   $(\mf F_{\lambda})_{\lambda\in \N}$ is as above relatively to $I$ and $W$, then the sequence  $(\Delta_{\zeta}^{*}(\mf F_{\lambda}))_{\lambda\in \N}$ satisfies the same properties relatively to  $J$ and $W^{\zeta}$, and thus 
       $$\chi_{\zeta}: H_{I,W}=\bigcup _{\lambda\in \N}  \restr{\mf F_{\lambda}}{\Delta(\ov{\eta})} =\bigcup _{\lambda\in \N}  \restr{\Delta_{\zeta}^{*}(\mf F_{\lambda})}{\Delta(\ov{\eta})}=H_{J,W^{\zeta}}
       $$ is $\on{Gal}(\ov F/F)^{J}$-equivariant. 
    
\noindent {\bf Proof of \propref{prop-a-b-c}.  }    The properties a) and b) were already explained.     
Property c) comes from the fact that   the Hecke-finite  part of 
\eqref{I-vide-intro} consists exactly of cuspidal  automorphic forms, i.e.  $$\Big(C_{c} (G(F)\backslash G(\mb A)/K_N \Xi,E)\Big)^{\mr{Hf}}=C_{c}^{\rm{cusp}}(G(F)\backslash G(\mb A)/K_N \Xi,E).$$

\noindent {\bf Proof of $\supset$. } Any   cuspidal function is Hecke-finite  
 because  the $\mc O_{E}$-module $ C_{c}^{\rm{cusp}}(G(F)\backslash G(\mb A)/K_N \Xi,\mc O_{E}) $
  is   of finite type and is stable by all the Hecke operators $T(f)$ for $f\in C_{c}(K_{N}\backslash G(\mb A)/K_{N},\mc O_{E})$. 
 
 \noindent {\bf Proof of $\subset$. } We assume by contradiction that a   Hecke-finite  function  $f$  is not   cuspidal. Then there exists a parabolic subgroup $P\subsetneq G$, of Levi quotient $M$ and unipotent radical  $U$,  such that the constant term $f_{P}: g\mapsto \int_{U(F)\backslash U(\mb A)}f(ug)$ is non zero.   Let $v$ be a  place of $X\sm N$. 
       Since the ring of (finite dimensional) representations of   $\wh M$ is a module of finite type over  the ring of  representations of $\wh G$, $f_{P}$ is also Hecke-finite w.r.t. the   Hecke operators   for $M$ at $v$. 
  These operators admit as particular cases the translations by the elements of $Z_{M}(F_{v})$. 
  We have a degree map (relatively to  
       $M/Z$), from $U(\mb A)M(F)\backslash G(\mb A)/K_{N}\Xi$ to a free  $\Z$-module of finite type, on which $Z_{M}(F_{v})$ acts by non trivial  translations. 
      But the support of $f_{P}$ is included in a translation of a cone in this 
 free  $\Z$-module of finite type, 
 and this contradicts the fact that 
 $f_{P}$ belongs to a finite dimensional vector space stable by      $Z_{M}(F_{v})$. We refer to    proposition 8.23 of \cite{coh} for the details of the  proof.   \cqfd

    \section{Proof  of  \thmref{intro-thm-ppal}  with the help of   \propref{prop-a-b-c}} 
    \label{intro-idee-heurist}
    Here is an overview of the idea. 
           Thanks to  \eqref{egalite-vide-0-cusp} (which we had decued from   \propref{prop-a-b-c}),    we have  
    $$H_{\{0\},\mbf 1}=C_{c}^{\mr{cusp}}(G(F)\backslash G(\mb A)/K_{N}\Xi,E).$$
 To get the decomposition \eqref{intro1-dec-canonical}  it is thus equivalent to  construct (increasing  $E$ if necessary) 
   a   canonical decomposition
     \begin{gather}\label{intro3-dec-canonical}
  H_{\{0\},\mbf 1}=\bigoplus_{\sigma}
 \mf H_{\sigma}.\end{gather}
          This decomposition will be  obtained by   spectral decomposition  of a commutative family of endomorphisms of $H_{\{0\},\mbf 1}$, called excursion operators, that we will now    study with the help of  properties a) and  b) of   \propref{prop-a-b-c}. First we recall their construction, already given in  \defiref{defi-constr-excursion-intro}.

      Let $I$   be a finite  set and  $W$ be a  $E$-linear representation  of  
    $(\wh G)^{I}$. 
    We denote by  $\zeta_{I}:I\to \{0\}$ the obvious map, so that  
    $ W^{\zeta_{I}}$ is just  $W$ equipped with  the diagonal action   of $\wh G$. Let 
$x: \mbf 1\to W^{\zeta_{I}}$ and  $\xi :  W^{\zeta_{I}}\to  \mbf 1$
be  morphisms of  representations of  $\wh G$ (in other words  $x\in W$ and  $\xi\in W^{*}$  are invariant  under  the diagonal action  of  $\wh G$).  Let  $(\gamma_{i})_{i\in I}\in \on{Gal}(\ov F/F)^{I}$. 
 
 \begin{defi}We define  the operator  
  \begin{gather*}S_{I,W,x,\xi,(\gamma_{i})_{i\in I}}\in 
  \on{End}(H_{\{0\},\mbf 1})  \end{gather*}
    as the composition 
  \begin{gather}\label{excursion-def-intro}
  H_{\{0\},\mbf  1}\xrightarrow{\mc H(x)}
 H_{\{0\},W^{\zeta_{I}}}\isor{\chi_{\zeta_{I}}^{-1}} 
  H_{I,W}
  \xrightarrow{(\gamma_{i})_{i\in I}}
  H_{I,W} \isor{\chi_{\zeta_{I}}} H_{\{0\},W^{\zeta_{I}}}  
  \xrightarrow{\mc H(\xi)} 
  H_{\{0\},\mbf  1}.
    \end{gather}\end{defi}
    
   This operator will be called an ``excursion operator''. Paraphrasing \eqref{excursion-def-intro} it is the composition of 
    \begin{itemize}
    \item a    creation operator associated to $x$, whose effect is to create legs at the same (generic) point  of the curve, 
    \item a  Galois action, 
     which moves the legs  over the curve independently of each other, and then puts then back 
    at the same (generic) point  of the curve, 
       \item an annihilation operator  associated to $\xi$,  which annihilates the legs. 
         \end{itemize}

   The following lemma will be deduced from  properties a) and  b) of  \propref{prop-a-b-c}. 
   \begin{lem} The  excursion operators $S_{I,W,x,\xi,(\gamma_{i})_{i\in I}}$ satisfy 
    the following properties : 
     \begin{gather}
     \label{SIW-p0}      S_{I,W,x,{}^{t} u(\xi'),(\gamma_i)_{i\in I}}=S_{I,W',u(x),\xi',(\gamma_i)_{i\in I}} 
 \end{gather}
where $u:W\to W'$ is a   $(\wh G)^{I}$-equivariant  morphism and  $x\in W$ and $\xi'\in (W')^{*}$ are $\wh G$-invariant,  
  \begin{gather}   \label{SIW-p1}
     S_{J,W^{\zeta},x,\xi,(\gamma_j)_{j\in J}}=S_{I,W,x,\xi,(\gamma_{\zeta(i)})_{i\in I}},
     \\
      \label{SIW-p2}
 S_{I_{1}\cup I_{2},W_{1}\boxtimes W_{2},x_{1}\boxtimes x_{2},\xi_{1}\boxtimes \xi_{2},(\gamma^{1}_i)_{i\in I_{1}}\times (\gamma^{2}_i)_{i\in I_{2}}}= S_{I_{1},W_{1},x_{1},\xi_{1},(\gamma^{1}_i)_{i\in I_{1}}}\circ 
S_{I_{2},W_{2},x_{2},\xi_{2},(\gamma^{2}_i)_{i\in I_{2}}}, 
\\ 
\label{SIW-p3} 
S_{I,W,x,\xi,(\gamma_i(\gamma'_i )^{-1}\gamma''_i)_{i\in I}}=
    S_{I\cup I \cup I,W\boxtimes W^{*}\boxtimes W,\delta_{W} \boxtimes x,
    \xi \boxtimes \on{ev}_{W},
    (\gamma_i)_{i\in I} \times (\gamma'_i)_{i\in I} \times (\gamma''_i)_{i\in I}
    }
             \end{gather}
  where most of the   notations are obvious, 
    $I_{1}\cup I_{2}$ and  $I\cup I\cup I$ denote disjoint unions, and  $\delta_{W}: \mbf 1\to W\otimes W^{*}$ and $\on{ev}_{W}:W^{*} \otimes W \to \mbf 1$ are the natural morphisms.      
  \end{lem}
   
     \noindent {\bf Proof of  \eqref{SIW-p0}.}  
    We set  $x'=u(x)\in W'$ and $\xi={}^{t}u(\xi')\in W^{*}$. 
The diagram 
$$   \xymatrix{
     &  H_{\{0\},(W')^{\zeta_{I}}}  
    \ar[r]^{\chi_{\zeta_{I}}^{-1}}&  H_{I,W'} \ar[r]^{(\gamma_{i})_{i\in I}} 
     &
    H_{I,W'} \ar[r]^{\chi_{\zeta_{I}}} &  
    H_{\{0\},(W')^{\zeta_{I}}}  \ar[dr]^{\mc H(\xi')} 
     \\
  H_{\{0\},\mbf  1}  \ar[r]_{\mc H(x)}    \ar[ur]^{\mc H(x')}  &    
    H_{\{0\},W^{\zeta_{I}}}   \ar[u]_{\mc H(u)}
    \ar[r]_{\chi_{\zeta_{I}}^{-1}} & 
    H_{I,W}  \ar[u]_{\mc H(u)} \ar[r]_{(\gamma_{i})_{i\in I}} & 
   H_{I,W}  \ar[u]^{\mc H(u)} \ar[r]_{\chi_{\zeta_{I}}} &
     H_{\{0\},W^{\zeta_{I}}} \ar[u]^{\mc H(u)} \ar[r]_{\mc H(\xi)} & 
     H_{\{0\},\mbf  1}
    }$$           
is commutative. But the lower line is equal to   
   $S_{I,W,x,\xi,(\gamma_{i})_{i\in I}}
$ and the upper line  is equal to  
   $S_{I,W',x',\xi',(\gamma_{i})_{i\in I}}
$. \cqfd
       
   \noindent {\bf Proof  of  \eqref{SIW-p1}.}  
The diagram 
$$   \xymatrix{
     & &  H_{J,W^{\zeta}} \ar[rr]^{(\gamma_{j})_{j\in J}}  &&
    H_{J,W^{\zeta}} \ar[dr]^{\chi_{\zeta_{J}}}
     \\
  H_{\{0\},\mbf  1}  \ar[r]_{\mc H(x)}    &    
    H_{\{0\},W^{\zeta_{I}}}   \ar[ur]^{\chi_{\zeta_{J}}^{-1}} 
    \ar[r]_{\chi_{\zeta_{I}}^{-1}} & 
    H_{I,W}  \ar[u]_{\chi_{\zeta}} \ar[rr]_{(\gamma_{\zeta(i)})_{i\in I}} && 
   H_{I,W}  \ar[u]^{\chi_{\zeta}} \ar[r]_{\chi_{\zeta_{I}}} &
     H_{\{0\},W^{\zeta_{I}}}  \ar[r]_{\mc H(\xi)} & 
     H_{\{0\},\mbf  1}
    }$$           
is commutative. But the lower line is equal to   
   $S_{I,W,x,\xi,(\gamma_{\zeta(i)})_{i\in I}}
$ and the upper line  is equal to  
   $    S_{J,W^{\zeta},x,\xi,(\gamma_j)_{j\in J}}$.  \cqfd
       
   \noindent {\bf Proof  of  \eqref{SIW-p2}.}  
   The obvious map 
   $\{0\}\cup \{0\} \to \{0\}$ 
  gives an   isomorphism $ H_{\{0\}\cup \{0\},\mbf  1}\simeq  H_{\{0\},\mbf  1}$. If we denote by   $\zeta_{1}:I_{1}\to \{0\}$ and  $\zeta_{2}:I_{2}\to \{0\}$ the obvious maps,  the LHS of  \eqref{SIW-p2} is equal to the composition 
       \begin{gather*}
 H_{\{0\}\cup \{0\},\mbf  1}
 \xrightarrow{\mc H(x_{1}\boxtimes x_{2})}
 H_{\{0\}\cup \{0\},W_{1}^{\zeta_{1}}\boxtimes W_{2}^{\zeta_{2}}}
 \isor{\chi_{\zeta_{1}\times \zeta_{2}}^{-1}} 
  H_{I_{1}\cup I_{2},W_{1}\boxtimes W_{2}}\\
  \xrightarrow{(\gamma^{1}_{i})_{i\in I_{1}}\times (\gamma^{2}_{i})_{i\in I_{2}}}
  H_{I_{1}\cup I_{2},W_{1}\boxtimes W_{2}} \isor{\chi_{\zeta_{1}\times \zeta_{2}}} 
  H_{\{0\}\cup \{0\},W_{1}^{\zeta_{1}}\boxtimes W_{2}^{\zeta_{2}}}  
  \xrightarrow{\mc H(\xi_{1}\boxtimes \xi_{2})} 
  H_{\{0\}\cup \{0\},\mbf  1}. 
   \end{gather*}
Putting together  $x_{1},\chi_{\zeta_{1}}^{-1},(\gamma^{1}_{i})_{i\in I_{1}} ,\chi_{\zeta_{1}},\xi_{1}$ on one side   and 
$x_{2},\chi_{\zeta_{2}}^{-1},(\gamma^{2}_{i})_{i\in I_{2}} ,\chi_{\zeta_{2}},\xi_{2}$ on the other side we get  the RHS. 
We are allowed to do this because in the following diagram   (where we write 
$\gamma^{1}=(\gamma^{1}_{i})_{i\in I_{1}}$ and $\gamma^{2}=(\gamma^{2}_{i})_{i\in I_{2}}$) 
all squares and triangles commute. 
\cqfd

   {  \resizebox{14cm}{!}{ $$\!\!\!\!\!\!\!\!\!\!\!\!\!\!\!    \xymatrix{
     {  H_{\{0\}\cup \{0\},\mbf 1} }
           \ar[d]_{{\mc H(x_{1}\boxtimes 1)} }
           \ar[rd]^-{{\mc H(x_1\boxtimes x_2)}}
             \\
      {    H_{\{0\}\cup \{0\},W_{1}^{\zeta_{1}}\boxtimes \mbf 1}}
           \ar[d]_{{\chi_{\zeta_1\times \Id}^{-1}}}
           \ar[r]^<<<{{\mc H(\Id\boxtimes x_{2} )} }
       &  { H_{\{0\}\cup \{0\},W_{1}^{\zeta_{1}}\boxtimes W_{2}^{\zeta_{2}}}}
        \ar[d]_{{\chi_{\zeta_1\times \Id}^{-1}} }
        \ar[rd]^-{{\chi_{\zeta_1\times \zeta_{2}}^{-1}}}
            \\
    {     H_{I_{1}\cup \{0\},W_{1}\boxtimes \mbf 1} }
         \ar[d]_{{\gamma^{1}\times 1}}
         \ar[r]^-{{ \mc H(\Id\boxtimes x_{2} )} }
       &   { H_{I_{1}\cup \{0\},W_{1} \boxtimes W_{2}^{\zeta_{2}}}}
        \ar[d]_{{\gamma^{1}\times 1}}
        \ar[r]^-{{\chi_{\Id\times \zeta_2}^{-1}} }
       &{  H_{I_{1}\cup I_{2},W_{1} \boxtimes W_{2} }  }
         \ar[d]_{{\gamma^{1}\times 1} }
         \ar[rd]^-{ {\gamma^{1}\times \gamma^{2}}}
   \\
   { H_{I_{1}\cup \{0\},W_{1}\boxtimes \mbf 1} }
    \ar[d]_{{\chi_{\zeta_1\times \Id}}}
  \ar[r]^-{{\mc H(\Id\boxtimes x_{2} )} }&  
  { H_{I_{1}\cup \{0\},W_{1} \boxtimes W_{2}^{\zeta_{2}}} }
    \ar[d]_{{\chi_{\zeta_1\times \Id}}}
  \ar[r]^-{{\chi_{\Id\times \zeta_2}^{-1}} }
  & 
 {   H_{I_{1}\cup I_{2},W_{1} \boxtimes W_{2}} }
    \ar[d]_{{\chi_{\zeta_1\times \Id}}}
   \ar[r]^-{{ 1\times  \gamma^{2}}  }
    &{  H_{I_{1}\cup I_{2},W_{1} \boxtimes W_{2}}   }
     \ar[d]_{{\chi_{\zeta_1\times \Id}}}
   \ar[rd]^-{{\chi_{\zeta_1\times \zeta_{2}}}       }
 \\
  {  H_{\{0\}\cup \{0\},W_{1}^{\zeta_{1}}\boxtimes \mbf 1} }
  \ar[d]_{{\mc H(\xi_{1}\boxtimes 1)}  }
  \ar[r]^<<<{{\mc H(\Id\boxtimes x_{2} )} }
  &  { H_{\{0\}\cup \{0\},W_{1}^{\zeta_{1}} \boxtimes W_{2}^{\zeta_{2}}}}
     \ar[d]_{{\mc H(\xi_{1}\boxtimes \Id)} } 
  \ar[r]^-{{\chi_{\Id\times \zeta_2} ^{-1}}}
  & { H_{\{0\}\cup I_{2},W_{1}^{\zeta_{1}} \boxtimes W_{2}}}
    \ar[r]^-{{ 1\times  \gamma^{2}}  }
     \ar[d]_{{\mc H(\xi_{1}\boxtimes \Id)} } &
 {  H_{\{0\}\cup I_{2},W_{1}^{\zeta_{1}} \boxtimes W_{2}}}
      \ar[d]_{{\mc H(\xi_{1}\boxtimes \Id)}  }
       \ar[r]^-{{\chi_{\Id\times \zeta_{2}}}   }
  & {  H_{\{0\}\cup \{0\},W_{1}^{\zeta_{1}} \boxtimes W_{2}^{\zeta_{2}} }}
  \ar[d]_{{\mc H(\xi_{1}\boxtimes \Id)}  }
 \ar[rd]^-{{\mc H(\xi_{1} \boxtimes \xi_2)} }
 \\           
 { H_{\{0\}\cup \{0\},\mbf 1} }
  \ar[r]^-{{ \mc H(1\boxtimes x_{2} )}} &   
 { H_{\{0\}\cup \{0\},\mbf 1  \boxtimes W_{2}^{\zeta_{2}}} }
  \ar[r]^-{{ \chi_{\Id\times \zeta_2} ^{-1}}}&  
{  H_{\{0\}\cup I_{2},\mbf 1  \boxtimes W_{2}} }
   \ar[r]^-{{ 1\times  \gamma^{2}}   } &  
 { H_{\{0\}\cup I_{2},\mbf 1  \boxtimes W_{2}} }
 \ar[r]^-{{\chi_{\Id\times \zeta_{2}}}   }
&{  H_{\{0\}\cup \{0\},\mbf 1  \boxtimes W_{2}^{\zeta_{2}}} }
\ar[r]^-{{\mc H(1\boxtimes \xi_{2} )}} & 
{ H_{\{0\}\cup \{0\},\mbf 1}} }$$}
  }
          
  \vskip1mm        
  \noindent{\bf Proof  of   \eqref{SIW-p3}.}  
   For every  $(g_i)_{i\in I}\in (\wh G)^{I}$, 
\begin{itemize}
\item   $ \xi \boxtimes \on{ev}_{W}$     is  invariant by $(1)_{i\in I}\times (g_i)_{i\in I} \times  (g_i)_{i\in I} $
\item $\delta_{W} \boxtimes x$   is  invariant by      $(g_i)_{i\in I} \times (g_i)_{i\in I} \times  (1)_{i\in I} $.
    \end{itemize}    
Therefore for every  $(\alpha_{i})_{i\in I}$ and $(\beta_{i})_{i\in I}$ in $\on{Gal}(\ov F/F)^{I}$, the RHS of \eqref{SIW-p3} is equal to \begin{gather}\label{SIW-p3-alpha-beta}
 S_{I\cup I \cup I,W\boxtimes W^{*}\boxtimes W,\delta_{W} \boxtimes x,
    \xi \boxtimes \on{ev}_{W},
    (\gamma_i\beta_{i})_{i\in I} \times (\alpha_{i}\gamma'_i\beta_{i})_{i\in I} \times (\alpha_{i}\gamma''_i)_{i\in I}
   . } \end{gather}
 To prove it in a formal way we factorize  the RHS of \eqref{SIW-p3} through  
  $$H_{I, \mbf 1} \xrightarrow{ \mc H(\delta_{W})}
     H_{ I,(W\boxtimes W^{*})^{\zeta}} \text{ \  and  \   } H_{ I,(W^{*}\boxtimes W)^{\zeta}} \xrightarrow{ \mc H(\on{ev}_{W})}
    H_{I, \mbf 1},$$ where  $\zeta:I\cup I\to I$ is the obvious map,   
and we use the fact that   $\on{Gal}(\ov F/F)^{I}$ acts trivially on   $ H_{I,\mbf 1}\simeq  H_{\emptyset,\mbf 1}$. We take  $\alpha_{i}=\gamma_i (\gamma'_i)^{-1}$ and 
$\beta_{i}=(\gamma'_i)^{-1}\gamma''_i $. Then   \eqref{SIW-p3-alpha-beta} is equal to  
           \begin{gather}\label{eq-lem-b-2-ggg-intro}
                 S_{I\cup I \cup I,W\boxtimes W^{*}\boxtimes W,\delta_{W} \boxtimes x,\xi \boxtimes \on{ev}_{W},
                 (\gamma_i(\gamma'_i )^{-1}\gamma''_i)_{i\in I} \times 
                 (\gamma_i(\gamma'_i )^{-1}\gamma''_i)_{i\in I} \times
                  (\gamma_i(\gamma'_i )^{-1}\gamma''_i)_{i\in I}
                  }.                   \end{gather}
Applying   \eqref{SIW-p1} to  the obvious map  $\zeta:I\cup I \cup I\to I$, 
we see that \eqref{eq-lem-b-2-ggg-intro} is equal to 
  \begin{gather}\label{eq-lem-b-2-ggg2-intro}S_{I,W\otimes W^{*}\otimes W,\delta_{W} \otimes x,\xi \otimes \on{ev}_{W},(\gamma_i(\gamma'_i )^{-1}\gamma''_i)_{i\in I}}.\end{gather}
Lastly one shows that 
 \eqref{eq-lem-b-2-ggg2-intro} is equal  to the LHS of \eqref{SIW-p3}
by applying   \eqref{SIW-p0} 
to   the $(\wh G)^{I}$-linear injection   $$u:W=\mbf  1\otimes W
    \xrightarrow{  \delta_{W}    \otimes\Id_{W}} W\otimes W^{*}\otimes W, $$
         which satisfies   $\delta_{W} \otimes x =u( x)$ and  
                  ${}^{t}u(\xi\otimes\on{ev}_{W}) =\xi$, since 
                   \begin{gather}\label{zorro} \text{the composition   } W
    \xrightarrow{  \delta_{W}    \otimes\Id_{W}} W\otimes W^{*}\otimes W
     \xrightarrow{ \Id_{W}\otimes\on{ev}_{W}} W 
                \text{    is equal to     } \Id_{W} \end{gather}
                (this is an easy  lemma in the category of vector spaces, sometimes called Zorro lemma, and   it is  also one of the axioms of the tannakian categories).       
                       \cqfd

      We denote by  
$\mc B \subset \on{End}_{C_{c}(K_{N}\backslash G(\mb A)/K_{N},E)}( H_{\{0\},\mbf 1}) $ 
the subalgebra generated by all  the  excursion operators $S_{I,W,x,\xi,(\gamma_{i})_{i\in I}}$.
By  \eqref{SIW-p2},  $\mc B$ is commutative. 
   
   In the rest of this  section we consider  $\wh G$ as a group scheme  defined  over $E$. 
     
 \begin{observation} \label{observation-f-W-x-xi} The functions \begin{gather}\label{intro-def-f}
   f: (g_{i})_{i\in I}\mapsto \s{\xi, (g_{i})_{i\in I}\cdot x}\end{gather} that  we get by making   $W$, 
 $x$, and  $\xi $ vary, are exactly the regular functions   on  the coarse quotient     of  $(\wh G)^{I}$ by   left and right translations by diagonal $\wh G$, that we denote by  $\wh G\backslash (\wh G)^{I}/\wh G$. 
 \end{observation}
 
 \begin{lem}\label{f-W,x,xi}
The operator $S_{I,W,x,\xi,(\gamma_{i})_{i\in I}}$ depends only on   $I$, $f$, and $(\gamma_{i})_{i\in I}$, where $f$ is given  by \eqref{intro-def-f}. 
\end{lem}
\dem Let $W,x,\xi$ be as above and let  $f\in \mc O(\wh G\backslash (\wh G)^{I}/\wh G)$ be given  by \eqref{intro-def-f}. 
 We denote by  $W_{f}$ the  finite  dimensional $E$-vector subspace    of 
 $\mc O((\wh G)^{I}/\wh G)$ generated by the left  translates  of $f$ by $(\wh G)^{I}$. We set  $x_{f}=f\in W_{f}$ and we denote by  $\xi_{f}$ the linear form   on  $W_{f}$ given by the evaluation at $1\in (\wh G)^{I}/\wh G$. 
Then  $W_{f}$ is a  subquotient of $W$: if $W_{x}$ is the $E$-linear 
 $(\wh G)^{I}$-subrepresentation   of $W$ generated by $x$, 
 $W_{f}$ is the  quotient of $W_{x}$ by the biggest  $E$-linear $(\wh G)^{I}$-subrepresentation on which $\xi$ vanishes. 
 Thus we  have  the   diagrams $$W \overset{\alpha}{\hookleftarrow} W_{x} \overset{\beta}{\twoheadrightarrow} W_{f}, 
 \ \ x \overset{\alpha}{\longleftarrow\!\shortmid} x \overset{\beta}{\shortmid\!\longrightarrow} x_{f}, 
 \ \ \xi \overset{{}^{t}\alpha}{\shortmid\!\longrightarrow}  \restr{\xi}{W_{x}} \overset{{}^{t}\beta}{\longleftarrow\!\shortmid} \xi_{f}$$
 of  $(\wh G)^{I}$-representations, of 
    $\wh G$-invariant vectors  and of $\wh G$-invariant linear forms. 
 Applying \eqref{SIW-p0} 
  to $u=\alpha$ and $u=\beta$, we get 
      \begin{gather}\nonumber 
  S_{I,W,x,\xi,(\gamma_{i})_{i\in I}}=
  S_{I,W_{x},x,\xi |_{W_{x}},(\gamma_{i})_{i\in I}}
  =S_{I,W_{f},x_{f},\xi_{f},(\gamma_{i})_{i\in I}}.\end{gather}
This shows that  $S_{I,W,x,\xi,(\gamma_{i})_{i\in I}}$ depends only on   $I$, $f$, and $(\gamma_{i})_{i\in I}$. 
\cqfd

 The previous lemma allows to define the following simplified notation (already introduced in  \defiref{defi-constr-excursion-intro}). 
 
 \begin{defi}
 For every function  $f\in  \mc O(\wh G\backslash (\wh G)^{I}/\wh G)$ we  set 
 \begin{gather} S_{I,f,(\gamma_{i})_{i\in I}}=S_{I,W,x,\xi,(\gamma_{i})_{i\in I}}\in \mc B \end{gather}
  where $W,x,\xi$ are such that  $f$ satisfies 
 \eqref{intro-def-f}. 
 \end{defi}
 
 This  new  notation allows a more synthetic formulation of  properties     \eqref{SIW-p1}, \eqref{SIW-p2} and \eqref{SIW-p3}, 
 in the form of  properties (i), (ii), (iii) and (iv) of   \propref{prop-SIf-i-ii-iii} that we are now able to justify.

      \noindent {\bf Proof of properties (i), (ii), (iii) and (iv) of   \propref{prop-SIf-i-ii-iii}.}    The fact that    $S_{I,f,(\gamma_{i})_{i\in I}}$ depends only on the image of   
 $(\gamma_{i})_{i\in I}$ in   $\pi_{1}(X\sm N, \ov\eta)^{I}$ is proved in the next lemma. We deduce   (ii) from \eqref{SIW-p1}. To  prove   (i) we use 
\eqref{SIW-p2} with $I_{1}=I_{2}=I$,  and we apply   \eqref{SIW-p1} to the obvious map $\zeta :I\cup I\to I$. Property  (iii) comes from   \eqref{SIW-p3}, since  Zorro lemma \eqref{zorro} implies that  
    $$\s{ \xi \boxtimes \on{ev}_{W}, \big((g_{i})_{i\in I} 
  \boxtimes  (g'_{i})_{i\in I} \boxtimes  (g''_{i})_{i\in I} \big)\cdot ( \delta_{W} \boxtimes x)}
  =\s{\xi, (g_{i}(g'_{i})^{-1}g''_{i})_{i\in I} \cdot x}.$$
  Lastly (iv) comes from the fact that  $H_{I,W}$ is an inductive limit of  finite  dimensional continuous representations   of $\on{Gal}(\ov F/F)^{I}$.  \cqfd

The following lemma was found by     B\"ockle,    Harris,  Khare and  Thorne and is used in their  article \cite{boeckle-harris...}. 

\begin{lem}
\label{prop-harris} 
   For any  $I$ and $f\in \mc O(\wh G\backslash (\wh G)^{I}/\wh G)$,    $S_{I,f,(\gamma_{i})_{i\in I}}$ depends only on the image of    
 $(\gamma_{i})_{i\in I}$ in   $\pi_{1}(X\sm N, \ov\eta)^{I}$, and 
 $(\gamma_{i})_{i\in I}\mapsto S_{I,f,(\gamma_{i})_{i\in I}}$ is continuous from the profinite group $\pi_{1}(X\sm N, \ov\eta)^{I}$ to the  $E$-algebra of finite   dimension   
   $\mc B$ equipped with the  $E$-adic topology. 
\end{lem}
 \dem    
          Let $v$ be a  place of $X\sm N$. We fix  an embedding $\ov F\subset \ov {F_{v}}$, whence an  inclusion 
            $ \on{Gal}(\ov {F_{v}}/F_{v})\subset  \on{Gal}(\ov F/F)$.  Let $I_{v}=\on{Ker}( \on{Gal}(\ov {F_{v}}/F_{v})\to \wh \Z)$ be the inertia subgroup  at   $v$. Then 
for $I, W,x,\xi$ as in \eqref{excursion-def-intro}, the image of  
the  composition $  H_{\{0\},\mbf  1}\xrightarrow{\mc H(x)}
 H_{\{0\},W^{\zeta_{I}}}\isor{\chi_{\zeta_{I}}^{-1}} 
  H_{I,W}$ (which is the beginning  of \eqref{excursion-def-intro}) consists of elements invariant  by 
    $(I_{v})^{I}$, since  the creation  operators 
      are  morphisms of sheaves  over the whole $\Delta(X\sm N)$  (and  in particular  over $\Delta(v)$).   Thus  for $(\gamma_{i})_{i\in I}\in \on{Gal}(\ov F/F)^{I}$ and 
    $(\delta_{i})_{i\in I}\in (I_{v})^{I}$ we have 
    \begin{gather}\label{rel-S-gammai-deltai}S_{I,W,x,\xi,(\gamma_{i})_{i\in I}}=S_{I,W,x,\xi,(\gamma_{i}\delta_{i})_{i\in I}}.\end{gather} 
  This is true for any embedding $\ov F\subset \ov {F_{v}}$ 
    (in fact,    by \remref{rem-gamma-en-plus},  \eqref{rel-S-gammai-deltai}  for an embedding implies \eqref{rel-S-gammai-deltai} for all embeddings). But  $\pi_{1}(X\sm N, \ov\eta)$ is the topological quotient     of 
    $\on{Gal}(\ov F/F)$ by the closed subgroup generated  
        by the  $I_{v}$ for $v\in (X\sm N)$ and their conjugates. \cqfd

 We do not know   if $\mc B$ is reduced. Nevertheless we get  a spectral decomposition  
  (i.e.  a decomposition into  generalized eigenspaces) 
  \begin{gather}\label{dec-nu-intro-intro}
  H_{\{0\},\mbf 1} =\bigoplus_{\nu} \mf H_{\nu}\end{gather}
 where in the RHS the direct sum  is indexed  by the  characters  $\nu$ of  $\mc B$. Increasing  $E$ if necessary,  we assume that all the characters of $\mc B$ are defined  over $E$. 
 
  The  following  proposition allows to obtain   decomposition \eqref{intro3-dec-canonical}  from \eqref{dec-nu-intro-intro} because it associates to every character $\nu$ a Langlands parameter $\sigma$.

 \begin{prop} \label{intro-Xi-n}  For every character  $\nu$ of  $\mc B$ there exists  a    morphism 
       $\sigma:\pi_{1}(X\sm N, \ov\eta)\to \wh G(\Qlbar)$
such that  
\begin{itemize}
\item [] (C1) $\sigma$ takes its values   in  $\wh G(E')$, where $E'$ is a finite  extension  of $E$, and it  is continuous, 
\item [] (C2) $\sigma$ is semisimple, i.e.   if its image is included in a parabolic, 
it is included in an associated Levi
(since $\Qlbar$ has characteristic  $0$ this is equivalent to say that the Zariski closure of its image is reductive   \cite{bki-serre}),  
\item[] (C3) 
for every  $I$ and    $f\in \mc O(\wh G\backslash (\wh G)^{I}/\wh G)$,  we have  
$$\nu(S_{I,f,(\gamma_{i})_{i\in I}})= f\big((\sigma(\gamma_{i}))_{i\in I} \big).
$$
\end{itemize}
Moreover   $\sigma$ is unique up to conjugation  by   $\wh G(\Qlbar)$. 
\end{prop}
\noindent

   \noindent {\bf Proof.}   We refer to  the proof of   proposition 10.7 of \cite{coh} for  some additional details. 
   The proof uses only  properties (i), (ii), (iii) and (iv) of   \propref{prop-SIf-i-ii-iii}. 
    Let $\nu$ be a character of   $\mc B$.

For every $n\in \N$ 
 we denote by 
   $(\wh G)^{n}\modmod \wh G$ the  coarse quotient of  $(\wh G)^{n}$ by the action  of  $\wh G$  by  diagonal conjugation, i.e.    
   $$h.(g_{1},...,g_{n})=(hg_{1}h^{-1},...,hg_{n}h^{-1}).$$ 
   Then the morphism 
   \begin{gather*} (\wh G)^{n}\to (\wh G)^{\{0,...,n\}}, (g_{1},...,g_{n})\mapsto 
   (1,g_{1},...,g_{n})\end{gather*}
 induces  an  isomorphism  $$\beta: (\wh G)^{n}\modmod \wh G\isom \wh G\backslash (\wh G)^{\{0,...,n\}}/\wh G ,$$ whence an algebra   isomorphism     $$\mc O( (\wh G)^{n}\modmod \wh G)\isom \mc O(\wh G\backslash (\wh G)^{\{0,...,n\}}/\wh G ), \ \ f\mapsto f \circ \beta^{-1}. $$ 

We introduce   \begin{align} \nonumber \Theta_{n}^{\nu}: \mc O((\wh G)^{n}\modmod \wh G)&\to 
C(\pi_{1}(X\sm N, \ov\eta)^{n}, E)\\ \nonumber 
f& \mapsto [(\gamma_{1},...,\gamma_{n})\mapsto 
\nu(S_{I,f \circ \beta^{-1},(1,\gamma_{1},...,\gamma_{n})})]\end{align}

Condition (C3) that  $\sigma$ must satisfy can be reformulated in the following way : for every $n$ and for every  $f\in \mc O((\wh G)^{n}\modmod \wh G)$, 
\begin{gather}\label{rel-sigma-nu}\Theta_{n}^{\nu}(f)=[(\gamma_{1},...,\gamma_{n})\mapsto f((\sigma(\gamma_{1}),...,\sigma(\gamma_{n})))   ]. \end{gather}

We deduce immediately from  properties (i), (ii), (iii) and (iv) of    \propref{prop-SIf-i-ii-iii} that the sequence 
$(\Theta_{n}^{\nu})_{n\in \N^{*}}$ 
satisfies  the following properties
\begin{itemize} \item 
for every $n$, $\Theta_{n}^{\nu}$ is an algebra   morphism, 
                                                 \item  the sequence  $(\Theta_{n}^{\nu})_{n\in \N^{*}}$ is functorial relatively to  all the  maps between the sets  $\{1,...,n\}$, i.e.      for $m,n\in \N^{*}$, 
               $$\zeta: \{1,...,m\}\to \{1,...,n\}$$ arbitrary,   
               $f\in \mc O((\wh G)^{m}\modmod \wh G)$ and 
               $(\gamma_{1},...,\gamma_{n })\in \pi_{1}(X\sm N, \ov\eta)^{n }$, 
               we have  
             $$\Theta_{n}^{\nu}( f^{\zeta}) ((\gamma_{j})_{j\in \{1,...,n\}})=
             \Theta_{m}^{\nu}(f)((\gamma_{\zeta(i)})_{i\in \{1,...,m\}})$$ 
            where  $f^{\zeta}\in \mc O((\wh G)^{n}\modmod \wh G)$  is defined by 
               $$f^{\zeta}((g_{j})_{j\in \{1,...,n\}})=f((g_{\zeta(i)})_{i\in \{1,...,m\}}), $$
 \item     for $n\geq 1$, 
  $f\in \mc O((\wh G)^{n}\modmod \wh G)$  
    and   $(\gamma_{1},...,\gamma_{n+1})\in \pi_{1}(X\sm N, \ov\eta)^{n+1}$ we have  
   $$\Theta_{n+1}^{\nu}( \wh f)(\gamma_{1},...,\gamma_{n+1})=
   \Theta_{n}^{\nu}( f)(\gamma_{1},...,\gamma_{n}\gamma_{n+1}) $$
 where $\wh f\in  \mc O((\wh G)^{n+1}\modmod \wh G)$ is defined by 
   $$\wh f(g_{1},...,g_{n+1})=f(g_{1},...,g_{n}g_{n+1}). $$
\end{itemize}
 To justify  the last property, we apply   property 
   (iii) of  \propref{prop-SIf-i-ii-iii} to  $$I=\{0,...,n \}, 
   (\gamma_{i})_{i\in I}=(1,\gamma_{1},...,\gamma_{n}), (\gamma'_{i})_{i\in I}=(1)_{i\in I}, (\gamma''_{i})_{i\in I}=(1, ...,1,\gamma_{n+1})$$ and we use  (ii) to delete all the   $1$ except the first one in  
   $(\gamma_{i})_{i\in I}\times (\gamma'_{i})_{i\in I}\times (\gamma''_{i})_{i\in I}$.

We will see that these properties of 
 the sequence 
$(\Theta_{n}^{\nu})_{n\in \N^{*}}$ 
imply  the existence and the unicity of   $\sigma$ satisfying (C1), (C2)  and (C3) (i.e. 
\eqref{rel-sigma-nu}). 

For $G=GL_{r}$ the result is already known:  the sequence $(\Theta_{n}^{\nu})_{n\in \N^{*}}$  
is determined by $\Theta_{1}^{\nu}(\on{Tr})$ 
(which must be  the character of $\sigma$) 
and  
$\Lambda^{r+1}\mr{St}=0$ implies  the   pseudo-character  relation 
whence the existence of $\sigma$ by     \cite{taylor}.  We refer to  remark  11.8 of \cite{coh} for further details. 

For general $G$ we use  results of   \cite{richardson}. 
We say that a   $n$-uplet   $(g_{1},...,g_{n})\in
 \wh G(\Qlbar)^{n}$ is semisimple if 
 every parabolic subgroup    containing it 
 admits  an associated Levi  subgroup    containing it. Since $  \Qlbar$ is of characteristic $0$ this is equivalent  (see \cite{bki-serre}) to the condition that 
 the Zariski closure  
  $\ov{<g_{1},...,g_{n}>}$ of the  subgroup  $<g_{1},...,g_{n}>$ generated by   $g_{1},...,g_{n}$ is reductive. 
 By theorem 3.6 of  \cite{richardson} the $\wh G$-orbit  (by conjugation) of $(g_{1},...,g_{n})$ is closed in 
 $(\wh G)^{n}$ if and only if $(g_{1},...,g_{n})$ is semisimple. 
Therefore  the points over $\Qlbar$ of the coarse quotient $(\wh G)^{n}\modmod \wh G$ (which correspond to  the 
  closed   $\wh G$-orbits defined over    $\Qlbar$ in   $(\wh G)^{n}$)  are  in  bijection with  the   $\wh G(\Qlbar)$-conjugacy classes    of semisimple  $n$-uplets    $(g_{1},...,g_{n})\in
 \wh G(\Qlbar)^{n}$.

 For any $n$-uplet 
   $(\gamma_{1},...,\gamma_{n})\in \pi_{1}(X\sm N, \ov\eta)^{n}$ we denote by   
  $\xi_{n} (\gamma_{1},...,\gamma_{n})$ the point defined  over $\Qlbar$ of 
   the  coarse quotient $(\wh G)^{n}\modmod \wh G$ associated to  the character 
   $$\mc O((\wh G)^{n}\modmod \wh G)\to \Qlbar,  \ f\mapsto \Theta_{n}^{\nu}( f)(\gamma_{1},...,\gamma_{n}).$$ 
   We denote by  
 $\xi_{n}^{\mr{ss}}(\gamma_{1},...,\gamma_{n})$ the  conjugacy class of   
   semisimple $n$-uplets  corresponding to   $\xi_{n}(\gamma_{1},...,\gamma_{n})$ by the result of  \cite{richardson} stated above.  

The  relation  \eqref{rel-sigma-nu} is equivalent to the condition that for every  $n$ and for every 
$(\gamma_{1},...,\gamma_{n})$, 
$ (\sigma(\gamma_{1}),...,\sigma(\gamma_{n}))\in (\wh G(\Qlbar))^{n}$ 
(which is not always  semisimple)
lies over  
$\xi_{n}(\gamma_{1},...,\gamma_{n})$. 

\noindent {\bf Unicity of $\sigma$ (up to conjugation).}  We choose $n$ and $(\gamma_{1},...,\gamma_{n})$ such that $\sigma(\gamma_{1}),...,\sigma(\gamma_{n})$ generate a Zariski dense subgroup  in  $\ov{ \mr{Im}(\sigma)}$. 
Since $\sigma$ is assumed to be semisimple, $(\sigma(\gamma_{1}),...,\sigma(\gamma_{n}))$ is semisimple. 
We fix 
$(g_{1},...,g_{n})$ in $\xi_{n}^{\mr{ss}}(\gamma_{1},...,\gamma_{n})$. 
Then $(\sigma(\gamma_{1}),...,\sigma(\gamma_{n}))$ is conjugated to 
$(g_{1},...,g_{n})$ and by conjugating  $\sigma$ we can assume the latter is equal to the former. Then $\sigma$ is uniquely determined because for every $\gamma$, 
$\sigma(\gamma)$  belongs to the Zariski closure of the subgroup generated by $(g_{1},...,g_{n})$ and 
$(g_{1},...,g_{n},\sigma(\gamma)) \in 
\xi_{n+1}^{\mr{ss}}(\gamma_{1},...,\gamma_{n},\gamma)$, therefore the knowledge  of $\xi_{n+1} (\gamma_{1},...,\gamma_{n},\gamma)$ determines  $\sigma(\gamma)$ uniquely. 

\noindent {\bf Existence of $\sigma$.}  
  For every $n$ and every  $(\gamma_{1},...,\gamma_{n})\in \pi_{1}(X\sm N, \ov\eta)^{n}$ we choose $(g_{1},...,g_{n})\in \xi_{n}^{\mr{ss}}(\gamma_{1},...,\gamma_{n})$ (well defined  up to conjugation).  Then we   choose   $n$ and 
   $(\gamma_{1},...,\gamma_{n})\in \pi_{1}(X\sm N, \ov\eta)^{n}$ such that
       \begin{itemize}
   \item    (H1)
  the  dimension of 
  $\ov{<g_{1},...,g_{n}>}$ is the greatest  possible 
  \item     (H2) the centralizer   $C(g_{1},...,g_{n})$ of   $<g_{1},...,g_{n}>$ is the smallest   possible (minimal dimension and then  minimal number of connected components).  
  \end{itemize}
  We fix  $(g_{1},...,g_{n})\in \xi_{n}^{\mr{ss}}(\gamma_{1},...,\gamma_{n})$ for the rest  of the proof and we construct a map   $$\sigma: \pi_{1}(X\sm N, \ov\eta)\to \wh G(\Qlbar)$$  by asking  that for every   $\gamma\in \pi_{1}(X\sm N, \ov\eta)$, $\sigma(\gamma)$ is the unique element $g$ of  $\wh G(\Qlbar)$ such that  $(g_{1},...,g_{n},g)\in \xi_{n+1}^{\mr{ss}}(\gamma_{1},...,\gamma_{n},\gamma)$.  
  The existence and the unicity of such a $g$ are justified in the following way. 
  \begin{itemize}
  \item  {\bf A) Existence of $g$} : for $(h_{1},...,h_{n},h)\in \xi_{n+1}^{\mr{ss}}(\gamma_{1},...,\gamma_{n},\gamma)$, $(h_{1},...,h_{n})$ is necessarily semisimple    (because $(h_{1},...,h_{n})$ is over  $\xi_{n}(\gamma_{1},...,\gamma_{n})$ and $(g_{1},...,g_{n})\in \xi_{n}^{\mr{ss}}(\gamma_{1},...,\gamma_{n})$
  thus by  theorem 5.2 of \cite{richardson}, $\ov{<h_{1},...,h_{n}>}$ admits a Levi  subgroup  isomorphic to $\ov{<g_{1},...,g_{n}>}$, but 
  $\dim(\ov{<h_{1},...,h_{n}>})\leq \dim(\ov{<g_{1},...,g_{n}>})$
  by (H1)), 
  thus conjugating  $(h_{1},...,h_{n},h)$ we can assume  that $(h_{1},...,h_{n})=(g_{1},...,g_{n})$ and then we take  $g=h$.   
  \item {\bf B) Unicity of $g$} : we have $C(g_{1},...,g_{n},g)\subset C(g_{1},...,g_{n})$ and   equality holds by (H2), therefore $g$ commutes with $C(g_{1},...,g_{n})$ and since it was well defined up to conjugation by $C(g_{1},...,g_{n})$ it is unique. 
  \end{itemize}
  
Then we  show   that the map  $\sigma$ we have just constructed  is a morphism of groups. 
Indeed let  $\gamma, \gamma'\in \pi_{1}(X\sm N, \ov\eta)$. The same argument as in A) above shows that there exist $g,g'$ such that 
\begin{gather}\label{semisimple-g-g'}(g_{1},...,g_{n},g,g') \in \xi_{n+2}^{\mr{ss}}(\gamma_{1},...,\gamma_{n},\gamma,\gamma').\end{gather}
   Thanks to the   properties satisfied  by the sequence 
      $(\Theta_{n}^{\nu})_{n\in \N^{*}}$ we see that 
      $\xi_{n+1}(\gamma_{1},...,\gamma_{n},\gamma\gamma')$ is the image of  $\xi_{n+2}(\gamma_{1},...,\gamma_{n},\gamma,\gamma')$ by the morphism 
   $$(\wh G)^{n+2}\modmod \wh G \to (\wh G)^{n+1}\modmod \wh G, 
   (h_{1},...,h_{n},h,h')\mapsto (h_{1},...,h_{n},hh').$$
       From this we deduce  that 
        $   (g_{1},...,g_{n},gg')$ is over  
        $  \xi_{n+1} (\gamma_{1},...,\gamma_{n},\gamma\gamma')$. Moreover 
        $ (g_{1},...,g_{n},gg')$ is semisimple by the same argument as in A), because 
       \begin{gather*}  \dim(\ov{< g_{1},...,g_{n},gg' >}) 
        \leq  \dim(\ov{< g_{1},...,g_{n},g,g' >}) 
         \leq  \dim(\ov{< g_{1},...,g_{n}  >}) \end{gather*}
          (where the last  inequality comes from  (H1)). Therefore 
          $ (g_{1},...,g_{n},gg')$ belongs to 
        $  \xi_{n+1}^{\mr{ss}} (\gamma_{1},...,\gamma_{n},\gamma\gamma')$
          and 
           $gg'=\sigma(\gamma \gamma')$. 
        The same arguments show  that 
     \begin{gather}\label{gg'-gamma-gamma'}g=\sigma(\gamma)\text{ \ \   and \ \ }  g'=\sigma(\gamma').\end{gather} 
Endly we showed that $\sigma(\gamma \gamma') =\sigma(\gamma)\sigma(\gamma')$. 
  
  Thus  $\sigma$ is a group morphism with values  in $\wh G(E')$
  (where $E'$ is a finite  extension of $E$ such that $g_{1},...,g_{n}$ belong to $\wh G(E')$).    
The argument to prove that $\sigma$ is continuous is the following. 
We know  that  for every function $f$ over $(\wh G_{E'})^{n+1}\modmod \wh G_{E'}$, the map $$
\pi_{1}(X\sm N, \ov\eta)\to E',  \ \ 
\gamma\mapsto f(g_{1},...,g_{n}, \sigma(\gamma))= 
\Theta_{n+1}^{\nu}( f)(\gamma_{1},...,\gamma_{n},\gamma)$$
 is continuous. 
But  the morphism 
 \begin{align*}\mc O((\wh G_{E'})^{n+1}\modmod \wh G_{E'}) &\to  \mc O(\wh G_{E'}\modmod C(g_{1},...,g_{n})) \\ f & \mapsto  [g\mapsto f(g_{1},...,g_{n},g)]
\end{align*} is surjective
(because $(g_{1},...,g_{n})$ is semisimple, therefore its orbit by conjugation is an affine closed subvariety  of 
$(\wh G_{E'})^{n}$, 
isomorphic to $\wh G_{E'}/C(g_{1},...,g_{n})$). 
Moreover, if we denote by 
$D(g_{1},...,g_{n})$ the centralizer  of $C(g_{1},...,g_{n})$ (which contains  the image of $\sigma$), the restriction morphism  
$$\mc O(\wh G_{E'}\modmod C(g_{1},...,g_{n})) =
\mc O(\wh G_{E'} )^{C(g_{1},...,g_{n})} \to \mc O(D(g_{1},...,g_{n}))$$ is surjective  because  the   restriction $\mc O(\wh G_{E'}) \to \mc O(D(g_{1},...,g_{n}))$  is  obviously   surjective, 
and remains so when we take   the invariants by the  reductive  group    $C(g_{1},...,g_{n})$ (acting by conjugation), and that 
$C(g_{1},...,g_{n})$ acts trivially on  $\mc O(D(g_{1},...,g_{n}))$.  
Thus for every function $h\in \mc O(D(g_{1},...,g_{n}))$, 
the map 
$$
\pi_{1}(X\sm N, \ov\eta)\to E',  \ \ \gamma\mapsto h(\sigma(\gamma))$$
 is continuous, and we proved that $\sigma$ is continuous. 
 
 It remains to  prove \eqref{rel-sigma-nu}, i.e.   for 
   $m\in \N^{*}$, $f\in \mc O((\wh G)^{m}\modmod \wh G)$ and  $(\delta_{1},...,\delta_{m})\in \pi_{1}(X\sm N, \ov\eta)^{m}$, we have  
\begin{gather}\nonumber f(\sigma(\delta_{1}),...,\sigma(\delta_{m}))=
 \big(\Theta^{\nu}_{m}(f)\big)(\delta_{1},...,\delta_{m}).\end{gather}
By the same  arguments as  for \eqref{semisimple-g-g'} and \eqref{gg'-gamma-gamma'} we show  that  
 $$(g_{1},...,g_{n},\sigma(\delta_{1}),...,\sigma(\delta_{m}))\in \xi_{n+m}^{ss}(\gamma_{1},...,\gamma_{n},\delta_{1},...,\delta_{m} ). $$
 Therefore 
 $(\sigma(\delta_{1}),...,\sigma(\delta_{m}))$ is over   $\xi_{m}(\delta_{1},...,\delta_{m} )$. \cqfd

Thus we have obtained the decomposition \eqref{intro1-dec-canonical}. This concludes  the proof of \thmref{intro-thm-ppal}, provided we admit the following result,  which will be justified   in   section \ref {subsection-intro-decomp}. 

 \noindent{\bf Result admitted for the moment} (\propref{S-non-ram-concl-intro}). The decomposition \eqref{intro1-dec-canonical}  is compatible with the Satake isomorphism at all the places of $X\sm N$.

            \section{Creation and annihilation morphisms}
           \label{subsection-crea-annihil-intro}
            From now one and until  section \ref{subsection-intro-decomp} we will justify              \begin{itemize}
             \item   the two results admitted   in  the section \ref{section-esquisse-abc}
 and  which were used for the proof  of   \propref{prop-a-b-c}
       \item and   the result admitted at the end of the previous section  and  which was used to finish the proof of  \thmref{intro-thm-ppal}. 
                     \end{itemize}

    In this section our goal is to use the   coalescence isomorphisms   \eqref{intro-isom-coalescence} to construct   creation and annihilation morphisms, and then to rewrite the   Hecke operators at the  places of  $X\sm N$ as the composition   
 \begin{itemize}
 \item of a creation morphism,
  \item of  the action of a   partial Frobenius morphism,
  \item of an   annihilation morphism,
  \end{itemize}
and to use this to extend  the  Hecke operators \eqref{defi-Tf} to morphisms of sheaves over the whole   $(X\sm N)^{I}$  and to obtain  the Eichler-Shimura relations.

  Let   $I$ and $J$ be finite sets. 
  We will now define  the creation and annihilation morphisms, in the following way. The legs indexed by $I$ will remain unchanged and we will create (or annihilate)  the legs indexed by $J$ at the  same point of the curve (indexed by a set with one element, that we will denote by  $\{0\}$).

  We have obvious maps  
   $$\zeta_{J} :J \to \{0\},  \ 
  \zeta_{J}^{I}=(\Id_{I},\zeta_{J}): I\cup J\to I\cup \{0\} \text{ \   and   \ } \zeta_{\emptyset}^{I}=(\Id_{I},\zeta_{\emptyset}): I \to I \cup \{0\}.$$ 
  Let $W$ and $U$ be   finite dimensional $E$-linear representations of 
$(\wh G)^{I}$ and   $(\wh G)^{J}$ respectively. 
We recall that   $U^{\zeta_{J}}$ is 
the representation of $\wh G$ obtained by restricting 
  $U$ to the diagonal    $\wh G\subset (\wh G)^{J}$. 
Let  $x\in U$ and $\xi\in U^{*}$ be invariant under  the diagonal action  of $\wh G$.  
Then $W\boxtimes U$ is a representation of  $(\wh G)^{I\cup J}$ and 
 $W\boxtimes U^{\zeta_{J}}$ and  $W\boxtimes \mbf 1$ are  representations of  $(\wh G)^{I\cup \{0\}}$ linked  by the   morphisms 
 $$W\boxtimes \mbf 1\xrightarrow{\Id_{W}\boxtimes x}W\boxtimes U^{\zeta_{J}}\text{ \  and \ }  W\boxtimes U^{\zeta_{J}}\xrightarrow{\Id_{W}\boxtimes \xi}W\boxtimes \mbf 1.$$ 
 
 We denote by   $\Delta:X\to X^{J}$ the   diagonal morphism and we  write  $E_{X\sm N}$  for  the constant sheaf  over  $X\sm N$. 
         \begin{defi}
 We define   the  creation morphism
  $
\mc C_{x}^{\sharp}$ as the  composition   \begin{gather*}
  \mc H _{N,I,W}^{\leq\mu} \boxtimes E_{(X\sm N)} 
  \isor{\chi_{\zeta_{\emptyset}^{I}}}
 \mc H _{N,I\cup\{0\},W\boxtimes \mbf 1}^{\leq\mu}  \\
 \xrightarrow{\mc H(\Id_{W}\boxtimes x) }
 \mc H _{N,I\cup\{0\},W\boxtimes U^{\zeta_{J}}}^{\leq\mu} 
 \isor{ \chi_{\zeta_{J}^{I}}^{-1} }
 \restr{ \mc H _{N,I\cup J,W\boxtimes U}^{\leq\mu}}{(X\sm N)^{I}\times \Delta(X\sm N)}  
   \end{gather*}
   where $\chi_{\zeta_{\emptyset}^{I}}$ and $\chi_{\zeta_{J}^{I}}$ are  the   coalescence isomorphisms   of \eqref{intro-isom-coalescence}. 
   Similarly we define the  annihilation morphism $
\mc C_{\xi}^{\flat}$  as the  composition 
  \begin{gather*}
  \restr{ \mc H _{N,I\cup J,W\boxtimes U}^{\leq\mu}}{(X\sm N)^{I}\times \Delta(X\sm N)} 
  \isor{ \chi_{\zeta_{J}^{I}}}
 \mc H _{N,I\cup\{0\},W\boxtimes U^{\zeta_{J}}}^{\leq\mu} \\
 \xrightarrow{ \mc H(\Id_{W}\boxtimes \xi) }
  \mc H _{N,I\cup\{0\},W\boxtimes \mbf 1}^{\leq\mu} 
 \isor{\chi_{\zeta_{\emptyset}^{I}}^{-1} }
   \mc H _{N,I,W}^{\leq\mu} \boxtimes E_{(X\sm N)} . 
    \end{gather*}
    \end{defi}
    
   All the  morphisms above are    morphisms of sheaves over  
    $(X\sm N)^{I} \times (X\sm N)$.

 Now we will use these  morphisms with $J=\{1,2\}$. 
      Let $v$  be a  place in   $|X|\sm |N|$.    We consider  $v$ also  as a  subscheme of  $X$ and we denote by  $E_{v}$ the constant sheaf  over   $v$.    Let $V$ be an irreducible representation of $\wh G$.    As previously we denote by   $\mbf 1\xrightarrow{\delta_{V}} V\otimes V^{*}$ and $ V\otimes V^{*}\xrightarrow{\on{ev}_{V}} \mbf 1$ the  natural morphisms. 
        
      For  $\kappa$ big enough  (in function of $\deg(v),V$), we define   $S_{V,v}$ as the  composition 
  \begin{gather}\label{def-SVv-intro1}
 \mc H _{N,I,W}^{\leq\mu} \boxtimes E_{v} \\ \label{def-SVv-intro2}
  \xrightarrow{ \restr{\mc C_{
    \delta_{V}}^{\sharp}}{(X\sm N)^{I}\times v}}
 \restr{ \mc H _{N,I\cup\{1,2\},W\boxtimes V\boxtimes V^{*}}^{\leq\mu}}{(X\sm N)^{I}\times \Delta(v)} \\ \label{def-SVv-intro3}
 \xrightarrow{ \restr{(F_{\{1\}})^{\deg(v)} }{(X\sm N)^{I}\times \Delta(v)}}
 \restr{ \mc H _{N,I\cup\{1,2\},W\boxtimes V\boxtimes V^{*}}^{\leq\mu+\kappa}}{(X\sm N)^{I}\times \Delta(v)} \\ \label{def-SVv-intro4}
  \xrightarrow{ \restr{\mc C_{
    \on{ev_{V}}}^{\flat }}{(X\sm N)^{I}\times v}}
     \mc H _{N,I,W}^{\leq\mu+\kappa}  \boxtimes E_{v} . 
 \end{gather}
In other words we create two new  legs at  $v$ with the help  of  
 $ \delta_{V}:\mbf 1\to V\otimes V^{*}$, we apply the  partial Frobenius morphism (to the power $\deg(v)$) to the first one, and then we annihilate them  with the help  of  $ \on{ev_{V}}: V\otimes V^{*}\to\mbf 1$. 
  
  As a   morphism of  constructible sheaves over 
     $(X\sm N)^{I}\times v$,  $S_{V,v}$ commutes with the natural action  of the partial   Frobenius  morphism on  $E_{v}$ in  \eqref{def-SVv-intro1} and  \eqref{def-SVv-intro4}, since 
     \begin{itemize}
     \item the  creation and annihilation morphisms intertwine this  action with the action  of $F_{\{1,2\}}$ over  \eqref{def-SVv-intro2} and \eqref{def-SVv-intro3}, by  \propref{rem-coalescence-frob-compat}, 
     \item  $F_{\{1\}}$ and therefore   $F_{\{1\}}^{\deg(v)}$ commute  with  $F_{\{1,2\}}=F_{\{1\}}F_{\{2\}}$. 
     \end{itemize}
  \begin{defi} By abuse of notations we still write    
        $$S_{V,v}: \mc H _{N,I,W}^{\leq\mu}\to 
         \mc H _{N,I,W}^{\leq\mu+\kappa}$$  for the  morphism 
    of sheaves over   $(X\sm N)^{I}$ obtained by descent relatively to the action of 
    $\Z/\deg(v)\Z$ (i.e. by taking  the invariants by  the natural action  of the partial Frobenius morphism  on  $E_{v}$ in  \eqref{def-SVv-intro1} and  \eqref{def-SVv-intro4}). 
     \end{defi}
  
 \begin{prop}\label{prop-coal-frob-cas-part-intro} 
  The restriction of  $S_{V,v}$  to  $(X\sm (N \cup v))^{I}$
is equal, as  a  morphism of sheaves over  $(X\sm (N \cup v))^{I}$,  to the Hecke operator $$T(h_{V,v}):  \restr{\mc H _{N,I,W}^{\leq\mu}}{(X\sm (N \cup v))^{I}}\to 
       \restr{  \mc H _{N,I,W}^{\leq\mu+\kappa}}{(X\sm (N \cup v))^{I}}
       . $$
       \end{prop}
\noindent
 It is enough to prove this  when   $V$  and $W$ are irreducible. The proof is of geometric nature.   We sketch it here in a simple case, where it is reduced to the intersection of two smooth substacks  in a smooth Deligne-Mumford stack   and where this intersection happens to be transverse. The proof  is more complicated in  general because of the singularities. We refer to  the proof of  proposition 6.2 of \cite{coh} for the general case (but an alternative  solution could consist to reduce to the  situation of a smooth transverse   intersection with the help     of Bott-Samelson resolutions). 
  
\noindent {\bf Proof when  $V$ is minuscule  and  $\deg(v)=1$.} 
We recall that an  irreducible representation of $\wh G$ is said to be  minuscule if all its weights are conjugated by the  Weyl group. This is equivalent to the property that the corresponding orbit  in the affine grassmannian is closed  (and therefore it implies that the  corresponding closed stratum    is smooth). We denote by  $d$ the dimension of the orbit associated to  $V$.

Thanks to  the hypothesis that  $\deg(v)=1$ we can delete 
$\boxtimes E_{v}$ everywhere. 
We consider the   Deligne-Mumford stack   
$$\mc Z^{(\{1\},\{2\}, I)}=\restr{\Cht_{N,I \cup \{1,2\},W \boxtimes V\boxtimes V^{*}}^{(\{1\},\{2\},I)}}
 {(X\sm (N\cup v))^{I}\times \Delta(v)}.$$
We will construct two closed substacks   $\mc Y_{1}$ and $\mc Y_{2}$ in $\mc Z^{(\{1\},\{2\}, I)}$,    equipped with   morphisms 
 $\alpha_{1}$ and $\alpha_{2}$ towards   $$ \mc Z^{(I)}=\restr{\Cht_{N,I ,W }^{(I)}}
 {(X\sm (N\cup v))^{I}}$$ in such a way that  
 \begin{itemize}
 \item {\bf A)} the  restriction to  $(X\sm (N\cup v))^{I}$ of the  composition 
\eqref{def-SVv-intro1}$\to$\eqref{def-SVv-intro2}$\to$\eqref{def-SVv-intro3}
of the   creation morphism
 and of the action of the partial Frobenius morphism  is realized by a  cohomological correspondence   supported by the correspondence $\mc Y_{2}$ from   $\mc Z^{(I)}$ to     $\mc Z^{(\{1\},\{2\}, I)}$, and whose  restriction to the open smoothness locus is determined by its support, with a corrective  coefficient of  $q^{-d/2}$, 

  \item {\bf B)} the  restriction to $(X\sm (N\cup v))^{I}$ of the annihilation morphism    \eqref{def-SVv-intro3}$\to$\eqref{def-SVv-intro4}  is realized by a   cohomological correspondence   supported by    the   correspondence $\mc Y_{1}$ from  $\mc Z^{(\{1\},\{2\}, I)}$ to $\mc Z^{(I)}$, and whose  restriction to the open smoothness locus is determined by its support.   
 \end{itemize}
Therefore $S_{V,v}$ will be  realized by a  cohomological correspondence supported by   the product  
    $\mc Y_{1}\times_{\mc Z^{(\{1\},\{2\}, I)}}\mc Y_{2}$ of these correspondences. We will see 
    \begin{itemize}
    \item that  the  product $\mc Y_{1}\times_{\mc Z^{(\{1\},\{2\}, I)}}\mc Y_{2}$ is nothing but  the   Hecke correspondence   $\Gamma^{(I)}$ of  
$ \mc Z^{(I)}$ to itself  (which is a finite correspondence, and even étale)
\item that $S_{V,v}$,  which is  therefore a  cohomological correspondence supported by $\Gamma^{(I)}$ is in fact  equal to  the  obvious cohomological correspondence   supported by $\Gamma^{(I)}$ with a corrective  coefficient of  $q^{-d/2}$ (which realizes 
  $T(h_{V,v})$ since $V$ is minuscule). \end{itemize}
  Thanks to  the hypothesis that  $V$ is minuscule it will suffice to do this computation over the open  smoothness locus, and the computation  will be immediate because we will see that over this open  smoothness locus  the intersection $\mc Y_{1}\times_{\mc Z^{(\{1\},\{2\}, I)}}\mc Y_{2}$ 
   is a  transverse intersection of two smooth substacks.

 We construct now all these objects. 
 The Hecke correspondence  $\Gamma^{(I)}$ is the  stack 
classifying the data  of   $(x_i)_{i\in I}$ and of a diagram 
\begin{gather} \label{intro-diag-Gamma}
 \xymatrix{
 (\mc G', \psi') \ar[r]^-{\phi'} & 
 (\ta{\mc G'}, \ta \psi')   \\
  (\mc G, \psi)  \ar[r]^-{\phi}  \ar[u]_-{\kappa}  &
 (\ta{\mc G}, \ta \psi) \ar[u]_-{\ta \kappa}
 } \end{gather}
 such that   
 \begin{itemize}
 \item the lower line  $\big( (x_i)_{i\in I}, (\mc G, \psi) \xrightarrow{\phi}   (\ta{\mc G}, \ta \psi)
\big)$
and  the upper line   $\big( (x_i)_{i\in I}, (\mc G', \psi') \xrightarrow{\phi'}   (\ta{\mc G'}, \ta \psi')
\big)
$ belong to  $\mc Z^{(I)}$, 
\item  $\kappa:\restr{\mc G}{(X\sm v)\times S}\isom \restr{\mc G'}{(X\sm v)\times S}$  is an   isomorphism 
such that the   relative  position of $\mc G$ w.r.t. $\mc G'$ at   $v$ is {\it equal} to the  dominant weight  of $V$ (we recall that  $V$ is minuscule),
\item the restriction of $\kappa$ to $N\times S$,  which is an isomorphism,  intertwines the level structures  $\psi$ and $\psi'$.   
\end{itemize}
Moreover the two   projections 
 $\Gamma^{(I)}\to \mc Z^{(I)}$ are the  morphisms  which keep the lower and upper lines of   \eqref{intro-diag-Gamma}. 
  
 Since the  legs indexed by   $I$ vary in   $X\sm (N\cup v)$ and remain disjoint from the  legs $1$ and   $2$ fixed  at  $v$, we can replace  the partition $(\{1\},\{2\}, I)$ by  $(\{1\},I,\{2\})$ and therefore we have  
    \begin{gather*}\mc Z^{(\{1\},\{2\}, I)}
=
 \restr{\Cht_{N,I \cup \{1,2\},W \boxtimes V\boxtimes V^{*}}^{(\{1\},I,\{2\})}}
 {(X\sm (N\cup v))^{I}\times \Delta(v)}.\end{gather*}
 In other words  the stack  $\mc Z^{(\{1\},\{2\}, I)}$ classifies the data of  $(x_i)_{i\in I}$ and of a diagram 
  \begin{gather} \label{intro-diag-W}
 \xymatrix{
 & (\mc G_{1}, \psi_{1}) \ar[r]^-{\phi'_{2}} \ar[d]^-{\phi_{2}}& 
 (\mc G'_{2}, \psi'_{2}) \ar[d]^-{\phi'_{3}} & 
 (\ta{\mc G_{1}}, \ta \psi_{1}) \\
 (\mc G_{0}, \psi_{0}) \ar[ru]^-{\phi_{1}} &
 (\mc G_{2}, \psi_{2}) \ar[r]^-{\phi_{3}}    &
 (\ta{\mc G_{0}}, \ta \psi_{0})\ar[ru]^-{\ta \phi_{1}} &
 } \end{gather}
 with 
\begin{gather*}\big( (x_i)_{i\in I}, (\mc G_{0}, \psi_{0}) \xrightarrow{\phi_{1}}   (\mc G_{1}, \psi_{1}) \xrightarrow{\phi_{2}} (\mc G_{2}, \psi_{2}) \xrightarrow{\phi_{3}}    (\ta{\mc G_{0}}, \ta \psi_{0})
\big)\\
\in \restr{\Cht_{N,I \cup \{1,2\},W \boxtimes V\boxtimes V^{*}}^{(\{1\},\{2\},I)}}
 {(X\sm (N\cup v))^{I}\times \Delta(v)}\end{gather*} and  \begin{gather*}\big( (x_i)_{i\in I}, (\mc G_{0}, \psi_{0}) \xrightarrow{\phi_{1}}   (\mc G_{1}, \psi_{1}) \xrightarrow{\phi'_{2}} (\mc G'_{2}, \psi'_{2}) \xrightarrow{\phi'_{3}}    (\ta{\mc G_{0}}, \ta \psi_{0})
\big)\\
\in \restr{\Cht_{N,I \cup \{1,2\},W \boxtimes V\boxtimes V^{*}}^{(\{1\},I,\{2\})}}
 {(X\sm (N\cup v))^{I}\times \Delta(v)}. \end{gather*}  
 The   oblique, vertical  and horizontal  arrow  of  the diagram  \eqref{intro-diag-W} are respectively the   modifications associated to the  leg  $1$, to the leg $2$ and to the  legs indexed by   $I$. The arrow $\ta \phi_{1}$ on the right of  diagram \eqref{intro-diag-W} is determined by  $\phi_{1}$ , but we draw it because it will be used to define $\mc Y_{2}$ below.

 We denote by   $\mc Y_{1} \overset{\iota_{1}}{\hookrightarrow} 
 \mc Z^{(\{1\},\{2\}, I)}$ the closed substack  defined  by  the condition that in the diagram  \eqref{intro-diag-W},  $\phi_{2}\phi_{1}: 
 \restr{\mc G_{0}}{(X-v)\times S}\to \restr{\mc G_{2}}{(X-v)\times S}$ 
 extends to  an  isomorphism over  $X\times S$. 
We have a   morphism $$\alpha_{1}: \mc Y_{1}\to \mc Z^{(I)}
 $$
 which sends    \begin{gather*}
 \xymatrix{
 & (\mc G_{1}, \psi_{1}) \ar[r]^-{\phi'_{2}} \ar[d]^-{\phi_{2}}& 
 (\mc G'_{2}, \psi'_{2}) \ar[d]^-{\phi'_{3}} & 
 (\ta{\mc G_{1}}, \ta \psi_{1}) \\
 (\mc G_{0}, \psi_{0}) \ar[ru]^-{\phi_{1}} \ar[r]^-{\sim}&
 (\mc G_{2}, \psi_{2}) \ar[r]^-{\phi_{3}}    &
 (\ta{\mc G_{0}}, \ta \psi_{0})\ar[ru]^-{\ta \phi_{1}} &
 } \end{gather*}
to the lower line, i.e.    
 \begin{gather}\label{intro-ligne-bas}\big( (x_i)_{i\in I}, (\mc G_{0}, \psi_{0}) \xrightarrow{\phi_{3}(\phi_{2}\phi_{1})}    (\ta{\mc G_{0}}, \ta \psi_{0})
\big). \end{gather}

The assertion  B)  above comes from a similar  statement involving 
the  Mirkovic-Vilonen sheaves. Indeed  
\begin{itemize}\item by   a) of \thmref{satake-geom-thm} the direct image of 
 $\mc S_{\{1,2\},V\boxtimes V^{*} }^{(\{1\}, \{2\})}$ (which is the constant sheaf  $E$ with a cohomological shift) by the   morphism
 (which forgets the intermediate  modification) 
 $\mr{Gr}_{\{1,2\},V\boxtimes V^{*} }^{(\{1\}, \{2\})}\to 
 \mr{Gr}_{\{1,2\},V\boxtimes V^{*} }^{(\{1,2\})}$ is equal to  
 $\mc S_{\{1,2\},V\boxtimes V^{*} }^{(\{1,2\})}$, 
 \item  by   c) of \thmref{satake-geom-thm} the restriction of  $\mc S_{\{1,2\},V\boxtimes V^{*} }^{(\{1,2\})}$ over  the diagonal  
 (and {\it a fortiori}  over  $\Delta(v)$) is equal to  
  $\mc S_{\{0\},V\otimes V^{*} }^{(\{0\})}$ that we  send to 
  the skyscraper sheaf   $\mc S_{\{0\},\mbf 1}^{(\{0\})}$ by $\on{ev}_{V}:V\otimes V^{*}\to \mbf 1$
\end{itemize}
and  by the proper base change theorem this gives rise to a  cohomological correspondence between $\restr{\mr{Gr}_{\{1,2\},V\boxtimes V^{*} }^{(\{1\}, \{2\})}}{\Delta(v)}$ and the point, and one can check that it  is the  obvious cohomological correspondence   supported by the smooth closed subscheme   of $\restr{\mr{Gr}_{\{1,2\},V\boxtimes V^{*} }^{(\{1\}, \{2\})}}{\Delta(v)}$ consisting of 
$ (\mc G_{0} \xrightarrow{\phi_{1}}  
\mc G_{1}\xrightarrow{\phi_{2}}
  \mc G_{2}  \isom G ) $ such that $\phi_{2}\phi_{1}$ is an isomorphism. 
  
  We denote by   
$\mc Y_{2} \overset{\iota_{2}}{\hookrightarrow} \mc Z^{(\{1\},\{2\}, I)}$ 
the  closed substack   defined  by  the condition that   in  the diagram  \eqref{intro-diag-W},   $\ta \phi_{1}\phi_{3}': 
 \restr{\mc G'_{2}}{(X-v)\times S}\to \restr{\ta{\mc G_{1}}}{(X-v)\times S}$ extends to an  isomorphism over $X\times S$. 
 We have a   morphism $$\alpha_{2}: \mc Y_{2}\to \mc Z^{(I)}
 $$ 
 which sends  \begin{gather}\label{intro-diag-element-Y2}
 \xymatrix{
 & (\mc G_{1}, \psi_{1}) \ar[r]^-{\phi'_{2}} \ar[d]^-{\phi_{2}}& 
 (\mc G'_{2}, \psi'_{2}) \ar[d]^-{\phi'_{3}} \ar[r]^-{\sim}& 
 (\ta{\mc G_{1}}, \ta \psi_{1}) \\
 (\mc G_{0}, \psi_{0}) \ar[ru]^-{\phi_{1}} &
 (\mc G_{2}, \psi_{2}) \ar[r]^-{\phi_{3}}    &
 (\ta{\mc G_{0}}, \ta \psi_{0})\ar[ru]^-{\ta \phi_{1}} &
 } \end{gather}
to the upper line, i.e.    
 \begin{gather}\label{intro-ligne-haut}\big( (x_i)_{i\in I}, (\mc G_{1}, \psi_{1}) \xrightarrow{(\ta \phi_{1}\phi'_{3})\phi'_{2}}    (\ta{\mc G_{1}}, \ta \psi_{1})
\big). \end{gather}
The  justification of  assertion A)  above is given 
\begin{itemize}
\item by an argument similar to the argument used to justify  B)
but involving   $\delta_{V}:\mbf 1\to V\otimes V^{*}$ and the stack
$\restr{\Cht_{N,I \cup \{1,2\},W \boxtimes V\boxtimes V^{*}}^{(I,\{2\},\{1\})}}
 {(X\sm (N\cup v))^{I}\times \Delta(v)}$
\item by the fact that  the restriction to   $ (X\sm (N\cup v))^{I}\times \Delta(v)$ of the partial Frobenius morphism
$$\on {Fr}_{\{1\}} ^{(\{1\},I,\{2\})}:  \Cht_{N,I \cup \{1,2\},W \boxtimes V\boxtimes V^{*}}^{(\{1\},I,\{2\})} \to 
 \Cht_{N,I \cup \{1,2\},W \boxtimes V\boxtimes V^{*}}^{(I,\{2\},\{1\})}$$
sends \eqref{intro-diag-W} to 
\begin{gather}  \nonumber 
 \xymatrix{
   (\mc G_{1}, \psi_{1}) \ar[r]^-{\phi'_{2}}  & 
 (\mc G'_{2}, \psi'_{2}) \ar[d]^-{\phi'_{3}} & 
 (\ta{\mc G_{1}}, \ta \psi_{1}) \\
      &
 (\ta{\mc G_{0}}, \ta \psi_{0})\ar[ru]^-{\ta \phi_{1}} &
 } \end{gather}
\end{itemize}

 The fact that we do not need to introduce signs in A) and B) above is justified in remark  6.9 of \cite{coh}, whose idea is as follows. 
 It is enough to find a situation involving a creation operator  and an annihilation operator of the same type (with created or annihilated legs appearing in the same order), and where we know how to compute their composition directly. 
 But the composition  $V\xrightarrow{\Id_{V}\otimes \on{\delta}_{V}} 
    V\otimes V^{*} \otimes V\xrightarrow{ \on{ev}_{V} \otimes \Id_{V}} V$ is the identity by Zorro lemma  \eqref{zorro}. The composition of the cohomological correspondences between  Hecke stacks 
   given by the 
   creation  and  annihilation operators associated to 
    $\Id_{V}\otimes \on{\delta}_{V}$ et $ \on{ev}_{V} \otimes \Id_{V}$ 
    can be easily computed and it coincides with the identity correspondence when we choose signs   as in  A) and  B) above, which proves that these sign choices were good (or at least that the product of the two normalizations was good, the normalization of each one is more subtle and useless for us).

On the other hand  we have a  canonical  isomorphism
\begin{gather}\label{intro-Gamma-produit-fibre}
\Gamma^{(I)}\simeq \mc Y_{1} \times_{\mc Z^{(\{1\},\{2\}, I)}}\mc Y_{2}. 
\end{gather}
Indeed  a  point  of  $\mc Z^{(\{1\},\{2\}, I)}$ belonging to   $\mc Y_{1} $ and  to  $\mc Y_{2}$ is given  by a  diagram 
  $$ \xymatrix{
 & (\mc G_{1}, \psi_{1}) \ar[r]^-{\phi'_{2}} \ar[d]^-{\phi_{2}}& 
 (\mc G'_{2}, \psi'_{2}) \ar[d]^-{\phi'_{3}} \ar[r]^-{\sim}& 
 (\ta{\mc G_{1}}, \ta \psi_{1}) \\
 (\mc G_{0}, \psi_{0}) \ar[ru]^-{\phi_{1}} \ar[r]^-{\sim}&
 (\mc G_{2}, \psi_{2}) \ar[r]^-{\phi_{3}}    &
 (\ta{\mc G_{0}}, \ta \psi_{0})\ar[ru]^-{\ta \phi_{1}} &
 } $$
Therefore it is nothing but a point    of  $ \Gamma^{(I)}$, because the contraction of the  two   isomorphisms in the previous  diagram produces    the diagram 
  \begin{gather*}
 \begin{CD} 
  (\mc G_{1}, \psi_{1}) 
  @>(\ta \phi_{1}\phi'_{3}) \phi'_{2}>> 
  (\ta{\mc G_{1}}, \ta \psi_{1})  \\
 @AA\phi_{1}A 
 @AA\ta \phi_{1}A \\
(\mc G_{0}, \psi_{0})
 @>\phi_{3}(\phi_{2}\phi_{1})>>  
 (\ta{\mc G_{0}}, \ta \psi_{0}) 
 \end{CD}
 \end{gather*}
that we identify to the  diagram \eqref{intro-diag-Gamma}. 

We have  natural morphisms from  the stacks  $\mc Z^{(\{1\},\{2\}, I)}, \mc Z^{(I)} , \mc Y_{1}, \mc Y_{2}$  and $\Gamma^{(I)}$ to  $ \mr{Gr}_{I,W}^{(I)}/
  G_{\sum n_{i} x_i}$. Since $V$ is minuscule 
  these morphisms are smooth. Therefore the open smoothness loci 
    ${}^{\circ}{\mc Z}^{(\{1\},\{2\}, I)}, {}^{\circ}{\mc Z}^{(I)},  {}^{\circ}{\mc Y}_{1}, {}^{\circ}{\mc Y}_{2}, {}^{\circ}{\Gamma}^{(I)}$
are the  inverse images of $ {}^{\circ}{\mr{Gr}}_{I,W}^{(I)}/
  G_{\sum n_{i} x_i}$ where $ {}^{\circ}{\mr{Gr}}_{I,W}^{(I)}$ denotes  the open smoothness locus  of 
$ \mr{Gr}_{I,W}^{(I)}$. 

A computation of tangent spaces  shows that 
 ${}^{\circ}\mc Y_{1}$ and ${}^{\circ}\mc Y_{2}$ are   smooth substacks in the 
  smooth Deligne-Mumford stack  ${}^{\circ}\mc Z^{(\{1\},\{2\}, I)}$ and  that their intersection is transverse, and moreover      \eqref{intro-Gamma-produit-fibre} implies that their  intersection is equal to 
${}^{\circ}\Gamma^{(I)}$.   
Thus we obtained an  equality of     cohomological correspondences
between  $S_{V,v}$ and $T(h_{V,v})$ over ${}^{\circ}\Gamma^{(I)}$ but since 
$\Gamma^{(I)}$ is an   étale  correspondence  between 
$\mc Z^{( I)}$ and itself the equality holds everywhere 
(indeed a morphism from the  perverse sheaf $\on{IC}_{\Gamma^{(I)}}$ to itself  is determined by its  restriction to ${}^{\circ}{\Gamma}^{(I)}$). 
\cqfd

    A  consequence of    \propref{prop-coal-frob-cas-part-intro} is that we have   for  every    $f\in C_{c}(K_{N}\backslash G(\mb A)/K_{N},E)$ and  $\kappa$ big enough,  
    a natural extension   of  the morphism $T(f)$ (introduced in  \eqref{defi-Tf})  to  a  morphism  
 $T(f):\mc H_{N,I,W}^{\leq\mu}\to 
 \mc H_{N,I,W}^{\leq\mu+\kappa}$ of  constructible sheaves over the whole  $  (X\sm N)^{I}$, in a way  compatible with the  composition of   Hecke operators. Indeed, if we write $K_{N}=\prod K_{N,v}$,  it is enough to do it  for each  place $v$ and for  every  $f\in C_{c}(K_{N,v}\backslash G(F_{v})/K_{N,v},E)$.
  There is nothing to do   if $v\in N$. If $v\not\in N$ 
  it is enough to consider the case  where $f=h_{V,v}$ and then 
  the extension is given by $S_{V,v}$ thanks to   \propref{prop-coal-frob-cas-part-intro}. For further details, we refer to  corollary 6.5 
  of \cite{coh}.

 For the  Shimura varieties over number fields such extensions  were defined  in many cases, in a modular way, by   Zariski closure or   with the help  of  nearby cycles (see \cite{deligne-bki-mod,faltings-chai, genestier-tilouine}). 

     Since $S_{V,v}$ is the extension  of 
   $T(h_{V,v})$, the following proposition  is exactly the 
  Eichler-Shimura relation . We use again $\{0\}$ to denote a set with one element  (indexing the leg to which the Eichler-Shimura relation  applies). 
 
 \begin{prop}\label{Eichler-Shimura-intro}
 Let $I,W$ be as above and $V$ be  an   irreducible representation of  $\wh G$. Then  $$F_{\{0\}}^{\deg(v)}: \varinjlim _{\mu}\restr{\mc H _{N, I\cup\{0\}, W\boxtimes V}^{\leq\mu}}{(X\sm N)^{I}\times v}\to \varinjlim _{\mu}\restr{\mc H _{N, I\cup\{0\}, W\boxtimes V}^{\leq\mu}}{(X\sm N)^{I}\times v}$$ is 
 killed by a polynomial of degree 
 $\dim(V)$ whose  coefficients are the 
  restrictions to  $(X\sm N)^{I}\times v$ of the  morphisms $S_{\Lambda^{i}V,v}$. 
More precisely we have  
  \begin{gather*}  \sum_{i=0}^{\dim V} (-1)^{i} (F_{\{0\}}^{\deg(v)})^{i}\circ \restr{S_{\Lambda^{\dim V-i}V,v}}{(X\sm N)^{I}\times v}=0  . \end{gather*} 
    \end{prop}
    We recall that  $S_{\Lambda^{i}V,v}$ extends  the Hecke operator 
 $T(h_{\Lambda^{i}V,,v})$ $$\text{from \ \ }  (X\sm (N\cup v))^{I\cup \{0\}}\text{ \ \  to \ \ }(X\sm N)^{I\cup \{0\}}$$ and we note that this  extension is absolutely necessary in order to take its  restriction to  $(X\sm N)^{I}\times v$. Thanks to  the definition of the  morphisms   $S_{\Lambda^{i}V,v}$ by  \eqref{def-SVv-intro1}-\eqref{def-SVv-intro4}, 
 the proof of   \propref{Eichler-Shimura-intro} is a simple computation of tensor algebra (inspired by a  proof  of the  Hamilton-Cayley theorem, and based uniquely on the fact that  $\Lambda^{\dim V+1}V=0$). 
 We refer to  chapter 7 of \cite{coh} for this proof. 
 
 \begin{rem} In \cite{xiao-zhu}, Liang Xiao and  Xinwen Zhu have defined
 (in a slightly different setting) a ring of cohomological  correspondences   
 between  $\restr{ \Cht_{ N,I\cup \{0\},W\boxtimes V}^{(I,\{0\})}}{(X\sm N)^{I}\times v}$ and itself, 
 in which the Eichler-Shimura relation results formally from the Hamilton-Cayley equality. 
\end{rem}
 
\section{ Constructible subsheaves stable under the action of the   partial Frobenius morphisms}
\label{section-sous-sheaves-constr}
The goal of this  section is to prove   \lemref{lem-Hf-union-stab} and    \propref{cor-action-Hf-intro},  which had been admitted  in   section \ref{section-esquisse-abc}. We refer to chapter 8 of \cite{coh} for more details.

  We recall that the we have fixed a  geometric point  $\ov\eta=\on{Spec}(\ov F)$  over  the  generic point  $\eta$ of $X$.   
   Let $I$ be a finite  set and  $W=\boxtimes_{i\in I}W_{i}$ be an   irreducible representation of $(\wh G)^{I}$. 
 We denote by   $\Delta:X\to X^{I}$ the  diagonal morphism. 
 We recall that we have fixed a  geometric point $\ov{\eta^{I}}$ over the   generic point $\eta^{I}$ of $X^{I}$ and a specialization arrow   
 $\on{\mf{sp}}: \ov{\eta^{I}}\to \Delta(\ov \eta)$.

 \begin{lem}\label{lem-Hf-union-stab}
 The space    $\Big( \varinjlim _{\mu}\restr{\mc H _{N, I, W}^{\leq\mu}}{\ov{\eta^{I}}}\Big)^{\mr{Hf}}$  is the union  of $\mc O_{E}$-submodules $\mf M=\restr{\mf G}{\ov{\eta^{I}}}$ where 
   $\mf G$ is a constructible   $\mc O_{E}$-subsheaf  of 
  $\varinjlim _{\mu}\restr{\mc H _{N, I, W}^{\leq\mu}}{\eta^{I}}$ 
     stable under the action of  the partial Frobenius morphisms. 
 \end{lem}
  \noindent {\bf Proof.} We refer to  the proof of    proposition 8.27 of \cite{coh} for more details. 
For every family  $(v_{i})_{i\in I}$ of closed points   of  $X\sm N$, we denote   ${\times_{i\in I} v_{i}}$ their product,  which is a finite union    of  closed points   of  $(X\sm N)^{I}$. 
  Let $\check{\mf M}$ be a  $\mc O_{E}$-submodule of finite type of 
  $\varinjlim _{\mu}\restr{\mc H _{N, I, W}^{\leq\mu}}{\ov{\eta^{I}}}$ stable by  
$\pi_{1}(\eta^{I},\ov{\eta^{I}})$ and   $C_{c}(K_{N}\backslash G(\mb A)/K_{N},\mc O_{E})$. 
We will  construct $\mf M\supset \check{\mf M}$ satisfying the properties of the statement of the lemma. 
Since  $\check{\mf M}$ is of finite type, we can find  $\mu_0$ such that  $\check{\mf M}$ is  included in  the image of  $\restr{\mc H _{N, I, W}^{\leq\mu_0}}{\ov{\eta^{I}}}$ in  $\varinjlim _{\mu}\restr{\mc H _{N, I, W}^{\leq\mu}}{\ov{\eta^{I}}}$. Increasing   $\mu_0$ if necessary,  we can assume that   $\check{\mf M}$ is a $\mc O_{E}$-submodule  of  $\restr{\mc H _{N, I, W}^{\leq\mu_0}}{\ov{\eta^{I}}}$. 
Let $\Omega_0$ be a dense open  subscheme of  $X^{I}$ over which   
$ \mc H _{N, I, W}^{\leq\mu_0}$ is smooth. We have a unique smooth 
$\mc O_{E}$-subsheaf  $\check{\mf G}\subset  \restr{\mc H _{N, I, W}^{\leq\mu_0}}{\Omega_0}$ over  $\Omega_0$ such that  $\restr{\check{\mf G}}{\ov{\eta^{I}}}=\check{\mf M}$. We choose  
 $(v_{i})_{i\in I}$ such that  ${\times_{i\in I} v_{i}}$ is included  in  $\Omega_0$. For every $i$, the  Eichler-Shimura relation   
 (\propref{Eichler-Shimura-intro}) 
 implies that  \begin{gather}\label{intro-ES-inclusion}(F_{\{i\}}^{\deg(v_{i})})^{\dim W_{i}}(\restr{\check{\mf G}}{\times_{i\in I} v_{i}})
   \subset 
   \sum_{\alpha=0}^{\dim W_{i}-1}  (F_{\{i\}}^{\deg(v_{i})})^{\alpha}(S_{\Lambda^{\dim W_{i}-\alpha}W_{i},v_{i}}(\restr{\check{\mf G}}{\times_{i\in I} v_{i}}))
  \end{gather} in  $
  \varinjlim_{\mu} \restr{\mc H _{N, I, W}^{\leq\mu}}{\times_{i\in I} v_{i}}$. 
Thanks to the smoothness  of 
$(\Frob_{\{i\}}^{\deg(v_{i})\dim W_{i}})^{*}(\check{\mf G})$
at ${\times_{i\in I} v_{i}}$, the inclusion   \eqref{intro-ES-inclusion} 
propagates to $\eta^{I}$, i.e.    
 \begin{gather*}
 F_{\{i\}}^{\deg(v_{i})\dim W_{i}}
 (\restr{(\Frob_{\{i\}}^{\deg(v_{i})\dim W_{i}})^{*}(\check{\mf G}}{\eta^{I}}) )
 \\   \subset 
   \sum_{\alpha=0}^{\dim W_{i}-1}  F_{\{i\}}^{\deg(v_{i})\alpha}(\Frob_{\{i\}}^{\deg(v_{i})\alpha})^{*}(\restr{S_{\Lambda^{\dim W_{i}-\alpha}W_{i},v_i}(\check{\mf G})}{\eta^{I}})
  \end{gather*} in  $
  \varinjlim_{\mu} \restr{\mc H _{N, I, W}^{\leq\mu}}{\eta^{I}}$. 
But   $\restr{\check{\mf G}}{\eta^{I}}$ is stable by  $S_{\Lambda^{\dim W_{i}-\alpha}W_{i},v_i}=T(h_{\Lambda^{\dim W_{i}-\alpha}W_{i},v_i})$ 
 since  
 $$h_{\Lambda^{\dim W_{i}-\alpha}W_{i},v_i}\in C_{c}(G(\mc O_{v_{i}})\backslash G(F_{v_{i}})/G(\mc O_{v_{i}}),\mc O_{E})\subset C_{c}(K_{N}\backslash G(\mb A)/K_{N},\mc O_{E}). $$
  Consequently 
 \begin{gather*}
 F_{\{i\}}^{\deg(v_{i})\dim W_{i}}
 (\restr{(\Frob_{\{i\}}^{\deg(v_{i})\dim W_{i}})^{*}(\check{\mf G}}{\eta^{I}}) )
   \subset 
   \sum_{\alpha=0}^{\dim W_{i}-1}  F_{\{i\}}^{\deg(v_{i})\alpha}(\Frob_{\{i\}}^{\deg(v_{i})\alpha})^{*}(\restr{\check{\mf G}}{\eta^{I}})
  \end{gather*} in  $
  \varinjlim_{\mu} \restr{\mc H _{N, I, W}^{\leq\mu}}{\eta^{I}}$. 
  Therefore   $$\mf G=\sum_{(n_{i})_{i\in I}\in \prod _{i\in I}\{0,...,\deg(v_{i})\dim(W_{i})-1\}}\restr{\prod _{i\in I}F_{\{i\}}^{n_{i}}\Big(\prod _{i\in I}\Frob_{\{i\}}^{n_{i}}\Big)^{*}(\check{\mf G})}{\eta^{I}}$$
 is  a  constructible $\mc O_{E}$-subsheaf   of 
  $\varinjlim _{\mu}\restr{\mc H _{N, I, W}^{\leq\mu}}{\eta^{I}}$ 
 which  is stable under  the action of the partial Frobenius morphisms. Since 
$\Big( \varinjlim _{\mu}\restr{\mc H _{N, I, W}^{\leq\mu}}{\ov{\eta^{I}}}\Big)^{\mr{Hf}}$  is the union  of $\mc O_{E}$-submodules $\check{\mf M}$ as at the beginning  of the proof and since $ \mf M =\restr{ \mf G }{\ov{\eta^{I}}}$ contains $\check{\mf M}$ we get the statement of the lemma. 
 \cqfd

\begin{prop}\label{cor-action-Hf-intro}  
The space $\Big( \varinjlim _{\mu}\restr{\mc H _{N, I, W}^{\leq\mu}}{\ov{\eta^{I}}}\Big)^{\mr{Hf}}$ is equipped with a  natural action   of  $\pi_{1}(\eta,\ov{\eta})^{I}$. More precisely it is a union of  $E$-vector subspaces   equipped with a  continuous  action of 
$\pi_{1}(\eta,\ov{\eta})^{I}$. 
 \end{prop}
\dem For every constructible 
   $\mc O_{E}$-subsheaf  $\mf G$   of 
  $\varinjlim _{\mu}\restr{\mc H _{N, I, W}^{\leq\mu}}{\eta^{I}}$ 
     stable under the action of  partial Frobenius morphisms, 
       Drinfeld's lemma   \ref{lem-Dr-intro} 
 provides  (thanks to  $\mf{sp}$ and \remref{rem-lem-Drinfeld}) a continuous action  of   $\pi_{1}(\eta,\ov{\eta})^{I}$ on 
$\mf M=\restr{\mf G}{\ov{\eta^{I}}}$. By \lemref{lem-Hf-union-stab}, 
  $\Big( \varinjlim _{\mu}\restr{\mc H _{N, I, W}^{\leq\mu}}{\ov{\eta^{I}}}\Big)^{\mr{Hf}}$ is the union of  such  $\mf M$. \cqfd

\begin{rem}\label{rem-WeilF} The action of $\pi_{1}(\eta,\ov{\eta})^{I}$
on $\Big( \varinjlim _{\mu}\restr{\mc H _{N, I, W}^{\leq\mu}}{\ov{\eta^{I}}}\Big)^{\mr{Hf}}$  is   uniquely  determined  by the actions of 
$\pi_{1}(\eta^{I},\ov{\eta^{I}})$ and   of the  partial Frobenius morphisms. This is a consequence  of \lemref{lem-Dr-intro}
but we can see it also in the following way (for more details  we refer to  chapter 8 of \cite{coh}). 
 Following  \cite{drinfeld78}, we shall define a group    $\on{FWeil}(\eta^{I},\ov{\eta^{I}})$ 
       \begin{itemize}
     \item  which   is a  extension of $\Z^{I}$ by   $\on{Ker}(\pi_{1}(\eta^{I},\ov{\eta^{I}})\to \wh \Z)$, 
     \item and which coincides,  when  $I$ is a   singleton, 
 with the usual Weil group 
   $\on{Weil} (\eta,\ov{\eta})=\pi_{1}(\eta,\ov{\eta})\times_{\wh \Z}\Z$.      
       \end{itemize}
       
   We denote by   $F^{I}$ the field of functions  of $X^{I}$, by $(F^{I})^{\mr{perf}}$ its perfectization  and by $\ov{F^{I}}$ the algebraic closure of $F^{I}$ such that 
  $\ov{\eta^{I}}=\on{Spec}(\ov{F^{I}})$. 
   Then we define   
          \begin{gather*}\on{FWeil}(\eta^{I},\ov{\eta^{I}})=
       \big\{\varepsilon \in \on{Aut}_{\ov\Fq}((\ov{F^{I}})), \exists (n_{i})_{i\in I}\in \Z^{I}, \restr{\varepsilon}{(F^{I})^{\mr{perf}}}=\prod_{i\in I}(\Frob_{\{i\}})^{n_{i}}\big\} 
  .\end{gather*}
The choice of  
  $\mf{sp}$ provides 
  a inclusion   $\ov{F}\otimes_{\ov\Fq} \cdots 
 \otimes_{\ov\Fq}\ov{F}\subset \ov{F^{I}} 
  $.  By restriction of the automorphisms, this gives a surjective
   morphism   
    \begin{gather}\label{morph-Weil-I}\on{FWeil}(\eta^{I},\ov{\eta^{I}})\to \big(\on{Weil} (\eta,\ov{\eta})\big)^{I}.\end{gather} The statement of  \propref{cor-action-Hf-intro}  can be reformulated by saying that the natural action  of $\on{FWeil}(\eta^{I},\ov{\eta^{I}})$ 
    on  $\Big( \varinjlim _{\mu}\restr{\mc H _{N, I, W}^{\leq\mu}}{\ov{\eta^{I}}}\Big)^{\mr{Hf}}$
    factorizes through the morphism  \eqref{morph-Weil-I}, and even through 
    $(\pi_{1}(\eta,\ov{\eta}))^{I}$.  
      \end{rem}

  \section{Specialization homorphisms   and  Hecke-finite cohomology} 
  \label{para-homom-spe-Hf}
   The goal of this  section is  to prove \corref{bijectivite-Hecke-fini},  which had been admitted in  section \ref{section-esquisse-abc}. Let $W=\boxtimes_{i\in I} W_{i}$ be an   irreducible  $E$-linear representation of $(\wh G)^{I}$.

  \begin{prop} \label{surjectivite-Hecke-fini-intro}
   The image of the specialization homomorphism  
  \begin{gather}\label{intro-sp}\on{\mf{sp}}^{*}: \varinjlim _{\mu}\restr{\mc H _{N, I, W}^{\leq\mu}}{\Delta(\ov{\eta})}\to 
  \varinjlim _{\mu}\restr{\mc H _{N, I, W}^{\leq\mu}}{\ov{\eta^{I}}}\end{gather}
 contains 
  $\Big( \varinjlim _{\mu}\restr{\mc H _{N, I, W}^{\leq\mu}}{\ov{\eta^{I}}}\Big)^{\mr{Hf}}$. 
  \end{prop}
   \noindent {\bf Proof.} We refer to  the proof of   proposition 8.31   of \cite{coh} for more details. 
   By  \lemref{lem-Hf-union-stab}, 
    $\Big( \varinjlim _{\mu}\restr{\mc H _{N, I, W}^{\leq\mu}}{\ov{\eta^{I}}}\Big)^{\mr{Hf}}$  is the union  of $\mc O_{E}$-submodules $\mf M=\restr{\mf G}{\ov{\eta^{I}}}$ where 
   $\mf G$ is a constructible  $\mc O_{E}$-subsheaf   of 
  $\varinjlim _{\mu}\restr{\mc H _{N, I, W}^{\leq\mu}}{\eta^{I}}$ 
     stable under the action of  the partial Frobenius morphisms. Therefore it is enough to prove that such a  $\mf M$ is included in the image of  \eqref{intro-sp}.  
  Let $\mu_0$ be big enough  so that   $\mf G$ is  a  $\mc O_{E}$-subsheaf  of    
  $\restr{\mc H _{N, I, W}^{\leq\mu_0}}{\eta^{I}}$.  Let $\Omega_0$ be a  dense  open subscheme of  $X^{I}$ such that   $\restr{\mc H _{N, I, W}^{\leq\mu_0}}{\Omega_0}$ is smooth. Then    $\mf G$ extends to  a  smooth $\mc O_{E}$-subsheaf  of 
  $\restr{\mc H _{N, I, W}^{\leq\mu_0}}{\Omega_0}$. 
    By the proof of lemma 9.2.1 of \cite{eike-lau}, 
     the set of the    $\prod_{i\in I}\Frob_{\{i\}}^{n_{i}}(\Delta(\eta))$ for $(n_{i})_{i\in I}\in \N^{I}$ is Zariski dense in $X^{I}$. Therefore we can find  $(n_{i})_{i\in I}\in \N^{I}$ such that 
  $\prod_{i\in I}\Frob_{\{i\}}^{n_{i}}(\Delta(\eta))$ belongs to   $\Omega_0$. 
 Then  
  $\restr{\mf G}{\prod_{i\in I}\Frob_{\{i\}}^{n_{i}}(\ov{\eta^{I}})}$ is included  in the image of    \begin{gather}\label{intro-sp-Frob}\on{\wt{\mf{sp}}}^{*}: \varinjlim _{\mu}\restr{\mc H _{N, I, W}^{\leq\mu}}{\prod_{i\in I}\Frob_{\{i\}}^{n_{i}}(\Delta(\ov{\eta}))}\to 
  \varinjlim _{\mu}\restr{\mc H _{N, I, W}^{\leq\mu}}{\prod_{i\in I}\Frob_{\{i\}}^{n_{i}}(\ov{\eta^{I}})}\end{gather}
 for every  specialization arrow $\wt{\mf{sp}}: \prod_{i\in I}\Frob_{\{i\}}^{n_{i}}(\ov{\eta^{I}})\to \prod_{i\in I}\Frob_{\{i\}}^{n_{i}}(\Delta(\ov{\eta}))$, and in particular  for the image of  $\mf{sp}$ by $\prod_{i\in I}\Frob_{\{i\}}^{n_{i}}$. Since   $\mf G$ is stable under the action  of the   partial Frobenius morphisms, we conclude that  
  $\mf M=\restr{\mf G}{\ov{\eta^{I}}}$ is included in the image of  \eqref{intro-sp}. \cqfd

\begin{prop}\label{injectivite-sp}
The specialization homomorphism 
 \begin{gather}\label{sp*-sans-Hf-intro2}\on{\mf{sp}}^{*}: 
 \varinjlim _{\mu}\restr{\mc H _{N, I, W}^{\leq\mu}}{\Delta(\ov{\eta})} \to
  \varinjlim _{\mu}\restr{\mc H _{N, I, W}^{\leq\mu}}{\ov{\eta^{I}}}\end{gather} is injective.   \end{prop}
\noindent{\bf Proof. }
Let $a$ be  in the kernel of \eqref{sp*-sans-Hf-intro2}. 
We choose $\mu_{0}$ and $\wt a\in \restr{\mc H _{N, I, W}^{\leq\mu_{0}}}{\Delta(\ov{\eta})}$ such that 
 $a$ is   the image of $\wt a$ in 
 $\varinjlim _{\mu}\restr{\mc H _{N, I, W}^{\leq\mu}}{\Delta(\ov{\eta})} $. 
Let $\Omega_{0}$ be a dense open subscheme  of $X\sm N$ over which  
$\Delta^{*}\big( \mc H _{N, I, W}^{\leq\mu_{0}}\big)$ is smooth. Let $v\in |\Omega_{0}|$. We set  $d=\deg(v)$ to shorten the formulas. Let $\ov v$ be a geometric point over  $v$. 
Let $\on{\mf{sp}}_{v}:\ov \eta\to \ov v$ be a specialization arrow. We still denote by     $\on{\mf{sp}}_{v}:\Delta(\ov \eta)\to \Delta(\ov v)$ the specialization arrow  equal to its image by $\Delta$. 
By the smoothness  of $\Delta^{*}\big( \mc H _{N, I, W}^{\leq\mu_{0}}\big)$  over 
 $\Omega_{0}$ we get  a unique element 
 $\wt b\in \restr{\mc H _{N, I, W}^{\leq\mu_{0}}}{\Delta(\ov{v})}$ such that $\wt a =\on{\mf{sp}}_{v}^{*}(\wt b)$. 
 We denote by  $b$ the image of  $\wt b$ in $\varinjlim _{\mu}\restr{\mc H _{N, I, W}^{\leq\mu}}{\Delta(\ov{v})}$, so that  $a$ is the image of $b$ by 
 $$\on{\mf{sp}}_{v}^{*}: \varinjlim _{\mu}\restr{\mc H _{N, I, W}^{\leq\mu}}{\Delta(\ov{v})}\to 
  \varinjlim _{\mu}\restr{\mc H _{N, I, W}^{\leq\mu}}{\Delta(\ov{\eta})} .$$
  
The action of the  partial Frobenius morphisms provide  for every $\mu$  and for every 
$(n_{i})_{i\in I}\in \N^{I}$  
a morphism of sheaves over $(X\sm N)^{I}$ 
\begin{gather}\label{action-Frob-partiels-proof}\prod_{i\in I}F_{\{i\}}^{dn_{i}}:(\prod_{i\in I}\Frob_{\{i\}}^{dn_{i}})^{*} ( \mc H _{N, I, W}^{\leq\mu})\to 
 \mc H _{N, I, W}^{\leq\mu+\kappa(\sum n_{i})}\end{gather} 
with  $\kappa $ big enough in function of $W$ and   $d$.   
Since $ \prod_{i\in I}\Frob_{\{i\}}^{dn_{i}} $ acts trivialy on $\Delta(v)$, 
 $\prod_{i\in I}F_{\{i\}}^{dn_{i}}$ acts on 
$\varinjlim _{\mu}\restr{\mc H _{N, I, W}^{\leq\mu}}{\Delta(\ov{v})}$. 

For every $(n_{i})_{i\in I}\in \N^{I}$ we write 
$$b_{(n_{i})_{i\in I}} =\prod_{i\in I}F_{\{i\}}^{dn_{i}}(b)\in 
\varinjlim _{\mu}\restr{\mc H _{N, I, W}^{\leq\mu}}{\Delta(\ov{v})}. $$
In particular  
$b_{(0)_{i\in I}}=b$. 

We set 
\begin{gather}\label{def-a-n-i}a_{(n_{i})_{i\in I}} =\on{\mf{sp}}_{v}^{*}(b_{(n_{i})_{i\in I}} )\in 
\varinjlim _{\mu}\restr{\mc H _{N, I, W}^{\leq\mu}}{\Delta(\ov{\eta})}, \end{gather}
and we note that   $a_{(0)_{i\in I}}=a$. 

The  sequence $a_{(n_{i})_{i\in I}}$ satisfies  the two properties stated in the following lemma. The first one claims that this  sequence is ``multirecurrent'', i.e.  recurrent in each  variable $n_{i}$, and the second one implies it vanishes ``almost everywhere''. We will deduce easily from  the conjonction of these  two properties that this  sequence vanishes everywhere, and   in particular  that 
$a=a_{(0)_{i\in I}}$ is zero. 

To state the second property we note that   $\on{\mf{sp}}^{*}(\on{\mf{sp}}_{v}^{*}(b))=\on{\mf{sp}}^{*}(a)=0$  
 in $\varinjlim _{\mu}\restr{\mc H _{N, I, W}^{\leq\mu}}{\ov{\eta^{I}}}
$. 
Therefore we can find  $\mu_{1}\geq \mu_{0}$ such that 
$\on{\mf{sp}}^{*}(\on{\mf{sp}}_{v}^{*}(\wt b))\in \restr{\mc H _{N, I, W}^{\leq\mu_{0}}}{\ov{\eta^{I}}}$
has a zero image  in $\restr{\mc H _{N, I, W}^{\leq\mu_{1}}}{\ov{\eta^{I}}}
$. In other words if we denote by 
$\wh b$ the image of $\wt b$ in 
$\restr{\mc H _{N, I, W}^{\leq\mu_{1}}}{\Delta(\ov{v})}$, 
we have 
$\on{\mf{sp}}^{*}(\on{\mf{sp}}_{v}^{*}(\wh b))=0$ in 
$\restr{\mc H _{N, I, W}^{\leq\mu_{1}}}{\ov{\eta^{I}}}$. 
Let $\Omega_{1}\subset (X\sm N)^{I}$ be a dense open subscheme over which $\mc H _{N, I, W}^{\leq\mu_{1}}$ is smooth. 

\begin{lem}\label{lem-a-b}
a) For every $j\in I$ and for every $(n_{i})_{i\in I}\in \N^{I}$, 
\begin{gather} \label{relation-ani}
\sum_{\alpha=0}^{\dim W_{j}}(-1)^{\alpha} S_{\Lambda^{\dim W_{j}-\alpha}W_{j},v} (a_{(n_{i}+\alpha \delta_{i,j})_{i\in I}})=0
\end{gather}
in $\varinjlim _{\mu}\restr{\mc H _{N, I, W}^{\leq\mu}}{\Delta(\ov{\eta})}$. 

b)  
For every $(n_{i})_{i\in I}\in \N^{I}$ such that $\prod_{i\in I}\Frob_{\{i\}}^{dn_{i}}
(\Delta(\ov\eta))\in \Omega_{1}$, we have $a_{(n_{i})_{i\in I}}=0$ in 
$\varinjlim _{\mu}\restr{\mc H _{N, I, W}^{\leq\mu}}{\Delta(\ov{\eta})} $. 
\end{lem}

\noindent{\bf Proof of a). }
The $b_{(n_{i})_{i\in I}} $ satisfy a relation identical to 
\eqref{relation-ani} 
(in $\varinjlim _{\mu}\restr{\mc H _{N, I, W}^{\leq\mu}}{\Delta(\ov{v})}$), namely  the  Eichler-Shimura relation  for the leg  $j$ (\propref{Eichler-Shimura-intro}). 
Then \eqref{relation-ani} is obtained by applying 
$\on{\mf{sp}}_{v}^{*}$ to this  relation (we are allowed to do this because  the  $S_{\Lambda^{\dim W_{j}-\alpha}W_{j},v}$ are   morphisms of sheaves).  \cqfd

\noindent{\bf Proof of b). } Let $(n_{i})_{i\in I}$ satisfy the hypothesis  of b). 
Since 
\eqref{action-Frob-partiels-proof} 
is  a morphism of sheaves over $(X\sm N)^{I}$, 
we can inverse the order of  the specialization homomorphisms and the partial Frobenius morphisms. 
In other words 
we have a commutative diagram 
  \begin{gather*} \xymatrixcolsep{5pc} \xymatrix{
     \restr{ \mc H _{N, I, W}^{\leq\mu_{1}}}{\Delta(\ov v)}
=
\restr{(\prod_{i\in I}\Frob_{\{i\}}^{dn_{i}})^{*} ( \mc H _{N, I, W}^{\leq\mu_{1}})}{\Delta(\ov v)}\ar[d]^-{ \on{\mf{sp}}_{v,(n_{i})_{i\in I}}^{*} } \ar[r]^-{\prod_{i\in I}F_{\{i\}}^{dn_{i}}}
 & \varinjlim _{\mu}\restr{\mc H _{N, I, W}^{\leq\mu}}{\Delta(\ov{v})} \ar[d]^-{\on{\mf{sp}}_{v}^{*}} 
 \\
\restr{(\prod_{i\in I}\Frob_{\{i\}}^{dn_{i}})^{*} ( \mc H _{N, I, W}^{\leq\mu_{1}})}{\Delta(\ov \eta)}
 \ar[r]^-{\prod_{i\in I}F_{\{i\}}^{dn_{i}}}
& \varinjlim _{\mu}\restr{\mc H _{N, I, W}^{\leq\mu}}{\Delta(\ov{\eta})}  }\end{gather*}
 where the notation $ \on{\mf{sp}}_{v,(n_{i})_{i\in I}}^{*}$ indicates that the specialization homomorphism associated to the arrow $\on{\mf{sp}}_{v}:\Delta(\ov \eta)\to \Delta(\ov v)$ {\it is applied to the sheaf } $(\prod_{i\in I}\Frob_{\{i\}}^{dn_{i}})^{*} ( \mc H _{N, I, W}^{\leq\mu_{1}})$ (and not to  
 $\mc H _{N, I, W}^{\leq\mu_{1}}$). 
 The previous diagram gives rise to 
    \begin{gather*} \xymatrixcolsep{5pc} \xymatrix{
\wh b \ar@{|->}[d]^-{\on{\mf{sp}}_{v,(n_{i})_{i\in I}}^{*}} \ar@{|->}[r]^-{\prod_{i\in I}F_{\{i\}}^{dn_{i}}}
 &b_{(n_{i})_{i\in I}} \ar@{|->}[d]^-{\on{\mf{sp}}_{v}^{*}} 
 \\
\on{\mf{sp}}_{v,(n_{i})_{i\in I}}^{*}(\wh b) \ar@{|->}[r]^-{\prod_{i\in I}F_{\{i\}}^{dn_{i}}}
&a_{(n_{i})_{i\in I}}  }\end{gather*}

Therefore to prove $a_{(n_{i})_{i\in I}} =0$ (and finish the proof of b))it suffices to prove that \begin{gather}\label{fibre-Frob-Delta}
  \on{\mf{sp}}_{v,(n_{i})_{i\in I}}^{*}(\wh b) \in 
\restr{ \mc H _{N, I, W}^{\leq\mu_{1}}}{(\prod_{i\in I}\Frob_{\{i\}}^{dn_{i}})(\Delta(\ov \eta))}=
\restr{(\prod_{i\in I}\Frob_{\{i\}}^{dn_{i}})^{*} ( \mc H _{N, I, W}^{\leq\mu_{1}})}{\Delta(\ov \eta)}\end{gather}
is zero. 
 But \eqref{fibre-Frob-Delta} may also be considered as  the image of $\wh b$ by a   specialization homomorphism for the sheaf 
 $\mc H _{N, I, W}^{\leq\mu_{1 }}$ associated to  a specialization arrow $(\prod_{i\in I}\Frob_{\{i\}}^{dn_{i}})(\Delta(\ov \eta))\to \Delta(\ov v)$. Therefore  \eqref{fibre-Frob-Delta} is zero because 
 \begin{itemize}
 \item 
  $\prod_{i\in I}\Frob_{\{i\}}^{dn_{i}}
(\Delta(\ov\eta))$ belongs to  $ \Omega_{1}$ by hypothesis 
\item for every geometric point $\ov x$ of 
$\Omega_{1}$
and every  specialization arrow $\on{\mf{sp}}_{\ov x}: \ov x\to \Delta(\ov v)$,     $\on{\mf{sp}}_{\ov x}^{*}(\wh  b )$ vanishes  in 
$\restr{\mc H _{N, I, W}^{\leq\mu_{1 }}}{\ov x}$. 
\end{itemize}
 This  last assertion comes from the fact that    $\mc H _{N, I, W}^{\leq\mu_{1}}$ is smooth over $\Omega _{1}$ 
and that   the image of $\wh  b$ by every specialization homomorphism to  
$\restr{\mc H _{N, I, W}^{\leq\mu_{1}}}{\ov{\eta^{I}}}
$ is zero (since it is the case  of $\on{\mf{sp}}^{*}(\on{\mf{sp}}_{v}^{*}(\wh  b))$ and  
$\pi_{1}(\eta^{I},\ov{\eta^{I}})$ acts transitively on the specialization arrows from  $\ov{\eta^{I}}$ to  $\Delta(\ov v)$).   \cqfd

\noindent{\bf End of the proof of   \propref{injectivite-sp}. }
 Since $\prod _{i\in I}\Frob_{\{i\}}$ is the  total Frobenius, 
 $\prod_{i\in I}F_{\{i\}}^{dn }$ acts bijectively on 
 $\varinjlim _{\mu}\restr{\mc H _{N, I, W}^{\leq\mu}}{\Delta(\ov{\eta})} $ and sends $a_{(n_{i} )_{i\in I}}$ to  $a_{(n_{i}+n)_{i\in I}}$. From this and from  a) of \lemref{lem-a-b} we deduce easily that to  prove that  $a=a_{(0)_{i\in I}}$ is zero
   (and even that the whole sequence $a_{(n_{i})_{i\in I}}$ vanishes) it suffices to find  $(n_{i})_{i\in I}\in \N^{I}$ such that 
 \begin{gather*} 
 a_{(n_{i}+\alpha_{i})_{i\in I}}=0 \text{ for every } (\alpha_{i})_{i\in I}\in \prod _{i\in I} \{0,..., \dim W_{i}-1\}.
 \end{gather*}
This is possible by b) of \lemref{lem-a-b}, because the density of  the open subscheme $\Omega_{1}$ implies that we can find  $(n_{i})_{i\in I}\in \N^{I}$ such that 
\begin{gather*} 
 \prod_{i\in I}\Frob_{\{i\}}^{d(n_{i}+\alpha_{i})}
(\Delta(\ov\eta))\in \Omega_{1} \text{ for every } (\alpha_{i})_{i\in I}\in \prod _{i\in I} \{0,..., \dim W_{i}-1\}.
 \end{gather*} 
 This ends  the proof of  \propref{injectivite-sp}. \cqfd
 
Propositions  \ref{surjectivite-Hecke-fini-intro} and \ref{injectivite-sp}   imply the following corollary. 
 
 \begin{cor}\label{bijectivite-Hecke-fini}
The specialization homomorphism 
 \begin{gather}\label{sp*-sans-Hf2}\on{\mf{sp}}^{*}: \Big( \varinjlim _{\mu}\restr{\mc H _{N, I, W}^{\leq\mu}}{\Delta(\ov{\eta})} \Big)^{\mr{Hf}}\to 
 \Big( \varinjlim _{\mu}\restr{\mc H _{N, I, W}^{\leq\mu}}{\ov{\eta^{I}}}\Big)^{\mr{Hf}}\end{gather}
 is a bijection.     \end{cor}
\dem
The injectivity comes from  \propref{injectivite-sp}. Here is  the proof of the surjectivity. Let $c\in  \Big( \varinjlim _{\mu}\restr{\mc H _{N, I, W}^{\leq\mu}}{\ov{\eta^{I}}}\Big)^{\mr{Hf}}$. By proposition \ref{surjectivite-Hecke-fini-intro} we can find  $a\in  \varinjlim _{\mu}\restr{\mc H _{N, I, W}^{\leq\mu}}{\Delta(\ov{\eta})}$ such that $\on{\mf{sp}}^{*}(a)=c$. The injectivity of $\on{\mf{sp}}^{*}$ that we proved in  \propref{injectivite-sp} implies that  $a$ is Hecke-finite. \cqfd

    \section{Compatibility with the Satake isomorphism at unramified  places} 
    \label{subsection-intro-decomp}
    The goal of this section is to prove the   following  lemma, as well as    \propref{S-non-ram-concl-intro}  which was admitted at the end  of section \ref{intro-idee-heurist}.

  The following lemma  shows assertion  (v)  of   \propref{prop-SIf-i-ii-iii}, 
  namely that the 
    Hecke operators at unramified places are particular cases of excursion operators. 
    
    Let 
      $v$ be a place in  $X\sm N$.   We fix an embedding 
       $\ov F\subset \ov F_{v}$. 
   As previously  $\mbf 1\xrightarrow{\delta_{V}} V\otimes V^{*}$  and   $ V\otimes V^{*}\xrightarrow{\on{ev}_{V}} \mbf 1$ are the natural morphisms.

     \begin{lem} \label{S-non-ram-intro}     For every $d\in \N$ and every   $\gamma\in \on{Gal}(\ov {F_{v}}/F_{v})\subset  \on{Gal}(\ov F/F)$ such that $\deg(\gamma)=d$, $
   S_{\{1,2\},V \boxtimes V^{*},\delta_{V},\on{ev}_{V},(\gamma,1)}$  depends only on    $d$, 
  and if    $d=1$ it is equal to  $T(h_{V,v})$. 
     \end{lem}
       \noindent{\bf Proof.}
 We fix a  geometric point  $\ov v$ over   $v$ and a  specialization arrow 
    $\on{\mf{sp}}_{v}:\ov \eta\to \ov v$, associated to the embedding   $\ov F\subset \ov {F_{v}}$ chosen above. We still denote by      $\on{\mf{sp}}_{v}$ the  specialization arrow    
    $\Delta(\ov \eta)\to \Delta(\ov v)$ equal to  its image   by $\Delta$. In order to reduce the size of the following diagram  we set  $I=\{1,2\}$ and $W=V\boxtimes V^{*}$. 
    The diagram  
       $$  \xymatrix{  C_{c}^{\mr{cusp}}(G(F)\backslash G(\mb A)/K_N \Xi,E)\ar[d]_-{\restr{\mc C_{\delta_{V}}^{\sharp}}{\ov v}} \ar[dr]^{ \mc C_{\delta_{V}}^{\sharp}} & &
       \\ \Big( \varinjlim _{\mu} \restr{\mc H _{N, I, W}^{\leq\mu}}{\Delta(\ov v)}\Big)^{\mr{Hf}} \ar[r]^-{\mf{sp}_{v}^{*}} \ar[d]^{F_{\{1\}}^{\deg(v)d}} 
       & 
       \Big( \varinjlim _{\mu} \restr{\mc H _{N, I, W}^{\leq\mu}}{\Delta(\ov\eta)}\Big)^{\mr{Hf}} \ar@{=}[r]  
       &  H_{I,W} \ar[d]^{(\gamma,1)}
       \\  \Big( \varinjlim _{\mu} \restr{\mc H _{N, I, W}^{\leq\mu}}{\Delta(\ov v)}\Big)^{\mr{Hf}} \ar[r]^-{\mf{sp}_{v}^{*}}\ar[d]_-{\restr{\mc C_{\on{ev}_{V}}^{\flat}}{\ov v}}   & 
        \Big( \varinjlim _{\mu} \restr{\mc H _{N, I, W}^{\leq\mu}}{\Delta(\ov\eta)}\Big)^{\mr{Hf}}   \ar@{=}[r]  \ar[dl]^{ \mc C_{\on{ev}_{V}}^{\flat}} & 
     H_{I,W}
       \\ C_{c}^{\mr{cusp}}(G(F)\backslash G(\mb A)/K_N \Xi,E)  & &
       } $$
         is commutative (the commutativity of the big rectangle is proven in lemma 10.4 of \cite{coh}). 
   But   $S_{\{1,2\},V \boxtimes V^{*},\delta_{V},\on{ev}_{V},(\gamma,1)}$  is equal  by definition to  the composition by the  rightmost path.    Thus     it  is equal to the  composition given by the left column. Consequently   it depends only on    $d$. 
   When $d=1$ the composition given by the left column is equal  by definition to  $S_{V,v}$, and thus  to $T(h_{V,v})$  by    \propref{prop-coal-frob-cas-part-intro}. \cqfd

       \begin{rem} We computed   the composition given by   the left column only for $d=1$ because  we can prove that for other values  of $d$ it  is equal to  a combination of $S_{W,v}$ with $W$ irreducible representation of $\wh G$, and therefore it does not bring anything new. 
       \end{rem}
       
   The  following proposition  claims that     the decomposition \eqref{intro1-dec-canonical} is compatible with 
the Satake isomorphism  at the places of $X\sm N$.           

\begin{prop}\label{S-non-ram-concl-intro} 
 Let $\sigma$ be a parameter occuring 
 in \eqref{intro1-dec-canonical}  and  
let $v$ be a place    of $X\sm N$. Then  
 for every irreducible representation $V$ of $\wh G$, 
       $T(h_{V,v})$ acts on   $\mf H_{\sigma}$
       by multiplication by the scalar  $\chi_{V}(\sigma(\Frob_{v}))$, where $\chi_{V}$ is the character of $V$ and $\Frob_{v}\in \pi_{1}(X\sm N, \ov\eta)$ is a   Frobenius element at $v$. 
 \end{prop}
 \dem 
We adopt again  the notations of \lemref{S-non-ram-intro}. 
 Since 
    $  \s{\on{ev}_{V}, (\sigma(\gamma), 1) . \delta_{V} }=\chi_{V}(\sigma(\gamma))$ this lemma implies that 
        for every    irreducible    representation $V$ of $\wh G$, 
   $\mf H_{\sigma}$ is included  in the generalized eigenspace  of $T(h_{V,v})$ for the eigenvalue 
             $\chi_{V}(\sigma(\Frob_{v}))$. But we know that  the Hecke operators at unramified places are diagonalizable 
      (they are normal operators on  the hermitian space of   cuspidal automorphic   forms  with coefficients in $\C$). Therefore  $T(h_{V,v})$ acts on $\mf H_{\sigma}$  by multiplication by the scalar  
      $\chi_{V}(\sigma(\Frob_{v}))$. 
\cqfd

 This ends   the proof  of \thmref{intro-thm-ppal}.

      \section{Link  with the geometric Langlands program}
     \label{subsection-link-langl-geom}
  Il is obvious that    coalescence   and    permutation of legs 
are linked to factorization   structures   introduced by  
  Beilinson and Drinfeld
 \cite{chiral} and indeed our article uses in an essential way    the fusion product  on the affine grassmannian of Beilinson-Drinfeld in the geometric Satake  equivalence       \cite{mv,ga-de-jong}.  
Moreover the idea of spectral decomposition  is familiar  in the geometric Langlands program, see \cite{beilinson-heisenberg} and especially 
 corollary   4.5.5 of 
 \cite{dennis-laumon}  which claims (in the setting of the geometric Langlands program for   $D$-modules where the curve $X$ is defined over an algebraically closed field  of characteristic $0$) that  the DG-category of $D$-modules  on 
 $\Bun_{G}$ ``lies over'' the stack of 
 $\wh G$-local systems. Curiously we remark that we do not know how to  formulate an analoguous statement  with the  $\ell$-adic sheaves when $X$ is over $\Fq$ (the vanishing conjecture proved by Gaitsgory in \cite{ga-vanishing} appears as the top of the iceberg), and however our  article may be  considered as  a ``classical'' or rather  ``arithmetical'' version   of such a  statement. 
  
 In fact the link is much more direct that a simple analogy: 
we shall  see   that the conjectures of the $\ell$-adic geometric Langlands program allow to understand the excursion operators and provide  
a very enlightening explanation of our approach  thanks to a construction of   Braverman and Varshavsky \cite{brav-var} which generalizes the fact that a sheaf on  $\Bun_{G}$ gives rise by traces of Frobenius to a function on   $\Bun_{G}(\Fq)$.    We take here  $N$ empty, i.e.   $K_{N}=G(\mathbb O)$ but we could consider any     level $N$ (and even non-split  reductive groups). 

   The conjectures of  the geometric Langlands   program involve the following   Hecke functors :  for every representation $W$ of $(\wh G)^{I}$  the Hecke functor  
     $$\phi_{I,W}:D^{b}_{c}(\Bun_{G},\Qlbar)\to  D^{b}_{c}(\Bun_{G}\times X^{I}, \Qlbar)$$      
    is given   by 
     $$\phi_{I,W}(\mc F)=q_{1,!}\big(q_{0}^{*}(\mc F)\otimes \mc F_{I,W}\big)$$
 where $\Bun_{G}\xleftarrow{q_{0}}\Hecke_{I,W}^{(I)}\xrightarrow{q_{1}}\Bun_{G}\times X^{I}$ is the Hecke   correspondence   and  
 \begin{itemize}
 \item when   $W$ is irreducible, $\mc F_{I,W}$ is equal,  up to a shift,  to the  IC-sheaf  of $\Hecke_{I,W}^{(I)}$   
 \item in general  it is defined, functorially in $W$,   as  the inverse image of 
 $\mc S_{I,W}^{(I)}$ by the smooth natural morphism  
 $\Hecke_{I,W}^{(I)}\to \mr{Gr}_{I,W }^{(I)}/G_{\sum n_{i}x_{i}}$ (where the $n_{i}$ are big enough).  
 \end{itemize}
Let $\mc E$  be a $\wh G$-local system over $X$. Then  $\mc F\in D^{b}_{c}(\Bun_{G},\Qlbar)$ is said to be an eigensheaf  for  $\mc E$  if we have,  for every finite  set $I$ and every  representation $W$ of 
  $(\wh G)^{I}$,  an isomorphism 
  $\phi_{I,W}(\mc F)\isom \mc F\boxtimes W_{\mc E}$,   functorial in   $W$, and compatible with exterior products and with   fusion (i.e.   with the inverse image by the   diagonal morphism $X^{J}\to X^{I}$ associated to every map $I\to J$). The  conjectures of the geometric Langlands    program 
  imply the existence of  an eigensheaf    $\mc F$  for  $\mc E$   (which  satisfies an additional  Whittaker normalization condition which prevents it in particular to be zero). In the geometric Langlands    program $X$ and $\Bun_{G}$ are usually  defined  over an algebraically closed  field but here we work  over $\Fq$. 
 
Let $\mc F$ be an eigensheaf for  $\mc E$. We denote by   $f\in C(\Bun_{G}(\Fq),\Qlbar)$ the function associated to  $\mc F$, i.e.    for $x\in \Bun_{G}(\Fq)$, $f(x)=\on{Tr}(\Frob_{x}, \restr{\mc F}{x})$. 
Let $\Xi\subset Z(F)\backslash Z(\mathbb A)$ be a lattice. We assume that  $\mc F$ is $\Xi$-equivariant, so that  $f\in  C(\Bun_{G}(\Fq)/\Xi,\Qlbar)$ (decreasing  $\Xi$ if necessary this is implied  by an adequate  condition   on $\mc E$, in fact on its  image by the morphism from $\wh G$ to its abelianization). 
It is well-kwown that  $f$ is an eigenvector w.r.t.  Hecke operators: for every place $v$ and every   irreducible  representation $V$ of $\wh G$, 
$T(h_{V,v})(f)=\on{Tr}(\Frob_{v},\restr{V_{\mc E}}{v})f$, where 
$\Frob_{v}$ is  a Frobenius element at $v$. 

The following proposition  (which relies on a result that has not yet been written) formulates the compatibility between  the geometric Langlands program and the decomposition \eqref{intro1-dec-canonical}. 

\begin{prop}\label{prop-compat-langl-geom}
Let  $\mc F$ be a $\Xi$-equivariant eigensheaf for $\mc E$   such that   the function $f$ associated  to $\mc F$ is   cuspidal. Then  $f$ belongs to  $\mc H_{\sigma}$  where 
 $\sigma:\pi_{1}(X,\ov\eta)\to 
\wh G(\Qlbar)$ is the Galois   representation   corresponding to the local system  $\mc E$. 
\end{prop}
\dem 
In  \cite{brav-var}, Braverman and Varshavsky use a very general trace  morphism, and the fact that  
$\Cht_{I,W}^{(I)}$ is the intersection of the  Hecke correspondence  with the  graph of the  Frobenius endomorphism   of $\Bun_{G}$, to  construct  a  morphism of sheaves over $X^{I}$
\begin{gather}\label{trace-brav-var-piIW}\pi^{\mc F,\mc E}_{I,W}: \varinjlim_{\mu}\mc H_{N,I,W}^{\leq\mu}\to W_{\mc E}.\end{gather}
These morphisms are functorial   in  $W$, and compatible  with the  coalescence of legs and with the action of the  partial Frobenius morphisms (this last point has not yet been written). 
Moreover   $\pi^{\mc F,\mc E}_{\emptyset,\mbf 1}:C_{c}(\Bun_{G}(\Fq)/\Xi,\Qlbar)\to \Qlbar$ is nothing but  $h\mapsto \int_{\Bun_{G}(\Fq)/\Xi}     fh$. 
The properties of these   morphisms $\pi^{\mc F,\mc E}_{I,W}$ imply that for every $I,W,x,\xi,(\gamma_{i})_{i\in I}$, 
we have 
$$S_{I,W,x,\xi,(\gamma_{i})_{i\in I}}(f)=
\s{\xi, (\sigma(\gamma_{i}))_{i\in I}.x}f.$$ 
This finishes  the proof of   \propref{prop-compat-langl-geom}. \cqfd

\section{Relation with previous works}\label{intro-previous-works}

      The methods used in this work are completely different from the methods based on the trace formulas  which were developped notably   by Drinfeld 
   \cite{drinfeld78,Dr1,drinfeld-proof-peterson,drinfeld-compact}, Laumon, Rapoport and Stuhler \cite{laumon-rapoport-stuhler}, 
  Laumon \cite{laumon-drinfeld-modular,laumon-cetraro},    Laurent  Lafforgue \cite{laurent-asterisque,laurent-jams,laurent-inventiones,laurent-tata}, 
 Ngô Bao Châu \cite{ngo-jacquet-ye-ulm,ngo-modif-sym},  Eike Lau \cite{eike-lau,eike-lau-duke}, Ngo Dac Tuan \cite{ngo-dac-ast,ngo-dac-09,ngo-dac-11},  Ngô Bao Châu and  Ngo Dac Tuan  \cite{ngo-ngo-elliptique}, Kazhdan and  Varshavsky~\cite{kvar,var-SANT} and 
 Badulescu and Roche \cite{badulescu}. 
  
      Nevertheless  the action on the  cohomology of the permutation groups of the legs of the chtoucas occurs already in the works  of  Ngô Bao Châu,  Ngo Dac Tuan and Eike Lau that we have just quoted. These actions of  the permutation groups also play an essential role   in the geometric Langlands program, and notably in the proof of the vanishing conjecture by Gaitsgory  \cite{ga-vanishing}. 
Moreover 
 the coalescence of legs  appears  in the thesis of Eike Lau \cite{eike-lau} and it is also used  in the preprint   \cite{brav-var} of Braverman and Varshavsky  
 (in order to prove the non-vanishing of the  morphisms \eqref{trace-brav-var-piIW}).  
 The article \cite{var} of Varshavsky about   the stacks  
 of $G$-chtoucas and the very enlightening preprint   \cite{brav-var} of Braverman and  Varshavsky,   were essential for us.  Lastly we repeat the link, already mentionned in the previous section,  with the   corollary  4.5.5 of 
 \cite{dennis-laumon} (which by the way generalizes  the vanishing conjecture).

 In the rest of this  section 
 (which does not bring any new result and can be skipped by the reader) 
 we explain a few additional  arguments which enable us to give a new  proof of the inductive step of 
  \cite{laurent-inventiones} as a consequence  of \thmref{intro-thm-ppal}. 
  We take  care to avoid any circularity and do not use 
    the results which are now well-known but are  consequences  of  \cite{laurent-inventiones}.  Although we have in mind the case of $GL_{r}$, it is more natural to state the following lemma for arbitrary $G$.

                \begin{lem}\label{lem-sigma-dans-coho} 
 Let $\sigma$ be a parameter 
 occuring in  decomposition 
      \eqref{intro1-dec-canonical} (i.e. such that $\mf H_{\sigma}\neq 0$). 
     Let $V$ be an  irreducible representation of $\wh G$ and 
     $V_{\sigma}=\oplus_{\tau} \tau \otimes \mf V_{\tau}$ be the  decomposition of the semisimple   representation  $V_{\sigma}$ indexed by the 
     isomorphism classes   of    irreducible representations $\tau $ of $\pi_{1}(\eta, \ov\eta)$. Then if  $\mf V_{\tau}\neq 0$, $\tau\boxtimes \tau ^{*}$ occurs  as a subquotient of the   representation 
     \begin{gather}\label{rep-12-VV*}  H_{\{1,2\}, V\boxtimes V^{*}}=    \Big(\varinjlim _{\mu}\restr{\mc H _{N, \{1,2\}, V\boxtimes V^{*}}^{\leq\mu}}{\ov{\eta^{\{1,2\}}}}\Big)^{\mr{Hf}}\end{gather}  of 
  $(\pi_{1}(\eta,\ov\eta))^{2}$. 
        Moreover $\tau$ is $\iota$-pure for every  isomorphism $\iota:\Qlbar \isom \C$.      \end{lem}

\begin{rem} Of course the previous  assertion is not a new  result because  theorem VII.6 of   \cite{laurent-inventiones} implies that every     irreducible  representation (defined over a finite  extension of $\Ql$ and continuous) of $\pi_{1}(X\sm N, \ov\eta)$ is $\iota$-pure for every $\iota$. 
\end{rem}

   \noindent{\bf Proof. } Increasing  $E$ if necessary, we assume that    $\sigma$ and $\mf H_{\sigma}$ are defined  over $E$. Let $h\neq 0$ be in the subspace of $\mf H_{\sigma}$ over which  $\mc B$ acts by the character $\nu$ associated to $\sigma$ by \eqref{relation-fonda} (we know this subspace is nonzero although we do not know whether  $\mc B$ is reduced). 
   Let $\check h\in   C_{c}^{\rm{cusp}}(G(F)\backslash G(\mb A)/K_{N}\Xi,E)$ be such that \begin{gather}\label{int-check-h-h}\int_{G(F)\backslash G(\mb A)/K_{N}\Xi} \check h \,  h=1.\end{gather} 
  We denote by  
  \begin{itemize}
  \item $f$ the element of \eqref{rep-12-VV*} equal to the image of $h$ by the composition 
 $$C_{c}^{\rm{cusp}}(G(F)\backslash G(\mb A)/K_{N}\Xi,E)=
  H_{\{0\},\mbf  1}\xrightarrow{\mc H(\delta_{V})}
 H_{\{0\},V\otimes V^{*}}\isor{\chi_{\zeta_{\{1,2\}}}^{-1}} 
 H_{\{1,2\}, V\boxtimes V^{*}},$$
  \item $\check f$ the linear form    over  \eqref{rep-12-VV*}
  equal to  the composition of  
 $$H_{\{1,2\}, V\boxtimes V^{*}} \isor{\chi_{\zeta_{\{1,2\}}}} H_{\{0\},V\otimes V^{*}}  \xrightarrow{\mc H(\on{ev}_{V})} 
  H_{\{0\},\mbf  1} =C_{c}^{\rm{cusp}}(G(F)\backslash G(\mb A)/K_{N}\Xi,E)$$
     and of the linear form 
     $$C_{c}^{\rm{cusp}}(G(F)\backslash G(\mb A)/K_{N}\Xi,E)\to E, \ \ g\mapsto \int_{G(F)\backslash G(\mb A)/K_{N}\Xi} \check h g.$$
  \end{itemize}
        Then $f$  and $\check f$ 
      are invariant under  the diagonal action  of $\pi_{1}(\eta,\ov\eta)$. 
    For every   $(\gamma, \gamma')\in (\pi_{1}(\eta,\ov\eta))^{2}$ we have  
        \begin{gather*}\s{\check f , (\gamma, \gamma') \cdot f}=
       \int_{G(F)\backslash G(\mb A)/K_{N}\Xi} \check h S_{\{1,2\}, V\boxtimes V^{*},\delta_{V},\on{ev}_{V},(\gamma,\gamma')}(h)\\
    =\nu(S_{\{1,2\}, V\boxtimes V^{*},\delta_{V},\on{ev}_{V},(\gamma,\gamma')})       =
        \chi_{V}(\sigma(\gamma\gamma'^{-1}))=  \chi_{V_{\sigma}}(\gamma\gamma'^{-1}),\end{gather*}  where 
   \begin{itemize}
   \item  the first equality  comes from the definition of the  excursion operators given in \eqref{excursion-def-intro}, 
      \item the second equality comes from    the hypothesis that $h$ is an eigenvector for $\mc B$ w.r.t. the  character $\nu$, and from  \eqref{int-check-h-h}, 
      \item the third  equality comes from the fact that  $\nu$ is associated to $\sigma$ by \eqref{relation-fonda}. 
        \end{itemize}
       The   quotient of the   representation of $(\pi_{1}(\eta,\ov\eta))^{2}$ generated by  $f$ by the biggest  subrepresentation on which  $\check f$ vanishes is then  isomorphic to the subrepresentation  generated by 
       $\chi_{V_{\sigma}}$  in  
     $C( \pi_{1}(\eta,\ov\eta), E)$ equipped with the action  by left and right translations by $(\pi_{1}(\eta,\ov\eta))^{2}$. By \cite{Bki-A8} chapter 20.5 theorem 1, this   representation is isomorphic to 
   $\oplus_{\tau, \mf V_{\tau}\neq 0} \tau\boxtimes \tau ^{*}$.   
    We have shown that  if   $\mf V_{\tau}\neq 0$, $\tau\boxtimes \tau ^{*}$ is a  quotient of a subrepresentation of  \eqref{rep-12-VV*}.      
    From this we deduce that   $\tau$ is $\iota$-pure.  Since $\tau\boxtimes \tau^{*}$ is a subquotient of \eqref{rep-12-VV*}, 
         it results from  Weil II \cite{weil2} that  $\tau\boxtimes \tau ^{*}$ is 
  $\iota$-pure of weight $\leq 0$ as a representation of $\pi_{1}((X\sm N)^{2}, \Delta(\ov\eta))$. 
 Therefore for almost all  place $v$ the eigenvalues of $\tau(\Frob_{v})$ have equal    $\iota$-weights, and they are determined  by the $\iota$-weight of $\det(\tau)$. 
     \cqfd
       
   From now on we take   $G=GL_{r}$.  
                              We recall that in \cite{laurent-inventiones} the Langlands correspondence  is obtained by induction on $r$,  with the help  of   the ``induction principle'' of Deligne,  which combines  
  \begin{itemize}
  \item
    the functional equations of the $L$-functions due to Grothendieck \cite{sga5}, 
    \item the product formula   of Laumon \cite{laumon-produit}, 
   \item  the multiplicity  one theorems   \cite{piat71,shalika}  
and  the  converse theorems   of Hecke, Weil,  Piatetski-Shapiro and Cogdell \cite{inverse-thm}.
\end{itemize}  
   The induction is explained  in   section 6.1 and appendix  B of 
    \cite{laurent-inventiones} and  the induction step is 
    \begin{gather}\label{hyp-laurent}
  \text{the hypothesis of   proposition VI.11 (ii) of   \cite{laurent-inventiones}}
  \end{gather}
  namely that to every    cuspidal automorphic representation  $\pi$ for $GL_{r}$ of level $N$ we can associate   $\sigma: \pi_{1}(X\sm N, \ov\eta)\to GL_{r}(\Qlbar) $ defined  over a finite  extension of  $\Ql$, continuous, 
    pure of weight $0$ and corresponding to $\pi$ in the sense  of Satake at all the  places of $X\sm N$. Our theorem
    \ref{intro-thm-ppal} provides a new  proof of  \eqref{hyp-laurent}, thanks to the   following lemma. 
    
    \begin{lem} We take $G=GL_{r}$. 
    We assume that  the Langlands correspondence  is known for  $GL_{r'}$ for every $r'<r$. 
         Then every  $\sigma$ 
  occuring  in the decomposition \eqref{intro1-dec-canonical} is irreducible and pure  of weight  $0$. 
    \end{lem}
    \begin{rem} This lemma does not bring any new result because {\it a posteriori} it results from   \cite{laurent-inventiones}. \end{rem}
 
 \noindent{\bf Proof (extracted from  \cite{laurent-inventiones}). }   
 Let $(H_{\pi},\pi)$ be a cuspidal automorphic    representation such that  $(H_{\pi})^{K_{N}}$ is non zero and appears in $\mf H_{\sigma}$. 
      We denote by   
   \begin{gather}\label{dec-sigma}\sigma=\oplus_{\tau} \tau \otimes \mf V_{\tau}\end{gather} 
   the decomposition of the  semisimple representation $\sigma$ indexed by   equivalence classes of irreducible representations  $\tau $ of $\pi_{1}(X\sm N, \ov\eta)$. 
    We assume by contradiction that   $\sigma$ is not irreducible. 
   Any  representation   $\tau$ such that  $\mf V_{\tau}\neq 0$ is therefore of rank   $r_{\tau}<r$ and,  by 
    hypothesis, there exists  a   cuspidal automorphic representation    $\pi_{\tau}$ for  $GL_{r_{\tau}}$  associated to   $\tau$ by the Langlands correspondence  for $GL_{r_{\tau}}$.  We choose a finite  set $S$ of places  outside of which  the representations $\pi,\sigma$, and $\tau $, $\pi_{\tau}$ such that  $\mf V_{\tau}\neq 0$ (whose number is finite)  are unramified and correspond to each other  by the Satake isomorphism. 
 By the  method of Rankin-Selberg for $GL_{r}\times GL_{r_{\tau}}$ (due to  Jacquet, Piatetski-Shapiro, Shalika \cite{jacquet-shalika-euler,rankin-selberg}), 
$L(\check{\pi} \times \pi_{\tau},Z)$ is a polynom in $Z$  
hence {\it  a fortiori}  
 $L_{S}(\check{\pi} \times \pi_{\tau},Z)$ is a polynom in  $Z$ 
 (since the local factors may have  poles but never zeros). 
  Thus 
$$L_{S}(\check{\pi}  \times \pi,Z)=L_{S}(\check{\pi} \times \sigma,Z)=\prod_{\tau }L_{S}(\check{\pi} \times \tau,Z) ^{\dim \mf V_{\tau}}=\prod_{\tau }L_{S}(\check{\pi} \times \pi_{\tau},Z)^{\dim \mf V_{\tau}}$$  has no pole. 
However by  theorem B.10 of \cite{laurent-inventiones} (due to Jacquet, Shahidi, Shalika), 
 $L_{S}(\check{\pi}  \times \pi,Z)$ has a pole at  $Z=q^{-1}$, and we get  a contradiction. 
    We have proven that  $\sigma$ is irreducible. For every $\iota$ we know by   \lemref{lem-sigma-dans-coho} that $\sigma$ is $\iota$-pure and  the knowledge of its determinant (by class field theory) implies then that  the $\iota$-weight is zero. Thus $\sigma$ is pure of weight $0$. \cqfd
  
  \begin{rem} We mention  for the reader  that  the already very important consequences  of the Langlands correspondence  for $GL_{r}$, explained in the chapter VII of \cite{laurent-inventiones},  were widened in recent works  of Deligne and Drinfeld 
  \cite{deligne-finitude,esnault-deligne,drinfeld-del-conj}
 about  rationality 
  of the  Frobenius  traces  and about the independence of $\ell$ for the $\ell$-adic  sheaves  over    smooth varieties over $\Fq$. 
\end{rem}

\end{document}